\documentclass[11pt,letterpaper]{article}
\usepackage[utf8]{inputenc}

\usepackage[a4paper, total={6.2in, 6.2in}, top={1in}, bottom={1.2in}
]{geometry}

\usepackage{amstext,amsmath, amssymb, amsthm, mathrsfs, mathtools}
\usepackage{authblk}

\usepackage{thmtools}

\usepackage{enumerate}
\usepackage[dvipsnames]{xcolor}
\usepackage{comment}
\usepackage[colorlinks=true,
    linkcolor=black,
    citecolor=RoyalBlue,
    urlcolor=cyan, hypertexnames=false]{hyperref}
\usepackage[capitalise]{cleveref}
\usepackage{etoolbox}
\usepackage{stmaryrd}

\usepackage{mathrsfs}
\usepackage{dutchcal}

\usepackage{setspace}
\allowdisplaybreaks

\usepackage[titles]{tocloft}

\usepackage{relsize}

\usepackage[backend=bibtex, style=alphabetic]{biblatex}
\newtheorem{theorem}{Theorem}[section]
\newtheorem{lemma}[theorem]{Lemma}

\newtheorem{proposition}[theorem]{Proposition}
\newtheorem{corollary}[theorem]{Corollary}
\newtheorem{definition}[theorem]{Definition}

\theoremstyle{definition} 

\numberwithin{equation}{section}

\DeclareMathOperator{\nR}{\mathbb{R}}

\DeclareMathOperator{\nS}{\mathbb{S}}

\DeclareMathOperator{\K}{\operatorname{K}}

\newcommand{\dd}{{\rm d}}

\DeclareMathOperator{\R}{R}

\DeclareMathOperator{\B}{B}
\DeclareMathOperator{\C}{C}

\DeclareMathOperator{\LL}{\mathfrak{L}}
\DeclareMathOperator{\p}{\mathrm{p}}
\DeclareMathOperator{\q}{\mathrm{q}}

\DeclareMathOperator{\olNabla}{\overline{\nabla}}
\newcommand{\DD}{\mathrm{D}^2}

\DeclareMathOperator{\bb}{\boldsymbol{b}}

\DeclareMathOperator{\xx}{\boldsymbol{x}}
\DeclareMathOperator{\yy}{\boldsymbol{y}}
\DeclareMathOperator{\zz}{\boldsymbol{z}}

\DeclareMathOperator{\vv}{\mathlarger{\mathlarger{\mathrm{v}}}}

\DeclareMathOperator{\nnu}{\boldsymbol{\mathlarger{\mathlarger{\nu}}}}

\newcommand{\G}{\mathfrak{g}}
\newcommand{\HH}{\operatorname{H}}

\newcommand{\DDD}{\mathrm{D}}
\newcommand{\ttt}{\mathcal{t}}

\newcommand{\TT}{\mathsf{T}}
\newcommand{\XX}{\mathsf{X}}

\DeclareMathOperator\arctanh{artanh}
\DeclareMathOperator\arccosh{arcosh}

\usepackage{lipsum}
\usepackage{appendix}

\usepackage{thmtools}

\title{\textbf{Entire boosted hypersurfaces of vanishing mean curvature in  
Schwarzschild spacetime}}

\makeatletter

\renewcommand{\AB@authnote}[1]{}
\renewcommand{\AB@affilnote}[1]{}

\renewcommand{\@author}{\AB@authlist}

\AtEndDocument{
    \par\bigskip
    \noindent\AB@affillist\par
}

\makeatother

\makeatletter
\renewcommand{\AB@affilsepx}{\protect\\[10pt]\protect\Affilfont}
\makeatother

\newcommand{\emailbreak}{\protect\\[6pt]}

\author{ \Large Albachiara Cogo}

\affil{\small 
University of T\"ubingen, Fachbereich Mathematik, Auf der Morgenstelle 10, 72076, T\"{u}bingen, Germany }
\affil{
\textit{Current affiliation:}  Centro di Ricerca Matematica Ennio de Giorgi, Scuola Normale Superiore di Pisa, Piazza dei Cavalieri 3, 56126, Pisa, Italy
\emailbreak
 \indent \texttt{albachiara.cogo@sns.it}}

\date{}

\begin{document}

\maketitle

\begin{abstract}
   \noindent We prove the existence of entire maximal hypersurfaces,
asymptotic to suitable coordinate-dependent boosts in the two asymptotically flat ends of maximally extended Schwarzschild spacetime.
As a corollary, this generates a foliation of maximal hypersurfaces in exterior Schwarzschild, which corresponds to a geometric choice of time function associated with any boost parameter. 
    As a crucial step in the proof, we introduce a novel rotationally symmetric time function in maximally extended Schwarzschild, which is
    Cauchy
    and 
    has properties of independent interest.
\end{abstract}

\setcounter{tocdepth}{2}

\vspace{0.4cm}

\renewcommand{\contentsname}{\large Contents}
\vspace{0.25cm}
\tableofcontents
\addtocontents{toc}{\protect\thispagestyle{empty}}


\section{Introduction}

A spacelike hypersurface is called \emph{maximal} if its mean curvature $\HH$ vanishes identically, that is, it is a solution to the Euler--Lagrange equation of the variational problem of area maximization over (weakly) spacelike hypersurfaces in a Lorentzian manifold. 

The existence of foliations by maximal hypersurfaces has been essential in \emph{General Relativity}, in particular in the study of its dynamic aspects, as well as in the related numerical computation. In fact, setting $\HH = 0$ drastically simplifies the \emph{Einstein constraint equations}, which describes admissible initial data to be evolved into solutions to the Einstein equations; this has been crucial in the \emph{conformal method} approach initiated in \cite{Lichnerowicz_1, York_OMurchadha}.
Moreover, such condition reduces the so-called \emph{Dominant Energy condition} for admissible initial data to the intrinsic geometric condition of non-negative scalar curvature $\R_g \geq 0$ and was assumed in 
the first proofs of the \emph{Positive Mass Conjecture} \cite{BrillDeser_PMT, ChoquetMarsden_PMT}. Finally, the preservation of $\HH = 0$ along the evolution, provides an elliptic PDE for the lapse function, which was essential for the development of numerical methods \cite{Wheeler, York_Numerics1} and in the proof of \emph{global nonlinear stability of Minkowski} \cite{Christodoulou}. 
Motivated by the above physically and numerically relevant problems, the study of the existence of \emph{maximal hypersurfaces} and the corresponding Dirichlet problem was pioneered by 
Choquet-Bruhat \cite{Choquet-Bruhat_Max} and Brill--Flaherty \cite{BrillFlaherty_MaximizingProperties} (see also \cite{Bancel, EardleySmarr, MT, Moncrief1981-rk, IseMon_Max3Torus}), but under the restrictive hypotheses on specific spacetimes  admitting compact slices and on the boundary data; 
preliminary results on asymptotically Minkowskian spacetimes in vacuum were also proposed 
in \cite{CantorFisherMarsden_MaxAsyFlat}. While the existence of weakly spacelike solutions to the variational problem is ensured by the early work of Avez \cite{Avez}, a regularity theory does not directly promote them to classical solutions of $\HH = 0$, in contrast for example to the case of graphical minimal surfaces. This is a consequence of the fact that in local coordinates, the equation is a \emph{second-order quasilinear elliptic} partial differential equation, which \emph{degenerates} as the hypersurface becomes null. This analytical difficulty  is already evident in the flat case of Minkowski $\left(\nR^{1, n}, \, \eta \coloneqq -dt^2 + \sum_{i = 1}^3 dx_i^2\right)$, where the variational problem is established on Liptschitz functions 
with $| \DDD u | \leq 1$. 
On the other hand, $\HH = 0$ reads
\begin{equation}\label{eq: maximal Minkowski}
    0 = a^{ij}(\DDD u) \,\DD_{ij} u, \qquad a^{ij}(\DDD u) = \frac{1}{\sqrt{1-|\DDD u|^2}}\bigg{(} \delta^{ij} + \frac{ \DDD^iu \;\DDD^ju}{1-|\DDD u|^2} \bigg{)}\,, 
\end{equation}
which is elliptic only if $|\DDD u | < 1$, uniformly. 

The first systematic theory in full generality for maximal surfaces in Minkowski was developed  by Bartnik--Simon
\cite{BartSim}, who also made clear the link between classical solution to \eqref{eq: maximal Minkowski} and the solutions to the variational problem. The interest in variational solutions (e.g.~\cite{Eck, Hong_AreaMaximizing}) and classical solutions over exterior domains \cite{Hong_Mink_ExteriorDomains} in Minkowski is still an active field of study. The non-flat case was known to be more complicated, since pathological spacetimes could prevent the existence of maximal hypersurfaces, as in the work of Brill \cite{Brill} (see also the non-existence results for constant mean curvature surfaces \cite{Bartnik_NoCMC, Chrusc_Gluing_NoCMC, LingOhan_CMC}). Exhaustive results on the solvability of the Dirichlet
problem were proposed at the same time by \cite{Ger} and Bartnik \cite{Bart}; the latter, in addition, developed techniques to prove existence of maximal hypersurfaces in asymptotically Minkowskian spacetimes with prescribed asymptotics, under certain assumptions on the spacetime (we call them \emph{Bartnik interior condition} and \emph{time-dependent Bartnik conditions}) and on its time slices (which have to be suitably decaying mean curvature). The relation with the variational solutions was studied by Bartnik in \cite{Bartnik_VariationalMax}; analogous results for
asymptotically anti-de Sitter spacetimes were proven in \cite{Akutagawa}, while for non-zero constant mean curvature surfaces in \cite{Iriondo}.

Before any existence result, the equation \eqref{eq: maximal Minkowski} in the flat case was studied in the context of rigidity for entire solutions to elliptic partial differential equations. The \emph{Calabi--Bernstein theorem} proven by Calabi \cite{Calabi1968ExamplesOB} and Cheng--Yau \cite{ChenYau} asserts that, in any dimensions, an entire maximal hypersurface in Minkowski must be an affine plane. Beyond its striking mathematical interest, such theorem also has a deep physical impact: it established a one-to-one correspondence between 
entire maximal hypersurfaces and affine planes, which represent the time slices of an inertial observer. In turn, any inertial observer corresponds to an isometry of the Minkowski spacetimes, called \emph{boost}, given by the map $(t, \xx) \mapsto (\widetilde{t}, \widetilde{\xx})$ with
\begin{equation}\label{def: boostMinkowski}
    \widetilde{t} =\frac{t - \bb \cdot \xx }{\sqrt{1 - \bb^2}} \,, \qquad \widetilde{\xx} = \frac{\xx -  t\bb }{\sqrt{1 - \bb^2}} \,
\end{equation}
for a fixed parameter $\bb \in \nR^3, \,|\bb| < 1$. Generalizing the Calabi--Berstein theorem to asymptotically Minkowskian spacetimes (see \cite{Bart_OpenProb}) would provide a canonical geometric choice of time functions. While non-flat spacetimes do not display the symmetries of Minkowski, if they are asymptotically Minkowskian,
they admit an asymptotic notion of boosts; on the other hand, such notion will depend on the choice of asymptotically Minkowskian coordinates. For a classification of the candidates for such a Bernstein-type statement, one needs to prove the existence of maximal hypersurfaces prescribing boosted slices as asymptotics.
Such existence was still unproved in the simplest nontrivial asymptotically Minkowskian solution to the vacuum Einstein equations, that is the maximally extended Schwarzschild spacetime. In this work we provide a resolution to this problem and we show the existence of entire solutions, which, in particular,  pass through the bifurcation
sphere. Both the result by Bartnik \cite{Bart} and by Bartnik--Chruściel--O’Murchadha \cite{Chr1} do not directly apply to our case, although they are important for our strategy; it is, in particular, remarkable that the latter made progress in dealing with the construction of barriers in the \emph{asymptotically flat ends}.

Complementary to ours,  
new insights on non-zero constant mean curvature in Schwarzschild spacetime have also been provided by recent works of Tam and
Li--Shi--Tam \cite{Tam2023SpacelikeCS, LiShiTam}.

\subsection{Outline of the results}

In the following, $\left( \LL_m, \mathfrak{g}_m \right)$ denotes the \emph{maximally extended Schwarzschild spacetime} of mass-parameter $m > 0$, which is described by the \emph{Kruskal--Szekeres coordinates} $(\TT, \XX, \theta, \varphi)$. Its \emph{region} $\mathrm{I}$ and \emph{region} $\mathrm{I}'$ are isometric to two copies of \emph{exterior Schwarzschild} $\left( (\LL_m)_{\mathrm{ext}}, \mathfrak{g}_m \right)$, which is asymptotically Minkowskian; see the definition \eqref{def: KruskalExtension} and the properties in \cref{subsec: Schwarzschild}. The so-called \emph{asymptotically flat ends} are, in fact contained in the \emph{region} $\mathrm{I}$ and \emph{region} $\mathrm{I}'$.
We say that a spacelike hypersurface   
is \emph{entire} if it can be written as a global graph over a Cauchy hypersurface, given a Cauchy splitting of the spacetime.

Before stating the main Theorem and its Corollary, we need to introduce the notion of \emph{coordinate boost}.  
Given any choice of asymptotically Minkowskian coordinates in the \emph{regions} $\mathrm{I}$ and $\mathrm{I}'$, with respect to which the metric decays to the Minkowski one, there is an asymptotic notion of boosted hypersurfaces.
To give a definition, we employ the so-called \emph{isotropic coordinates} (see the definition \eqref{def: Schw_Isotropic} and the properties in \cref{subsec: Schwarzschild}) and we consider graphs of linear functions with slope $\bb \in [0,1)$, allowing a deviation of logarithmic order. 
The logarithmic deviation is coherent with
the fact that there is a priori no preferred choice of asymptotically Minkowskian coordinates 
and with the following observation: while in Minkowski the radial geodesics lie on $\{ t = \mathcal{r}\}$ for $\mathcal{r} >0$ radial coordinate, in Schwarzschild they lie on $\left\{ t = r + 2m \log \left( \frac{r}{2m} -1\right)\right\}$ in \emph{Schwarzschildean coordinates} $t \in \nR, r \in (2m, +\infty)$; see the definition 
\eqref{def: Schw_Schwarzschildean} and the properties in \cref{subsec: Schwarzschild}.

\vspace{0.1cm}

\begin{definition}[\textbf{Coordinate boost}]\label{def: coordinate_boost}
    A (symmetric) 
    \emph{coordinate boost} of parameter $\bb \in [0, 1)$ 
    is a smooth spacelike hypersurface $\mathcal{S}_{\bb} \hookrightarrow \left( \LL_m, \mathfrak{g}_m \right)$ such that, outside a spatial compact set, is contained in  \emph{region} $\mathrm{I}$ and \emph{region} $\mathrm{I}'$. Moreover, 
    in \emph{region} $\mathrm{I}$, 
    in isotropic coordinates  
    it reads 
    $$\mathcal{S}_{\bb} \cap \{ \vert \xx \vert > Km\} = \left\{ (t, \xx) \in (\LL_m)_{\mathrm{ext}} \, : \,  t = u_{\bb}(\xx) \right\}$$
    for some $K > \frac{1}{2}$ and a smooth
    function $
    u_{\bb} : \nR^3 \setminus B_{Km} \to \nR$, which satisfies
    \begin{equation}\label{eq: coordinate_boost}
         u_{\bb}(\xx) = \bb x_1 + O\left(\log(  \vert \xx \vert)\right)\, \qquad \text{ as } \vert \xx \vert \to + \infty.
    \end{equation}
    The same property has to hold 
    in \emph{region} $\mathrm{I}'$, after applying the isometry $(\XX, \TT) \mapsto (-\XX, -\TT)$.

    We denote the class of symmetric coordinate boosts of parameter $\bb \in [0, 1)$ by $\mathcal{B}_{\bb}$.
\end{definition}

\vspace{0.1cm}

\begin{theorem}[\textbf{Main Theorem}]\label{thm: Main intro}
    In the maximally extended Schwarzschild spacetime,
    for any $\bb \in [0, 1)$, there exists an entire  
    smooth maximal hypersurface 
    $$M_{\bb} \in \mathcal{B}_{\bb}\,, \qquad  \HH(M_{\bb}) = 0\,, \qquad \{ \XX = \TT  =0\} \subset M_{\bb} \,.
    $$
\end{theorem}

\vspace{0.1cm}

\begin{corollary}[\textbf{Geometric boosted time function}]\label{coroll: main}
    For any $\bb \in [0, 1)$, the exterior Schwarzchild spacetime $\left( (\LL_m)_{\mathrm{ext}},  \mathfrak{g}_m\right)$ admits a foliation by  
    maximal hypersurfaces
    in the class $\mathcal{B}_{\bb}$, generated by an entire maximal hypersurface as in \cref{thm: Main intro}.
    This provides a geometric choice of time function associated to the asymptotic boost of parameter $\bb \in [0, 1)$.
\end{corollary}

The proof of \cref{thm: Main intro} is based on the following broad strategy: we carefully construct a coordinate boost $ \mathcal{S}_{\bb} \in \mathcal{B}_{\bb}$ of parameter $\bb \in [0, 1)$, and prove the existence of a maximal hypersurface $M_{\bb}$, which is asymptotic to it in both \emph{regions} $\mathrm{I}$ and $\mathrm{I}'$. To this end, we consider a sequence $\{M_k\}_{k=1}^{\infty}$ of (compact) graphical solutions to a family of Dirichlet problems associated with the maximal surface equation $\HH = 0$ with boundary data $\{\Sigma_k\}_{k=1}^\infty$, corresponding to the boundaries of an exhaustion of  $ \mathcal{S}_{\bb} \in \mathcal{B}_{\bb}$. 
To set up the a priori height and gradient estimates, we need $\mathcal{S}_{\bb}$ to be zero time slice of a time function  satisfying 
the so-called \emph{time-dependent Bartnik  conditions} and \emph{Bartnik interior condition}; for the precise definitions, see \eqref{cond: Bartnik_RiemNorm} and \eqref{cond: Bartink_Interior} in \cref{sec: prelim}.
We then prove 
a priori height and gradient estimates with respect to such time function and we show that they are independent of the parameter $k$.
Since these estimates, in particular, control the ellipticity of the maximal surface operator, when combined with the De Giorgi--Nash theory and interior Schauder estimates, they imply that the sequence of solutions $\{M_k\}_{k=1}^{\infty}$ sub-converges in $C^{2, \alpha}$ on any compact set. Using a diagonal argument, one obtains an entire maximal hypersurface $M_{\bb}$. To ensure that $M_{\bb}$ is itself a coordinate boost in the sense of \cref{def: coordinate_boost}, we prove that it is asymptotic to the fixed coordinate boost $\mathcal{S}_{\bb} \in \mathcal{B}_{\bb}$ employing upper and lower barriers and showing that all $M_k$ lie between these;  such barriers are also crucial in proving the a priori height estimate itself.

As we have mentioned, a 
 major step of the proof consists in the construction of a suitable time function, as in the following Theorem.

\vspace{0.2cm}

\begin{theorem}[simplified version]\label{thm: new time t_b intro}
    For any $\bb \in [0, 1)$, the maximally extended Schwarzschild spacetime admits a time function $\ttt_{\bb} \in C^{\infty}(\mathfrak{L}_m)$ such that the following hold simultaneously.
    \begin{itemize}
        \item [$(i)$]  The \emph{time-dependent Bartnik  conditions} hold on  a region $ \{ r > \kappa_0m\}$ for any $\kappa_0 > 0$.
        \item [$(ii)$] Each time slice is a coordinate boost of parameter $\bb \in [0, 1)$ as in \cref{def: coordinate_boost}.
         Moreover, the spacelike hypersurface $\mathcal{S}_{\bb} = \{ \ttt_{\bb} = 0\} \in \mathcal{B}_{\bb}$ 
         admits upper and lower barriers,
         and it is invariant under the isometry $(\XX, \TT) \mapsto (-\XX, -\TT)$. 
        \item [$(iii)$] A relaxed version of \emph{Bartnik interior condition} 
        holds in $ \{ r > \kappa_0m\}$ for any $\kappa_0  > 0$. 
    \end{itemize}
\end{theorem}

\vspace{0.2cm}

\noindent
We remark that the fact that the properties hold in $\{ r > \kappa_0m\}$, where $r$ is the Schwarzschildean radial coordinate, is not restrictive because the spacetime singularity 
of Schwarzschild is
 \emph{crushing}. In the following, we briefly describe the role of the properties of the time function $\ttt_{\bb}$ as in \cref{thm: new time t_b intro}. The condition $(i)$, combined with the fact that the boundary data $\Sigma_k$ lie on the time slice $\mathcal{S}_{\bb} = \{ \ttt_{\bb} = 0\}$ (and with a control on the spacetime mean curvature of $\Sigma_k$), is essential to derive an a priori gradient estimate from an a priori height estimate; this is  a consequence of \cite[Theorem 3.1-$(iii)$]{Bart}. On the other hand, proving an a priori  height estimate heavily relies on the simultaneous validity of all the conditions in \cref{thm: new time t_b intro}. In fact, as in \cite[Theorem 5.3]{Bart}, the height estimate is obtained combining the barriers as in $(ii)$ and integral estimates. Such integral estimates make use of the fact that the spacetime is asymptotically Minkowskian and the \emph{Bartnik interior condition} as in $(iii)$, in combination with the pointwise gradient estimate in \cite[Theorem 3.1-$(iv)$]{Bart}. In turn, the latter holds under the \emph{time-dependent Bartnik conditions} as in $(i)$. We also remark that, while our height estimate will
follow the proof strategy of \cite[Theorem 5.3]{Bart}, there are major difficulties, as   we  work with the spatial isotropic coordinates  instead of finding explicit coordinates on the time slices of $\ttt_{\bb}$. The asymptotic decays of geometric quantities on the time slices are more delicate in our case: for example, the second fundamental form 
decays as $O\left( \frac{\log(\vert \xx \vert)}{\vert \xx \vert^2}\right)$,
rather than  $O\left(\frac{1}{\vert \xx \vert^2}\right)$.

The time function $\ttt_{\bb}$ of \cref{thm: new time t_b intro} is constructed by gluing on large annuli regions two time functions: the first is set up in the asymptotically flat ends so that $(ii)$ holds, while the other is defined on the rest of the spacetime and is tailored to satisfy $(i)$ and $(iii)$.

In the asymptotically flat ends, we consider coordinate boosts as candidates for the zero time slice. A crucial problem is given by the fact that the natural candidate $\{ t = \bb x_1\}$ in isotropic coordinates, for example, does not allow for a construction of barriers that adapts barriers from the Minkowski spacetime. On the other hand, the possibility of deviating from the linear graph by $O( \log(\vert \xx \vert))$ 
in \cref{def: coordinate_boost} enables to achieve this objective.
In particular, we correct the standard coordinate boost by a function that solves a linear elliptic equation and with at most logarithmic growth; its existence 
is based on \cite[Lemma 5, Theorem 2]{Meyers} and is analogous to \cite[Lemma 2.3]{Chr1}, although structurally different. The resulting time function on large radii will be $t-u_{\bb}(\xx)$, where $\{ t = u_{\bb}(\xx)\}$ is the coordinate boost we built.

In the remaining part of the spacetime $(\LL_m, \mathfrak{g}_m)$, we need to propose a novel time function. In fact, the \emph{time-dependent Bartnik conditions} fail in the interior of the maximally extended Schwarzschild spacetime when considering the time coordinate $\TT$ of Kruskal--Szekeres. 
We thereby introduce 
\begin{equation}\label{def: t_tilde intro}
        \widetilde{\TT} \coloneqq 4m \arctanh{\left( \frac{\TT}{\sqrt{1 + \XX^2}}\right)} \; \in \left( -\infty, + \infty \right)\,,
    \end{equation}
which satisfies the \emph{time-dependent Bartnik conditions} in the region $\{ r > \kappa_0m\}$ for any $\kappa_0 > 0$, as in \cref{thm: new time t_b intro}-$(i)$; this then holds also after gluing it to $t-u_{\bb}(\xx)$.
We highlight that the time function $\widetilde{\TT}$ is interesting on its own as
 its time slices are asymptotic to the ones of the static time $t$. Moreover, $\widetilde{\TT}$ seems to be the most natural \emph{Cauchy time function} describing the global hyperbolicity of this spacetime. 

 Finally, we can prove that the glued function $\ttt_{\bb}$ satisfies $(iii)$, which is highly nontrivial. The \emph{Bartnik interior condition} 
requires uniform control on difference of the chosen time function computed at any given point $\p$ lying at a fixed large radius and at points in the interior of the spacetime, excluding the causal future (or past) of $\p$.   
To this end, we found a monotonically decreasing function  that controls the required quantity with respect to $\widetilde{\TT}$.

   We remark that across all the arguments presented, we will keep explicit control on the mass-parameter $m$, which highlight the scalings of the quantities involved.

\vspace{0.5cm}
\noindent \textbf{Structure of the work.}
In \cref{sec: prelim} we present the preliminary definitions and results, in particular with regard to the Schwarzschild spacetime and to the maximal surface equation. In \cref{sec: barriers} we establish the existence of a coordinate boost admitting barriers in the asymptotically flat ends,  after discussing why the standard choice of coordinate boost does not allow the desired barriers construction. In \cref{sec: Cauchy times} we introduce the new Cauchy time function as in \eqref{def: t_tilde intro}, proving its properties; the lengthy computations that are necessary for such a proof are outlined in \cref{Appendix: Bartink conds computations}.
In the same \cref{sec: Cauchy times}, we show the monotonicity of the quantity controlling the \emph{Bartnik interior condition}; moreover, we glue the new time function to the one generated by the hypersurface of \cref{sec: barriers}. Finally, combining all these results, we establish \cref{thm: new time t_b intro}. In \cref{sec: height estimate and main proof} we construct a sequence of suitable boundary data for the sequence of compact Dirichlet problems and we prove the desired a priori height estimate. Finally, we 
prove the main \cref{thm: Main intro}. 

\vspace{0.4cm}

\textbf{Acknowledgments.} The content of this manuscript is a revised version of the main result presented in the author's PhD thesis \cite{CogoPhD}.
I am deeply thankful to Prof. Gerhard Huisken for suggesting the study of this problem, for his supervision and helpful discussions. I am grateful to Prof. Piotr Chruściel for many interesting remarks. I am thankful  
to Prof. Carla Cederbaum for her interest in this work and for suggesting some adjustments on the first version of \cref{prop: Bartik_spactime_cond_graphs_exterior}. Finally, I am grateful to Andoni Royo Abrego and Stephen Lynch for reading and providing comments on parts of this manuscript.

\vspace{0.3cm}

\textbf{Declarations.} This manuscript has been conceived and written without AI support. The contents are mostly as in the author's PhD thesis \cite{CogoPhD}, submitted in May 2025.


\section{Setup and preliminary results}\label{sec: prelim}

When not differently indicated, we will use notations and definitions mostly consistent with 
\cite{Bart}. The first part of this section is dedicated to preliminaries on general spacetimes. A \emph{spacetime} $(\mathfrak{L}, \mathfrak{g})$ is a $4$-dimensional smooth manifold endowed with a smooth Lorentzian metric. 

\vspace{0.1cm}
 
\begin{definition}
A codimension-one smooth submanifold $M \hookrightarrow \mathfrak{L}$ is said to be a \emph{spacelike hypersurface} if it is embedded and at each point, the unit normal to $M$ is a timelike vector.
\end{definition}

\begin{definition}\label{def: timeFct}
    A smooth function $\ttt : \mathfrak{L} \to \nR$ is called \emph{time function}\footnote{For simplicity, in our definition we include the properties required for the time function to be a \emph{temporal function}, see \cite{BernalSanchez_CauchyTimeFunction}. } if $\dd \ttt |_{\p} \neq 0$ and its gradient $\operatorname{grad}_{\mathfrak{g}} \ttt|_{\p}$ is a timelike vector for any $\p \in \mathfrak{L}$. 

    Its level sets $\mathcal{S}_{\ttt_o} = \{ \ttt = \ttt_o\}$ are smooth hypersurfaces called \emph{time slices}.

    The function $\ttt$ is called a \emph{Cauchy time function} if all the time slices $\mathcal{S}_{\ttt_o}$ are \emph{Cauchy surfaces}, namely they are intersected exactly once by any inextendible timelike curve.
\end{definition}

\noindent Notice that the vector field $\operatorname{grad}_{\mathfrak{g}} \ttt$ also defines a \emph{time orientation} on $(\mathfrak{L}, \mathfrak{g})$, so that $-\operatorname{grad}_{\mathfrak{g}} \ttt$ is future-pointing.
We define the \emph{lapse} of the time function $\ttt$ to be 
\begin{equation}\label{def: lapse}
    \alpha_{\ttt}^2 \coloneqq - \langle \,\operatorname{grad}_{\mathfrak{g}} \ttt, \,\, \operatorname{grad}_{\mathfrak{g}} \ttt \,\rangle^{-1}
\end{equation}
and this implies that the \emph{future-pointing unit normal} to the time slices $\mathcal{S}_{\ttt}$ is given by
\begin{equation}\label{def: futureUnitNormal}
    \mathcal{T} \coloneqq - \alpha_{\ttt} \, \operatorname{grad}_{\mathfrak{g}} \ttt \,.
\end{equation}

Throughout this work, we will use the \emph{Einstein summation convention} with indexes. The greek indexes will refer to spatial and time components ranging $1 \leq \,\lambda, \mu, \nu \,\leq 4$, the latin indexes will refer to spatial components ranging $1 \leq \, i, j, k, a, b, c \, \leq 3$ (in particular, we will use $a, b, c$ for local frames on spacelike hypersurfaces and $i, j, k$ specifically for frames on time slices),
while the capital indexes will range $1 \leq \, I, J, K, L \, \leq 2$ and often refer to components on submanifolds of spacelike hypersurfaces.

\vspace{0.3cm} \noindent     \textbf{The geometry of spacelike hypersurfaces.}
Let $M \hookrightarrow (\mathfrak{L}, \mathfrak{g})$ spacelike with induced metric $g$.  
Let 
$\{ \partial_a\}_{a= 1}^3$ a local coordinate frame on $M$ and 
$\{\mathbf{e}_a\}_{a= 1}^3$ a local $g$-orthonormal frame on $M$. Given a smooth function $f \in C^{\infty}(M)$, the \emph{gradient} of $f$ is
$$\operatorname{grad}_M f = 
g^{ab} \,\partial_b (f) \, \partial_a = 
\sum_{a=1}^3\mathbf{e}_a (f) \,  \mathbf{e}_a\,. $$
Given a vector field $\mathfrak{X} \in \Gamma(T \mathfrak{L})$, the convention for the operators $\operatorname{div}_M, \Delta_M$ of \emph{divergence} and \emph{Laplacian} on $M$ is given by 
\begin{align*}
    \operatorname{div}_M \mathfrak{X} & \coloneqq \operatorname{tr}_M (\nabla \mathfrak{X}) 
    = g^{ab} \; \mathfrak{g}\left( \nabla_{\partial_a} \mathfrak{X}, \,\partial_b \right) 
    = \sum_{a=1}^3 \mathfrak{g}\left(  \nabla_{ \mathbf{e}_a} \mathfrak{X},\,  \mathbf{e}_a \right) \,,
\\
   \Delta_M f& \coloneqq \operatorname{div}_M  (\operatorname{grad}_M f) 
 = g^{ab} \; \mathfrak{g}\left( \nabla_{\partial_a} \operatorname{grad}_M f, \,\partial_b \right)
 = \sum_{a=1}^3 \mathbf{e}_a \left(  \mathbf{e}_a (f) \right) - \left(\nabla_{ \mathbf{e}_a}  \mathbf{e}_a \right) (f).
\end{align*}
where $\nabla$  is the Levi-Civita connection of $\mathfrak{g}$ and $\operatorname{tr}_M$ is the trace operator with respect to $g$.

Consider the future-pointing unit normal $\nnu$ to the hypersurface $M \hookrightarrow (\mathfrak{L}, \mathfrak{g})$. We define its \emph{second fundamental form} acting on $X, Y \in \Gamma(TM)$  as
$\mathrm{K}(X, Y) \coloneqq \mathfrak{g}\left( X, \nabla_Y \nnu \right)  = - \mathfrak{g}\left( \nabla_X Y,\,  \nnu \right)$
and its \emph{mean curvature}  as
\begin{align}\label{def: meanCurv}
    \HH(M)
    &\coloneqq \operatorname{tr}_M \mathrm{K} = \operatorname{div}_M \boldsymbol{\mathlarger \nu}\\  \nonumber
    & = g^{ab} \; \mathfrak{g} \left( \nabla_{\partial_a} \boldsymbol{\mathlarger \nu}, \,\partial_b \right)
     = \sum_{a=1}^3 \mathrm{K}(\mathbf{e}_a, \, \mathbf{e}_a) \,,
\end{align}
In this work, such curvature tensors on time slices $\mathcal{S}_{\ttt}$ will be denoted by $\mathrm{K}^o$ and $\HH^o$, respectively.

Given smooth submanifolds $\Sigma \hookrightarrow (\mathfrak{L}, \mathfrak{g})$ of codimension-$2$, we define its \emph{spacetime mean curvature} as the section of the $2$-dimensional normal bundle of $\Sigma$ in $(\mathfrak{L}, \mathfrak{g})$ given by
\begin{align}\label{def: spacetimeMeanCurv}
    \vec{\mathcal{H}}(\Sigma) & \coloneqq \operatorname{tr}_{\Sigma} \vec{\mathcal{K}}\,, 
\end{align}
where $\operatorname{tr}_{\Sigma}$ is the trace operator with respect to the induced metric on $\Sigma$ and $\vec{\mathcal{K}}$ is the \emph{vector-valued second fundamental form}, which acts on $X, Y \in \Gamma(T\Sigma)$ as  $\vec{\mathcal{K}}(X, Y) \coloneqq \left( \nabla_X Y\right)^{\perp}$,
where  ${\mathfrak{X}}^{\perp}$ indicates the orthogonal projection of a vector $\mathfrak{X} \in T_{\p}\mathfrak{L}$ onto $T_{\p}^{\perp}\Sigma$.

\subsection{Maximal hypersurfaces}

\begin{definition}
    A spacelike hypersurface $M$ is said to be \emph{maximal} if $\HH(M) = 0$.
\end{definition}
\noindent
The vanishing of the mean curvature as geometric condition can be translated into a partial differential equation. To this end, once fixed a time function $\ttt$, we consider the \emph{height function} $u \in C^{\infty}(M)$ of a spacelike hypersurface $M$ to be  
\begin{align}\label{def: heightFct}
    u: M \to \nR\,, \qquad & \quad
     u(\p) \coloneqq \ttt |_M(\p)\,.
\end{align}
We then consider the quantity
\begin{equation}\label{def: v}
     \vv_{\ttt} := - \langle \boldsymbol{\mathlarger{\nu}}, \, \mathcal{T} \rangle\,,
\end{equation}
which is related to the gradient of $u$ in the following sense
\begin{equation}\label{eq: grad_u_v}
    \operatorname{grad}_M u = \operatorname{grad}_{\mathfrak{g}}  \ttt + \alpha_{\ttt}^{-1} \, \vv_{\ttt} \, \boldsymbol{\mathlarger{\nu}} = \alpha_{\ttt}^{-1} \left( \vv_{\ttt} \, \boldsymbol{\mathlarger{\nu}} -\mathcal{T} \right)\,,
\end{equation}
which, particular, gives $\vv_{\ttt} \geq 1$ as a consequence of 
\begin{equation}\label{eq: alpha_gradu_v}
    \vert \operatorname{grad}_M u \vert^2_{\mathfrak{g}} =  \alpha_{\ttt}^{-2} \left( \vv_{\ttt}^2 - 1\right)\,.
\end{equation}
Moreover, the identity \eqref{eq: grad_u_v}, combined with the definition \eqref{def: meanCurv}, implies the following expression of the \emph{maximal surface equation}, which will be computationally useful:   
\begin{equation}\label{eq: MeanCurv_forIntEstimate}
    0 = \HH(M)\, \vv_{\ttt}  = \operatorname{div}_M (\alpha_{\ttt} \, \operatorname{grad}_M u ) + \operatorname{div}_M \mathcal{T}\,,
\end{equation}
and after extending $u$ to $\LL$ by requiring $\partial_{\ttt} u = 0$, one can compute (see \cite[Section 2]{Bart})
\begin{equation}\label{eq: divT}
    \operatorname{div}_M \mathcal{T} = \HH^o + (\vv_{\ttt}^2 -1) \left(\K^o \left(\frac{ \DDD u}{\vert \DDD u\vert}, \frac{ \DDD u}{\vert  \DDD u\vert} \right) - \alpha_{\ttt}^{-1} \mathcal{T}(\alpha_{\ttt})\right) + \dd \alpha_{\ttt} ( \operatorname{grad}_M u )\,, 
\end{equation}
where we denoted by $\DDD u \coloneqq \operatorname{grad}_{\mathcal{S}_{\ttt}} u$ the gradient of $u$ on the time slices.

In order to analyze the geometric partial differential equation underlying $\HH(M) = 0$, it is useful to write the equation in terms of geometric quantities of the time slices.
We consider local coordinates $(\ttt, x^i)$ and the associated local frame $\{ \partial_{\ttt}, \partial_i \}_{i = 1}^3$, with vanishing shift, so that $\mathcal{T} = \alpha_{\ttt}^{-1} \partial_{\ttt}$ and the metric takes the form 
$\mathfrak{g} = - \alpha_{\ttt}^2 \dd \ttt^2  + {}^{\ttt}g_{ij}  \dd x^i \dd x^j$,
where ${}^{\ttt}g$ is the induced metric on the time slices. Since the maximal hypersurface, after extending $u$ to $\LL$ as above, can be written as a level set $M = \{ \ttt - u = 0\}$, its future-pointing unit normal vector field is the normalization of $- \operatorname{grad}_{\mathfrak{g}} (\ttt - u) = \alpha_{\ttt}^{-1} \mathcal{T} + \mathfrak{g}^{\mu i} \partial_i u \, \partial_{\mu} = \alpha_{\ttt}^{-1} \mathcal{T} +  \DDD u $. Hence, the unit normal to $M$ and the quantity \eqref{def: v} are given by 
\begin{equation}\label{splitting_Normal_UT}
    \boldsymbol{\mathlarger{\nu}} = \vv_{\ttt}  \left( \alpha_{\ttt} \,  \DDD u + \mathcal{T} \right)\,, \qquad \vv_{\ttt}  = \frac{1}{\sqrt{1 - \vert \alpha_{\ttt} \,  \DDD u  \vert^2}}\,.
\end{equation}
Now, given the tangent vectors $X_i = \partial_i + \partial_i u \, \partial_{\ttt}$, the induced metric on $M$ reads 
\[ g_{ij} = \mathfrak{g}(X_i, X_j) = {}^{\ttt} g_{ij} - \alpha_{\ttt}^2 \,\partial_i u \,\partial_j u\,, \qquad  g^{ij} =  {}^{\ttt} g^{ij} + \alpha_{\ttt}^2 \vv_{\ttt}^2(u) \;{}^{\ttt}\DDD^i u \, {}^{\ttt}\DDD^j u\]
and, consequently, the \emph{maximal surface equation} has the form
\begin{align}\label{eq: max surf eq}
\nonumber
   0 =  \HH(M) &= \operatorname{div}_M \left( \vv_{\ttt} (\alpha_{\ttt} \,  \DDD u + \mathcal{T})  \right) = -g^{ij}  \,\vv_{\ttt} \;\mathfrak{g}\left(  \alpha_{\ttt} \,  \DDD u + \alpha_{\ttt}^{-1} \partial_{\ttt}
    , \, \nabla_{X_i}  X_j \right) \\ 
    & = a^{ij}\left(u,  \DDD u\right)  \,\DDD_{ij}^2 u + \mathcal{F}\left( u,  \DDD u\right)
\end{align}
where $\DDD_{ij}^2 u = \partial_{ij}^2 u - {}^{\ttt} \Gamma_{ij}^k \partial_k u $ is the Hessian of $u$ and
\begin{equation}\label{coefficient maximal surf eq}
    a^{ij}\left(u,  \DDD u\right)  =    \,\alpha_{\ttt} \vv_{\ttt}(u)\left( {}^{\ttt} g^{ij} + \alpha_{\ttt}^2 \vv_{\ttt}^2(u) \;\DDD^i u \, \DDD^j u\right)\,. 
\end{equation} 
The function $\mathcal{F}$ depends on the geometry of the time slices; more precisely, it holds
\begin{align}\label{eq: MeanCurv_secondPart}
    \mathcal{F}\left( u,  \DDD u\right) = \vv_{\ttt} \HH^o + \vv_{\ttt}(1 + \vv_{\ttt}^2) \, \dd \alpha_{\ttt}\left(  \DDD u\right) + \vv_{\ttt}^3\left( \partial_{\ttt} \alpha_{\ttt} \vert   \DDD u \vert - \alpha_{\ttt}^2 \, \K^o_{ij} {}^{\ttt}\DDD^i u \, {}^{\ttt}\DDD^j u \right).
\end{align}
It is then clear that \eqref{eq: max surf eq} is a \emph{second order quasilinear elliptic} partial differential equation, if a uniform control holds for $\alpha_{\ttt}$ and  for the quantity $\vv_{\ttt}$: considering \eqref{coefficient maximal surf eq}, given $\zeta \in \nR$, it holds
\[ {}^{\ttt} g^{ij}\zeta_i \zeta_j \;\leq \;\left( {}^{\ttt} g^{ij} + \alpha_{\ttt}^2 \vv_{\ttt}^2(u) \;\DDD^i u \, \DDD^j u \right) \zeta_i \zeta_j \;\leq \; \vv_{\ttt}^2\; {}^{\ttt} g^{ij}\zeta_i \zeta_j\,. \]

\vspace{0.2cm} \noindent     \textbf{Gradient estimates and time-dependent Bartnik conditions.} 
A crucial result of Bartnik in \cite{Bart} consists in deriving an a priori control on $\vv_{\ttt}$ on a maximal hypersurface $M$ directly from an a priori height estimate, provided that certain conditions are satisfied by the chosen time function $\ttt$ and by the boundary $\partial M$. For a spacetime satisfying $\mathfrak{Ric} = 0$, such \emph{time-dependent Bartnik conditions} consist in having a uniform bound on the quantities
\begin{equation}\label{cond: Bartnik_RiemNorm}
    {}^{\ttt} \| \log \alpha_{\ttt} \|_1\,, \; {}^{\ttt}\| \nabla \mathcal{T}\|_1  \, \leq \,C_o \,,
\end{equation}
for $\alpha_{\ttt}$ and $\mathcal{T}$ as in \eqref{def: lapse} and \eqref{def: futureUnitNormal}, where ${}^{\ttt}\| \cdot \|$ is  
a \emph{time-dependent Riemannian norm}. This defined on any tensor field $\boldsymbol{B} \in \Gamma(T_l^m \mathfrak{L})$ by 
\begin{equation}\label{def: timeDep_RiemNorm}
\begin{split}
    {}^{\ttt}\| \boldsymbol{B} \| & := \sup_{\p \, \in \,\mathfrak{L}} \sqrt{ \sum\displaylimits_{\lambda_1, \dots, \lambda_l, \mu_1, \dots, \mu_m} \Big{|}\boldsymbol{B}|_{\p} \left(\mathrm{E}_{\lambda_1}, \dots, \mathrm{E}_{\lambda_l}, \mathrm{E}^*_{\mu_1}, \dots, \mathrm{E}^*_{\mu_m}\right)\Big{|}^2} \\
    {}^{\ttt}\| \boldsymbol{B} \|_N & := \sum_{k = 0}^N {}^{\ttt}\| \nabla^k \boldsymbol{B} \|\,,
    \end{split}
\end{equation}
where $\{ \mathrm{E}_{\mu} \}_{\mu = 1}^{4}$ is an \emph{adapted local orthonormal frame} given by a local orthonormal frame on the time slices $\{ \mathrm{E}_i \}_{i = 1}^{3}$ and the future-pointing unit normal $ \mathrm{E}_{4} \coloneqq \mathcal{T}$; the elements of the dual frame are denoted by $ \{ \mathrm{E}_{\mu}^* \}_{\mu = 1}^{4}$.

 On a  surface $M$ with $\HH(M) = 0$, as a consequence of  \cite[(2.7) and Proposition 2.1]{Bart}, the quantities $u$ and  $\vv_{\ttt}$ satisfy
\begin{equation*}
\begin{split}
    \Delta_M u & = \operatorname{div}_M(\mathcal{T}) \\
    \Delta_M \vv_{\ttt} & = \vv_{\ttt} \, \vert \mathrm{K} \vert^2_{\mathfrak{g}} + \frac{1}{2}  (\nabla_{\boldsymbol{\nu}} \mathcal{L}_{\mathcal{T}}\mathfrak{g})(\mathbf{e}_a, \mathbf{e}_a) -  (\nabla_{\mathbf{e}_a} \mathcal{L}_{\mathcal{T}}\mathfrak{g})(\boldsymbol{\mathlarger{\nu}}, \mathbf{e}_a) - \mathcal{L}_{\mathcal{T}}\mathfrak{g}(\mathbf{e}_a, \mathbf{e}_b) \cdot \K (\mathbf{e}_a,\mathbf{e}_b)\,,
\end{split}
\end{equation*}
where $\{\mathbf{e}_a\}_{a= 1}^3$ a local $g$-orthonormal frame on $M$ and $\mathcal{L}_{\mathfrak{X}}$ denotes the \emph{Lie derivative} in direction $\mathfrak{X} \in \Gamma(T\mathfrak{L})$. A uniform control on \eqref{cond: Bartnik_RiemNorm} allows to suitably estimate the terms on the right hand side and apply an interior maximum argument on the function $\mathrm{e}^{\mathcal{k} u} \vv_{\ttt}$, for a sufficiently large constant $\mathcal{k} = \mathcal{k}(C_o) \gg 0$. In particular, \cite[Theorem 3.1-$(iii)$]{Bart} shows that under an additional uniform control on the spacetime mean curvature \eqref{def: spacetimeMeanCurv} of the boundary, which has to lie on a time slice  
\begin{equation}\label{cond: boundary grad}
    {}^{\ttt}\| \mathcal{H}(\partial M)\| \,\leq \, C_o\,, \qquad u|_{\partial M} = 0\,,
\end{equation}
the maximum and minimum of $\mathrm{e}^{\mathcal{k} u} \vv_{\ttt}$ are, in fact, attained in the interior and  $\vv_{\ttt} \leq 2 \mathrm{e}^{\mathcal{k} \sup_M u}$. By a similar argument and under the same \emph{time-dependent Bartnik conditions}, \cite[Theorem 3.1-$(iv)$]{Bart} provides also a \emph{pointwise gradient estimate} for points $\p \in M$ satisfying the additional requirement 
\begin{equation}\label{cond: poitwise grad}
    u(\p) - \sup_{\partial M} u  \geq c_o > 0\,,
\end{equation}
yielding $\vv_{\ttt}|_{\p} \leq 2 \mathrm{e}^{\mathcal{k} \left(\sup_M u - u(\p) \right)}$ for a sufficiently large $\mathcal{k} = \mathcal{k}(C_o, c_o^{-1})$. Such pointwise estimate is, in turn, essential in the integral estimates for an a priori height control: it will be applied to points at a fixed large radius and the quantity $\sup_M u - u(\p)$ will be uniformly controlled by the \emph{Bartnik interior condition}, which has to hold with respect to the chosen time function.

The \emph{Bartnik interior condition} require a uniform control on the difference between the time function computed at a point $\p$ on a large radius in the asymptotically flat ends and any point in the interior region $\mathfrak{L}_{\mathrm{int}}$ with respect to that radius, excluding the \emph{causal future} $I^+(\p)$ (\emph{or past} $I^-(\p)$)\footnote{The \emph{causal future} $I^+(\p)$ of a point $ \overline \p$ is the collection of points for which there exists a future-directed causal curve from $\p$ to them. The definition of \emph{causal past} $I^-(\p)$ is analogous, considering past-directed causal curves.} of the point.  If $\ttt(\p) \geq 0$, such condition reads
\begin{equation}
    \label{cond: Bartink_Interior}
    \sup_{\mathrm{q} \in \mathfrak{L}_{\mathrm{int}} \setminus I^+(\p)} \, \ttt(\mathrm{q}) - \ttt(\p) \, \leq \,C_o \,,
\end{equation}
for a uniform constant $C_o > 0$ independent of $\p$.


\subsection{The Schwarzschild spacetime}\label{subsec: Schwarzschild}

The \emph{Schwarzschild spacetimes} are a family of exact non-flat static solutions to the Einstein Field Equations $\mathfrak{Ric} = 0$ in vacuum, with zero cosmological constant. In this work, we use three different coordinates that describe the Schwarzschild spacetime of mass parameter $m > 0$: the so-called \emph{Schwarzschildean} and \emph{Isotropic} coordinates capture asymptotic Minkowskianity of the exterior regions, while the \emph{Kruskal--Szekeres} cover the maximally extended spacetime.

We follow the references \cite[Chapters V.3, V.4.1, VI.3.2]{Choquet-Bruhat_BookRelativity} and \cite[Chapters 6.1, 6.4]{Wald}, to which we refer the reader for additional details. See also the original papers by Schwarzschild \cite{Schwarzschild}, Kruskal \cite{Kruskal} and Szekeres \cite{Szekeres}.

\vspace{0.3cm}
\noindent
\textbf{Schwarzschildean coordinates.}
The \emph{exterior Schwarzschild} solution is given by
\begin{equation}\label{def: Schw_Schwarzschildean}
    \begin{split}
        (\mathfrak{L}_m)_{\mathrm{ext}} & = \nR \times \left(2m, \infty \right) \times \nS^{2}\, \ni (t, r, \boldsymbol{\vartheta})\,, \\
        \mathfrak{g}_m & = - h_m(r) \dd t^2 + \frac{1}{h_m(r)} \dd r^2 + r^2 g_{\nS^{2}}\,, \qquad  \quad  h_m(r) = 1 - \frac{2m}{r}  > 0
    \end{split}
\end{equation}
where $g_{\nS^2}$ is the standard round sphere. The coordinate $t$ in exterior Schwarzschild is said to be the \emph{static time}, since the metric \eqref{def: Schw_Schwarzschildean} is a warped product with the positive warping function and the spacial Riemannian metric independent of time. Moreover, in these coordinates, the spacetime is also rotationally symmetric and asymptotically Minkowskian, since $h_m(r) = 1 + O_{\infty}(r^{-1})$.

Considering the radial coordinate in the interval $(0, 2m)$, we obtain the so-called \emph{interior Schwarzschild}
\begin{equation}\label{def: Interior Schw_Schwarzschildean}
    \begin{split}
        (\mathfrak{L}_m)_{\mathrm{int}} & = \nR \times \left(0, 2m\right) \times \nS^{2}\, \ni (t, r, \boldsymbol{\vartheta})\,, \\
        \mathfrak{g}_m & = - h_m(r) \dd t^2 + \frac{1}{h_m(r)} \dd r^2 + r^2 g_{\nS^{2}}\,, \qquad  \quad  h_m(r) = 1 - \frac{2m}{r}  < 0\,,
    \end{split}
\end{equation}
where the causal nature of $\partial_t$ and $\partial_r$ is swapped and the spacetime is no longer rotationally symmetric nor static.

The metric degenerates as $r \to 2m$ and as $r \to 0^+$. The former is just a coordinate singularity, as it will be clear by the construction of the maximally extended Schwarzschild. On the other hand, the latter $\{ r = 0 \}$ is a proper spacetime singularity.

For later use, we remark that in these coordinates the spacelike condition on a graph of a smooth function $u : \Omega \subset  (2m, \infty) \times \nS^2 \to \nR$ reads
\begin{equation}\label{spacelikeCond_Schw}
    (\partial_r u)^2 + \frac{g_{\nS^2}^{II}}{r^2 \, h_m} (\partial_I u)^2 < \frac{1}{h_m^2}\,.
\end{equation}

\vspace{0.4cm} \noindent     \textbf{Isotropic coordinates.} Given the Schwarzschildian coordinates on exterior Schwarzschild \eqref{def: Schw_Schwarzschildean}, we consider the change of coordinates
\begin{equation}\label{changeCoord_Iso_Schw}
    r = \vert \xx \vert \left( 1 + \frac{m}{2 \vert \xx \vert}\right)^2
\end{equation}
where $\vert \xx \vert = \sqrt{\sum_{i=1}^3 (x^i)^2}$. Considering $\vert \xx \vert > \tfrac{m}{2}$, we obtain the static spacetime
\begin{equation}\label{def: Schw_Isotropic}
    \begin{split}
        (\mathfrak{L}_m)_{\mathrm{ext}} &= \nR \times \left(\nR^3 \setminus \overline{B}_{\frac{m}{2}}\right) \; \ni (t, \xx) \\
    \mathfrak{g}_m &= - N_m(\vert \xx \vert)^2 \dd t^2 + \varphi_m(\vert \xx \vert)^4 \delta_{\nR^3}\,, \\
    & \qquad N_m(\vert \xx \vert) = \frac{1-\frac{m}{2\vert \xx \vert}}{1+\frac{m}{2\vert \xx \vert}} \,, \quad   \varphi_m(\vert \xx\vert) = 1+\frac{m}{2\vert \xx \vert} \,.
    \end{split}
\end{equation}
The fact that exterior Schwarzschild is static and rotationally symmetric is evident also in these coordinates. Moreover, $N_m(\xx) = 1 + O_{\infty}\left( \vert \xx \vert^{-1}\right)$ and $\varphi_m(\xx) = 1 + O_{\infty}\left( \vert \xx \vert^{-1}\right)$ give asymptotic Minkowskianity. 

For later reference, we remark that in these coordinates, the spacelike condition on a graph of a smooth function $u: \Omega \subset \nR^3 \setminus \overline{B}_{\frac{m}{2}} \to \nR$ reads 
\begin{equation}\label{spacelikeCond_iso}
    \sum_{i = 1}^3\left( \partial_iu(\xx)  \right)^2 < \frac{\varphi_m^4(\vert \xx \vert)}{N_m^2(\vert \xx \vert)}
\end{equation}
and, using \eqref{eq: max surf eq}-\eqref{coefficient maximal surf eq}-\eqref{eq: MeanCurv_secondPart}, its mean curvature is given by
\begin{equation}\label{eq: MeanCurv_isotropic}
    \HH(u) = N_m \vv_t \Big( \varphi_m^{-4} \delta^{ij} +  N_m^2 \vv_t^2 \, \DDD^iu \, \DDD^j u  \Big) \, \DD_{ij} u + \vv_t (1 + \vv_t^2) \; \dd N_m \left(  \DDD u\right) 
\end{equation}
where $\vv_t = \frac{1}{\sqrt{1 - N_m^2 \vert \DDD u \vert^2_g}}$ for $g = \varphi^{-4} \delta_{\nR^3}$.

\vspace{0.5cm} \noindent     \textbf{Kruskal--Szekeres coordinates.}
The \emph{maximally extended Schwarzschild spacetime} is the Loretzian manifold
\begin{equation}\label{def: KruskalExtension}
    \begin{split}
        \mathfrak{L}_m & = \left\{ (\TT, \XX) \in \nR^2 \,: \, \TT^2 - \XX^2 < 1\right\} \times \nS^2 \\
        \mathfrak{g}_m & = \xi_m(r) \left(-\dd \TT^2 + \dd \XX^2 \right) + r^2 g_{\nS^2}\,, \qquad \qquad  \xi_m(r) = \frac{32 m^3}{r} e^{-\frac{r}{2m}} > 0
    \end{split}
\end{equation}
where the $r \in (0, +\infty)$ is implicitly given by 

 \begin{equation}\label{def: changeCoord_Kruskal_r}
    \XX^2 - \TT^2 = 
    \varrho_m(r) \coloneqq 
    \left( \frac{r}{2m} - 1\right) e^{\frac{r}{2m}} \; \in (-1, +\infty)\,.
\end{equation}
Notice that it is well defined as $\rho_m(r)$ is a strictly increasing function in $r$ and hence it admits an inverse $r(\TT, \XX)$ with values in $(0, + \infty)$. 
Let us define the following open sets 
\begin{align*}
    \text{ \emph{region} }&  \mathrm{I}   & \left\{ (\TT, \XX) \in \mathfrak{L}_m \,:\, \TT^2 < \XX^2, \, \XX > 0\right\} &, \\
    \text{ \emph{black hole region / region} } & \mathrm{II}  &\left\{ (\TT, \XX) \in \mathfrak{L}_m \,:\, \TT^2 > \XX^2, \, \TT > 0\right\} &,\\
    \text{ \emph{region} } & \mathrm{I}'    &\left\{ (\TT, \XX) \in \mathfrak{L}_m \,:\, \TT^2 < \XX^2, \, \XX < 0\right\} &, \\
    \text{ \emph{white hole region /  region} } & \mathrm{II}'   &\left\{ (\TT, \XX) \in \mathfrak{L}_m \,:\, \TT^2 > \XX^2, \, \TT < 0\right\} &.
\end{align*}
It is easy to see that the  \emph{region} 
 $\mathrm{I}$  is isometric to the \emph{region} $\mathrm{I}'$, while  \emph{region} 
 $\mathrm{II}$  to the \emph{region} $\mathrm{II}'$ via the map $(\TT, \XX) \mapsto (-\TT, -\XX)$.
Considering, in combination with \eqref{def: changeCoord_Kruskal_r}, the change of coordinates on $\mathfrak{L}_m \setminus \{ \TT^2 = \XX^2 \}$ given by
\begin{equation}\label{def: changeCoord_Kruskal_t}
    t = 2m \log \left| \frac{\XX + \TT}{\XX - \TT}\right| = 
    \begin{cases}
          4m \arctanh{\frac{\TT}{\XX}} & \text{ in \emph{regions} } \mathrm{I}, \mathrm{I}' \\
          4m \arctanh{\frac{\XX}{\TT}} & \text{ in \emph{regions} } \mathrm{II}, \mathrm{II}'\,,
    \end{cases}
\end{equation}
one can check that \emph{region} $\mathrm{I}$ and \emph{region} $\mathrm{I}'$ are both isometric to the exterior Schwarzschild spacetime \eqref{def: Schw_Schwarzschildean} (and \eqref{def: Schw_Isotropic}), while \emph{region} $\mathrm{II}$ and \emph{region} $\mathrm{II}'$ are both isometric to the interior Schwarzschild spacetime \eqref{def: Interior Schw_Schwarzschildean}, where $r$ and $t$ are the radial Schwarzschildean coordinate and the static time, respectively. The set of points $\{ r = 2m\} = \{ \TT^2 = \XX^2 \}$, at which the metric is now well-defined, is called \emph{bifurcate horizon}: it is given by the union of the null hypersurfaces $\{ \TT = \XX \neq 0\}$ and $\{ \TT = - \XX \neq 0\}$ and the so-called \emph{bifurcation sphere} $\{ \TT = \XX = 0\}$. 

For the sake of the computations, we recall the inverse transformations to \eqref{def: changeCoord_Kruskal_r} and \eqref{def: changeCoord_Kruskal_t} between the Kruskal coordinates and the Schwarzschildean coordinates:
\begin{equation}\label{def: changeCoord_Kruskal_rt_TTXX}
\small
    \begin{split}
        \TT = &\sqrt{\frac{r}{2m}-1} \;e^{\frac{r}{4m}}\, \sinh{\left( \frac{t}{4m}\right)}\,,  \quad     \XX = \sqrt{\frac{r}{2m}-1} \;e^{\frac{r}{4m}}\, \cosh{\left( \frac{t}{4m}\right)}\, \qquad  \text{ in \emph{region} } \mathrm{I} \; \\
        \TT = &\sqrt{\frac{r}{2m}-1} \;e^{\frac{r}{4m}}\, \cosh{\left( \frac{t}{4m}\right)}\,,  \quad     \XX = \sqrt{\frac{r}{2m}-1} \;e^{\frac{r}{4m}}\, \sinh{\left( \frac{t}{4m}\right)}\, \qquad \text{ in \emph{region} } \mathrm{II} \,.
    \end{split}
\end{equation}
The usual map $(\TT, \XX) \mapsto (-\TT, -\XX)$ provides the change of coordinates in the \emph{regions} $\mathrm{I}'$ and  $\mathrm{II}'$. 

\vspace{0.4cm} \noindent     \textbf{Barriers near the spacetime singularities.}
The spacetimes singularities $\left\{ \TT = \pm \sqrt{1 + \XX^2} \right\}$ corresponding to $\{ r = 0\}$ in interior Schwarzschild \eqref{def: Interior Schw_Schwarzschildean}, are said to be \emph{crushing}, as in the definition by Eardley--Smarr \cite{EardleySmarr}. In fact, they are foliated in a neighborhood by the spacelike hypersurfaces $\mathrm{B}_{r_0} \coloneqq \{ r = r_0 > 0\}$ with mean curvature given by \begin{equation*}\label{meanCurv_Barrier_Crushing}
        \HH(\mathrm{B}_{r_0}) = 
        \begin{cases}
            \frac{2r_0 - 3m}{\sqrt{r_0^3(2m - r_0)}}\, & \quad \text{ in \emph{region} } \mathrm{II} \\
            -\frac{2r_0 - 3m}{\sqrt{r_0^3(2m - r_0)}}\, & \quad \text{ in \emph{region} } \mathrm{II}'
        \end{cases}
    \end{equation*} 
    as computed e.g. in \cite{MalecOMurchada} or \cite[Lemma 2.3.1]{CogoPhD}. It is easy to see that, given $0 < r_0 < \tfrac{3}{2}m$, in \emph{region} $\mathrm{II}$ their mean curvature is 
    $0 > \HH(\mathrm{B}_{r_0}) \searrow -\infty$ as $r_0 \searrow 0^+$, while  in \emph{region} $\mathrm{II}'$  there holds 
    $0 < \HH(\mathrm{B}_{r_0}) \nearrow +\infty$ as $r_0 \searrow 0^+$. In particular, labelling the hypersurfaces in \emph{region} $\mathrm{II}$ by $\mathrm{B}_{r_0}^+$ and in \emph{region} $\mathrm{II}'$ by $\mathrm{B}_{r_0}^-$, it is possible to show the following; see e.g. the proof of \cite[Proposition 2.3.2]{CogoPhD}.

\begin{proposition}\label{prop: chrushingSing_barriers}
    In the maximally extended Schwarzschild spacetime $\left( \mathfrak{L}_m, \mathfrak{g}_m \right)$, the hypersurfaces $ \{\mathrm{B}_{r_0}^+\}_{0 < r_0 < \tfrac{3}{2}m} $ in \emph{region} $\mathrm{II}$ and $ \{\mathrm{B}_{r_0}^-\}_{0 < r_0 < \tfrac{3}{2}m} $ \emph{region} $\mathrm{II}'$ are a foliation of  \emph{upper barriers} and \emph{lower barriers}, respectively.
\end{proposition} 

This means that if a compact maximal hypersurface $M$ has boundary $\partial M \subset \{ r > r_0\}$ for some $0 < r_0 < \tfrac{3}{2}m$, then $M \subset \{ r > r_0\}$. Since we are concerned with hypersurfaces with prescribed asymptotics in the asymptotically flat ends, it will not be restrictive to work in the region $\{ r > r_0\}$ for some $0 < r_0 < \frac{3}{2}m$.


\section{Barriers in the asymptotically flat ends}\label{sec: barriers}

Barriers in the asymptotically flat ends are necessary to ensure that the entire maximal hypersurface is asymptotic to a coordinate-dependent boost. 

In \cite[Section 5]{Bart}, such barriers are given by graphs over the zero-time slice of functions $w^{\pm}_{\Lambda}$ corresponding to rotationally symmetric hypersurfaces with mean curvature $\HH_{\eta}(w^{\pm}_{\Lambda}) = \pm \frac{\Lambda}{\mathcal r^{3}}$ in Minkowski $(\nR^{1, 3}, \eta)$, where $\mathcal{r}$ is the radial coordinate in the flat space and $\Lambda > 0$ is a constant. The leading terms in the estimates of the mean curvature of such graphs are the mean curvature in Minkowski $\HH_{\eta}(w^{\pm})$ and the mean curvature of the zero-time slice $\HH^o$: this is computed and discussed in  \cite[Section 5]{Bart} and \cite[Section 2.6]{CogoPhD}. Hence, for such hypersurfaces to be barriers, it is necessary that the mean curvature of the coordinate-dependent boost we aim to approach, which could be taken as zero-time slice, decays with an order higher than the one of $\HH_{\eta}(w^{\pm}_{\Lambda})$. Although the balance between these two quantities can be relaxed by finding hypersurfaces in Minkowski with mean curvature $\HH_{\eta}(w^{\pm}_{\Lambda}) = \pm \frac{\Lambda}{\mathcal r^{2 + \varepsilon}}$ for $\varepsilon > 0$, the mean curvature of the standard coordinate-boost $\{ t = \bb x_1\}$ in isotropic coordinates of the exterior Schwarzschild decays too slowly. On the other hand, given the Definition \ref{def: coordinate_boost} and a result by Meyers \cite[Lemma 5]{Meyers}, we will show that there exists a coordinate-boost that has the desired property, in the spirit of \cite[Lemma 2.3]{Chr1}.  On the other hand, the barriers we construct are different from the ones in \cite{Chr1}, as we want to control their height along fixed radii in the isotropic coordinates, in order to ensure the validity of the \emph{Bartnik interior condition} and of the height estimate.

We first discuss the existence and the properties of rotationally symmetric barriers in Minkowski. Consider a constant $\Lambda > 0$, a small parameter $0 < \varepsilon < 1$ and the positive function
\begin{equation}\label{eq: w}
    w^+_{\Lambda}(\mathcal r) \coloneqq \Lambda \mathcal r^{-\varepsilon}, \qquad \partial_{\mathcal r} w^+_{\Lambda}(\mathcal r) = - \Lambda \varepsilon \, \mathcal r^{-1-\varepsilon}, \qquad \partial^2_{\mathcal r} w^+_{\Lambda}(\mathcal r) =  \Lambda (1 + \varepsilon)\varepsilon \, \mathcal r^{-2-\varepsilon}
\end{equation}
which is spacelike for sufficiently large radii depending on $\Lambda$ and $\varepsilon$. The mean curvature of the correspondent graph in Minkowski computed with respect to the future-pointing unit normal is given by 
\begin{align} \nonumber
    \HH_{\eta}(w_{\Lambda}^+) 
    &=  \frac{1}{\sqrt{1 - (\partial_{\mathcal{r}} w_{\Lambda}^+)^2}} \left[ \frac{\partial_{\mathcal{r}}^2 w_{\Lambda}^+}{1- (\partial_r w_{\Lambda}^+)^2} + \frac{2\, \partial_{\mathcal{r}} w_{\Lambda}^+}{\mathcal r} \right] \\ \nonumber
    & = \Delta_{\delta} w_{\Lambda}^+ + O\left( \mathcal r^{-3 - 2\varepsilon} \right) < 0
\end{align}
where the sign is negative for sufficiently big radii, since the flat Laplacian reads
\begin{align}\label{computation: Laplace_w}
     \Delta_{\delta} w_{\Lambda}^+ = \partial_{\mathcal{r}}^2 w_{\Lambda}^+ + \frac{2 \partial_{\mathcal{r}} w_{\Lambda}^+}{\mathcal{ r}} = - \frac{\Lambda \varepsilon (1 - \varepsilon)}{\mathcal r^{2 + \varepsilon}} < 0
\end{align}
The corresponding lower barrier is simply the following 
\begin{equation}\label{eq: w_-}
    w_{\Lambda}^-(\mathcal r) \coloneqq - w_{\Lambda}^+(\mathcal r)\,, \qquad \HH_{\eta}(w_{\Lambda}^-) = - \HH_{\eta}(w_{\Lambda}^+) > 0\,.
\end{equation}

We now estimate the mean curvature of standard coordinate-boost $\{ t = f_{\bb}(\xx)\}$ with $f_{\bb} = \bb x_1$ in isotropic coordinates of the exterior Schwarzschild, resulting in
\begin{align}\label{eq: decay_MC_f_b} \nonumber
    \Bigg\vert & \HH(f_{\bb})  - \frac{\bb(1 + \bb^2)m}{(1-\bb^2)^{\frac{3}{2}}} \frac{x_1}{|\xx|^3}  \Bigg\vert \leq \left\vert \HH(f_{\bb}) -\frac{\bb}{(1-\bb^2)^{\frac{3}{2}}} \left[ (2 - \bb^2)  \partial_1 N_m- (1- \bb^2) \delta^{ij} \Gamma_{ij}^1 + \bb^2   \Gamma_{11}^1 \right]\right\vert \\
    &  + \left\vert \frac{\bb}{(1-\bb^2)^{\frac{3}{2}}} \left[ (2 - \bb^2)  \partial_1 N_m- (1- \bb^2) \delta^{ij} \Gamma_{ij}^1 + \bb^2   \Gamma_{11}^1 \right] - \frac{\bb(1 + \bb^2)m}{(1-\bb^2)^{\frac{3}{2}}} \frac{x_1}{|\xx|^3} \right\vert \leq   \frac{C_{\bb}m^2}{\vert \xx \vert^3} \,,
\end{align}
for a constant $C_{\bb} > 0$ depending only on the boost-parameter $\bb$.
It is clear that
the decay rate is slower than the one of $\HH_{\eta}(w^{\pm})$. To prove \eqref{eq: decay_MC_f_b}, we consider the 
the mean curvature in isotropic coordinates as in \eqref{eq: MeanCurv_isotropic} on the graph of the function $f_{\bb}$.
We then estimate, on $\{ \vert \xx \vert \geq m\}$, the following quantities
\begin{align}\label{estimates N phi} \nonumber
   0 < 1 - N_m   
   \leq \frac{m}{\vert \xx \vert}\,, &\qquad
 0 < \varphi_m - 1 =  \frac{m}{2 \vert \xx \vert } \,, \\
     \left\vert \partial_i N_m - \frac{mx_i}{|\xx|^3} \right\vert 
    \leq \frac{m^2}{\vert \xx \vert^3} \left( 1 + \frac{m}{4 \vert \xx \vert}\right) \leq  \frac{5m^2}{4\vert \xx \vert^3} \,,
    & \qquad   \partial_i \varphi_m 
    = - \frac{mx_i}{2|\xx|^3}  \,, \\ \nonumber
    \vert g_{ij} - \delta_{ij} \vert  \leq \vert \varphi_m^4 -1 \vert 
    \leq \frac{5m }{\vert \xx \vert } , & \qquad 
    \vert g^{ij}- \delta^{ij} \vert  \leq \vert \varphi_m^{-4} -1 \vert  
    \leq \frac{5m }{\vert \xx \vert } \,.
\end{align}
Finally, the Christoffel symbols of the metric $g$ of the static time slices can be estimated as follows
\begin{equation}\label{estimate christoffel}
    \left\vert \Gamma_{ij}^k +  \frac{m }{|\xx|^3} \left( x_i  \delta_j^k + x_j  \delta_i^k  - x^k  \delta_{ij} \right)\right\vert \leq \frac{3m}{\vert \xx \vert^2} \frac{\vert \varphi_m -1\vert}{ \varphi_m} \leq \frac{3m^2}{2 \vert \xx \vert^3}\,.
\end{equation}
Using that $0 < 1- \frac{N_m^2(\vert \xx \vert)}{\varphi_m^4(\vert \xx \vert)} \leq \frac{6m}{|\xx|}$, we estimate
\begin{align*}
    \left\vert \vv_t(f_{\bb}) - \frac{1}{\sqrt{1 - \bb^2}}\right\vert &= \left\vert \frac{1}{\sqrt{1 - \frac{N_m^2(\vert \xx \vert)}{\varphi_m^4(\vert \xx \vert)} \bb^2}} - \frac{1}{\sqrt{1 - \bb^2}}\right\vert 
    \leq \frac{\bb^2}{(1 - \bb^2)^{\frac{3}{2}}} \frac{6m}{\vert \xx \vert }
\end{align*}
Employing all the above estimates and $\partial_i f_{\bb} = \bb \delta_{i1}$, $\partial_{ij}^2 f_{\bb} = 0$, we obtain
\begingroup
\allowdisplaybreaks
\small
\begin{align}
    \label{eq: Hu_2}
    \Bigg\vert -\vv_t(f_{\bb}) N_m \,  g^{ij} \, \Gamma_{ij}^k  \partial_k f_{\bb} + \frac{\bb m}{\sqrt{1-\bb^2}} \frac{x_1}{|\xx|^3} \Bigg\vert 
    \leq \left\vert \frac{\bb }{\sqrt{1 - \bb^2}} \left(   \delta^{ij} \Gamma_{ij}^1 - \frac{m x_1}{|\xx|^3}  \right)  \right\vert +  \frac{C_{\bb} m^2}{\vert \xx \vert^3}
    & \leq \frac{C_{\bb} m^2}{\vert \xx \vert^3} \\
    \label{eq: Hu_3}
   \Bigg\vert - \vv_t^3(f_{\bb}) N_m^3   \DDD^if_{\bb} \, \DDD^jf_{\bb} \,\Gamma_{ij}^k \partial_k f_{\bb} -\frac{\bb^3m}{(1-\bb^2)^{\frac{3}{2}}} \frac{x_1}{|\xx|^3} \Bigg\vert 
   \leq \left\vert \frac{\bb^3  }{(1-\bb^2)^{\frac{3}{2}}}  \left(   \Gamma_{11}^1 + \frac{m x_1}{|\xx|^3}\right)\right\vert +  \frac{C_{\bb} m^2}{\vert \xx \vert^3}
   & \leq   \frac{C_{\bb} m^2}{\vert \xx \vert^3} \\
   \label{eq: Hu_1}
    \Bigg\vert \vv_t(f_{\bb})\left(1 + \vv^2_t(f_{\bb}) \right)  
    \dd N_m(\DDD f_{\bb})  -  \frac{(2 - \bb^2)\bb m}{(1-\bb^2)^{\frac{3}{2}}} \frac{x_1}{|\xx|^3}\Bigg\vert 
    \leq \Bigg\vert \frac{(2 - \bb^2)\bb }{(1-\bb^2)^{\frac{3}{2}}} \left( \partial_1 N_m   -   \frac{m x_1}{|\xx|^3} \right)\Bigg\vert +  \frac{C_{\bb}  m^2}{\vert \xx \vert^3} 
    & \leq  \frac{C_{\bb}  m^2}{\vert \xx \vert^3} 
\end{align}
\endgroup
where $C_{\bb}$ is a constant depending only on $\bb$.
Combining \eqref{eq: Hu_1}, \eqref{eq: Hu_2} and \eqref{eq: Hu_3} 
in \eqref{eq: MeanCurv_isotropic}, we obtain \eqref{eq: decay_MC_f_b}.

\subsection{Coordinate-dependent boost admitting barriers}
The main aim of this section is proving the following Proposition, which provides a construction of an hypersurface $\mathcal{S}_{\bb} \in \mathcal{B}_{\bb}$ admitting barriers, modifying the standard coordinate-dependent boost $\{ t = f_{\bb}(\xx)\}$ discussed above.

\vspace{0.2cm}

\begin{proposition}\label{prop: barriers_infty}
   For any $\bb \in [0, 1)$, there exists a  spacelike hypersurface $\mathcal{S}_{\bb} \in \mathcal{B}_{\bb}$ and a 
    function  $\bar f_{\bb} \in C^{\infty}(\nR^3 \setminus \overline{B}_{1})$ satisfying
    \small
    \begin{equation}\label{def: C_f}
         \sup_{\yy \in \nR^3 \setminus \overline{B}_{1}} \log (\vert \yy \vert )^{-1} \left( \vert \bar f_{\bb}(\yy) \vert +  \vert \yy \vert \vert \partial \bar f_{\bb}( \yy) \vert + \vert \yy \vert^2 \vert \partial^2\bar f_{\bb}( \yy) \vert + \vert \yy \vert^3  \vert \partial^3\bar f_{\bb}( \yy) \vert \right) = C_{\bar f_{\bb}} 
    \end{equation} 
    \normalsize
    for some constant $ 0 < C_{\bar f_{\bb}} < \infty$
    with the property that, in isotropic coordinates, for any  $0 < \varepsilon < 1$, there exists $K_{\bb} = K_{\bb}(\bb, \varepsilon, C_{\bar f_{\bb}}) \gg 1$ such that
    \begin{equation}\label{S_b u_b admitting barriers}
        \mathcal{S}_{\bb} \cap \{ \vert \xx \vert > K_{\bb} m\} = \{ t = u_{\bb}(\xx)\}\,, \qquad u_{\bb}(\xx) \coloneqq \bb x_1 + m\bar f_{\bb}\left(\frac{\xx}{m} \right)
    \end{equation}
    and such that for any $K \geq K_{\bb}$  
    there exist hypersurfaces $\mathcal{W}^{\pm}_{\bb, K, \varepsilon}$ satisfying
     \begin{align}\label{def: barriers_K} 
        \mathcal{W}^{\pm}_{\bb, K, \varepsilon} &= \{ t = w^{\pm}_{\bb, K, \varepsilon}(\xx)  \}  \; \text{ spacelike on }  \; \{ \vert \xx \vert \geq Km\}\,,  
        \\ 
        \nonumber
        \HH(w^+_{\bb, K, \varepsilon}) < 0\,, &\qquad \; 0 < w^{+}_{\bb, K, \varepsilon} - u_{\bb} = O(\vert \xx \vert^{-\varepsilon})\,, \qquad \; (w^{+}_{\bb, K, \varepsilon} - u_{\bb})\big|_{\{ \vert \xx \vert = Km\}} \geq K^{ \varepsilon} \, m \,,   \\ \nonumber
        \HH(w^-_{\bb, K, \varepsilon}) > 0\,, &\qquad \; 0 <  u_{\bb} - w^{-}_{\bb, K, \varepsilon}  = O(\vert \xx \vert^{-\varepsilon})\,,  \qquad \; (u_{\bb}- w^{-}_{\bb, K, \varepsilon} )\big|_{\{ \vert \xx \vert = Km\}}\geq K^{ \varepsilon} \, m \,. 
    \end{align}
    Moreover, there exists a constant $C_{\bb}$ independent of $K$ such that for $\zz(\xx) = \left( \sqrt{1 - \bb^2} x_1 , x_2, x_3 \right)$, it holds $\left\vert  \HH(w_{\bb, K, \varepsilon}^{\pm})  \pm  \frac{K^{2\varepsilon} \varepsilon (1- \varepsilon) m^{1 + \varepsilon}}{ \vert \zz (\xx) \vert^{2 + \varepsilon}} \right\vert \leq \frac{C_{\bb} K^{2\varepsilon} m^2 \log^2 \left( \frac{\vert \xx \vert}{m}\right)}{\vert \xx \vert^{3}}$
     as $\vert \xx \vert \to \infty$.
\end{proposition}

\vspace{0.3cm}

The proof of \cref{prop: barriers_infty} will rely on the following result, which is a consequence of \cite[Lemma 5, Theorem 2]{Meyers}.

\vspace{0.2cm}

\begin{lemma}[{\cite[Proposition 2.1]{Chr1}}]\label{lemma: Mayers}
    Let $\phi \in C^{k, \alpha}_{-2}(\nR^3 \setminus \overline{B}_1)$ \footnote{\label{def: weighted Holder} For the definitions of \emph{weighted H\"older spaces} see e.g. \cite[Definition A.22]{LeeGeometricRelativity}.}
     for some $k \geq 0$, $\alpha \in (0,1)$.  
    Then there exists a solution $\bar f \in C^{k+2, \alpha}_{loc}(\nR^3 \setminus \overline{B}_1)$ to 
    \[ \Delta_{\delta} \bar f = \phi\, \qquad  \text{ on } \nR^3 \setminus \overline{B}_1\]
    such that $\bar f = O_{k+2}( 1+ \log( \vert{\yy} \vert))$ for $\vert {\yy} \vert \to \infty$.
\end{lemma}

\vspace{0.1cm}

\begin{proof}
 The existence of a solution satisfying $\bar f = O (\log |\yy|)$ is ensured by  \cite[Lemma 5]{Meyers}, while the behavior of higher derivatives can be shown by rescaling and employing the standard \emph{interior Schauder's estimate}.
 For completeness and in the absence of the proof of \cite[Proposition 2.1]{Chr1}, we briefly outline the full argument in our specific setup.

We can assume $\phi$ vanishes in a neighborhood of $B_1$, as the solutions would just differ by a function that is harmonic  outside such a neighborhood and decaying as $O(\vert \yy \vert^{-1})$. We then consider the following solution to the above Poisson equation, for points $|\yy| > 1$
 \begin{equation}\label{eq: solution bar f}
     \bar f (\yy) = \int_{\vert \zz \vert > 1} \left[\frac{1}{4\pi \vert \zz - \yy \vert} - \frac{1}{4 \pi \vert \zz \vert}\right] \phi(\zz)\,\dd \zz 
 \end{equation}
  and let us show that it satisfies the desired decay.   
  Let us split the function into three parts, corresponding to integration over $\{1 < \vert \zz \vert < \frac{ \vert \yy \vert}{2}\}$, $\{ \frac{ \vert \yy \vert}{2} < \vert \zz \vert <  2 \vert \yy \vert \}$ and $\{ 2 \vert \yy \vert  < \vert \zz \vert \}$, respectively.
 In the following, we will use $C > 0$ to denote positive constants.
Using that in the first region $ \frac{  \vert \yy \vert}{2} < \vert \zz - \yy \vert$, we then estimate the first term 
 \begin{align*}
    \left\vert \;\int\displaylimits_{1 < \vert \zz \vert < \frac{ \vert \yy \vert}{2} } \frac{\phi(\zz)}{4\pi \vert \zz - \yy \vert}  \,\dd \zz  \right\vert    
    & \leq C  \frac{2}{\vert \yy \vert} \int\displaylimits_{ 0 < \vert \zz \vert < \frac{ \vert \yy \vert}{2} } \frac{\dd \zz}{\vert \zz  \vert^2}   \leq C \,.
 \end{align*}
 The first term in the second region is instead estimated as follows
 \begin{align*}
      \left\vert \;\int\displaylimits_{ \frac{ \vert \yy \vert}{2} < \vert \zz\vert < 2 \vert \yy \vert} \frac{\phi(\zz)}{4\pi \vert \zz - \yy \vert} \dd \zz  \right\vert  & \leq C  \frac{4}{\vert \yy \vert^2} \int\displaylimits_{\frac{ \vert \yy \vert}{2} < \vert \zz\vert < 2 \vert \yy \vert} \frac{\dd \zz}{\vert \zz - \yy \vert }  \leq C\,,
 \end{align*}
 as the integral of the fundamental solution on this interval is bounded by a constant times $\vert \yy \vert^2$, which is easily found by scaling. We can then estimate the second term of \eqref{eq: solution bar f} in the first two regions  as follows
 \begin{equation*}
    \int_{1 < \vert \zz \vert < 2\vert \yy \vert} \frac{\vert \phi(\zz) \vert }{4 \pi \vert \zz \vert} \dd \zz \leq C  \int_{1 < \vert \zz \vert < 2\vert \yy \vert} \frac{\dd \zz }{ \vert \zz \vert^3} \leq C \log(\vert {\yy} \vert )\,.
 \end{equation*}
 Finally, in the third region where $\frac{\vert \yy \vert }{ \vert \zz \vert } < \frac{1}{2} < 1$, the fundamental solution can be expanded in the following absolutely uniformly convergent series $\frac{1}{\vert \zz - \yy \vert} = \sum_{j = 0}^{\infty} \vert \yy\vert^j \vert \zz \vert^{-1 - j} P_j\left( \frac{\zz}{\vert \zz \vert} \cdot  \frac{\yy}{\vert \yy \vert} \right) $\,, where $P_j$ are the  \emph{Legendre polynomials} computed on the angular coordinates. We notice that since $P_0 = 1$, the term $\frac{1}{4\pi \vert \zz \vert}$ in \eqref{eq: solution bar f} is exactly the first term of the series.  Using that $\vert P_j \vert \leq 1 $, we  estimate
    \begin{align*}
      \int\displaylimits_{ 2 \vert \yy \vert < \vert \zz\vert } \left\vert  \frac{1}{4\pi \vert \zz - \yy \vert} - \frac{1}{4\pi \vert \zz \vert }\right\vert \vert \phi(\zz) \vert \dd \zz    &  \leq C \sum_{j = 1}^{\infty} \frac{1}{j \,2^j } \leq C \,.
 \end{align*}
 Combining all the above estimates, we obtain  $\bar f = O( \log( \vert{\yy} \vert))$ for $\vert {\yy} \vert \to \infty$.

 To show the decay of the higher derivatives, we consider any radius $R > 1$ and the function $\bar f_R(\yy) \coloneqq \bar f(R \yy)$ defined on $\vert \yy \vert > 1$. Such a rescaled function satisfies $\Delta_{\delta} \bar f_R(\yy) = R^2 \phi(R \yy)$ and by the \emph{interior Schauder's estimate} (see e.g. \cite[Theorem 6.2]{GT}),  
 on the annulus $A_1 = \{ 1 < \vert \yy \vert < 2\}$, there exists a constant $C = C(\alpha) > 0$ such that 
 \begin{align*}
     \| \bar f_R \|^*_{C^{2, \alpha}(A_1)} & \leq C \left( \| \bar f(R \; \cdot \;)\|_{L^\infty(A_1)} + \| R^2 \phi(R \;\cdot\;)\|^{(2)}_{C^{0, \alpha}(A_1)}   \right)  \leq  C (1 + \log R)\,,
 \end{align*}
 where the norms $\| \, \cdot \, \|^*_{C^{2, \alpha}(A_1)}$ and $ \|\, \cdot \,\|^{(2)}_{C^{0, \alpha}(A_1)} $  are as in \cite[(6.10) Section 6.1]{GT}.
 Since $\partial^l \bar f_R(\yy) = R^l \,\partial^l \bar f (R \yy)$ 
 and since $R > 1$ is arbitrary, we obtain, by iterating this procedure, $\bar f = O_{k +2}( 1+ \log( \vert{\yy} \vert))$ for $\vert {\yy} \vert \to \infty$. 
\end{proof}

\begin{proof}[Proof of \cref{prop: barriers_infty}]
    Consider a solution $\bar f_{\bb} \in C^{\infty}(\nR^3 \setminus \overline{B}_{m})$ 
    to
    \begin{equation}\label{def: tilde_f}
        \left(\sum_{i=1}^3 \frac{\partial^2}{\partial x_i^2} + \frac{\bb^2}{1-\bb^2} \frac{\partial^2}{\partial x_1^2}\right) \bar f_{\bb}\left( \frac{\xx}{m} \right)  = - \frac{\bb(1 + \bb^2)}{(1-\bb^2)} \frac{x_1}{|\xx|^3}\,,
    \end{equation}
    where the right-hand side is the leading term of $- \frac{\sqrt{1 - \bb^2}}{m} \,
    \HH(f_{\bb})$ as computed in \eqref{eq: decay_MC_f_b}. Such a solution exists with  $\bar f_{\bb}$ satisfying \eqref{def: C_f}
    as it is sufficient to change the coordinates 
    \begin{equation}\label{def: changeCoords_xz}
        z_1 = \sqrt{1-\bb^2} x_1, \quad z_2 = x_2, \quad z_3 = x_3
    \end{equation}
    and consider for $\yy = \frac{\zz}{m}$ the solution to the Poisson equation 
    \begin{align*}
        \Delta_{\delta} \bar f_{\bb}(\yy) = 
    - \frac{\bb(1 + \bb^2)}{\sqrt{1 - \bb^2}} \frac{y_1}{\left(y_1^2 + (1-\bb^2) (y_2^2 + y_3^2)\right)^{\frac{3}{2}}}\,,
    \end{align*}
    ensured by \cref{lemma: Mayers}. Consider now on large radii $\vert \xx \vert \gg m$ the following hypersurface
    \begin{equation}\label{def: S_b}
        \mathcal{S}_{\bb} \coloneqq \{ t = u_{\bb}(\boldsymbol{x})\}\,, \qquad  u_{\bb}(\boldsymbol{x}) \coloneqq f_{\bb}(\boldsymbol{x}) + m\bar f_{\bb}\left( \frac{\xx}{m} \right)
    \end{equation}
    where $f_{\bb} = \bb x_1$ as before and $\bar f_{\bb}$ as in \eqref{def: tilde_f}, both in isotropic coordinates. For any fixed $0 < \varepsilon < 1$, we claim that $\mathcal{S}_{\bb} \in \mathcal{B}_{\bb}$ is the desired hypersurface admitting barriers on $\{ \vert \xx \vert \geq Km\}$ for any big enough
    $K > K_{\bb}$ with $ K_{\bb} =  K_{\bb} (\bb, \varepsilon, C_{\bar f_{\bb}}) \gg 1 $ to be determined, given by
    \begin{equation}\label{def: barriers_infty}
        \mathcal{W}_{\bb, K, \varepsilon}^{\pm} \coloneqq \left\{ t = w_{\bb, K, \varepsilon}^{\pm}(\xx) \right\}\,, \qquad  w_{\bb, K, \varepsilon}^{\pm}(\xx) \coloneqq u_{\bb}(\xx) + 
        m\, w^{\pm}_{K^{2 \varepsilon}}\left( \frac{\vert{\zz (\xx)}\vert}{m} \right) \,,
    \end{equation}
    with $ w^{\pm}_{ K^{2 \varepsilon}}$ as in \eqref{eq: w}--\eqref{eq: w_-}  and the coordinates $\zz$ as in \eqref{def: changeCoords_xz}.

    As immediate consequences of the definition of $\bar f_{\bb}$ in \eqref{def: tilde_f} through \cref{lemma: Mayers}, the function $\bar f_{\bb}$ satisfies \eqref{def: C_f} and
    $\mathcal{S}_{\bb}$ defined in \eqref{def: S_b} is in the class $\mathcal{B}_{\bb}$ as in \cref{def: coordinate_boost}. It is also immediate to see that 
    \[ 0 < w^{+}_{\bb, K, \varepsilon} - u_{\bb} = O(\vert \xx \vert^{-\varepsilon})\,, \qquad 0 < u_{\bb} - w^{-}_{\bb, K, \varepsilon}  = O(\vert \xx \vert^{-\varepsilon})\,.\]
    Moreover, we can check that on $\{\vert \xx \vert = Km \}$, it holds
    \begin{equation*}
        ( u_{\bb} - w^{-}_{\bb, K, \varepsilon})\big|_{\{ \vert \xx \vert = Km\}} 
        = \frac{K^{2\varepsilon} \,m^{1+ \varepsilon}}{\vert \zz (\xx) \vert^{\varepsilon}} \Bigg|_{\{ \vert \xx \vert = Km\}}\geq \frac{K^{2\varepsilon} \,m^{1+ \varepsilon} }{\vert \xx \vert^{\varepsilon}} \Bigg|_{\{ \vert \xx \vert = Km\}} = K^{ \varepsilon}  m 
    \end{equation*}
    and the same estimate holds for $(w^{+}_{\bb, K, \varepsilon} - u_{\bb})$ on $\{ \vert \xx \vert = Km\}$.
     In the above inequalities we have used that $\vert {\boldsymbol{z}(\xx)} \vert^2  = z_1(\xx)^2 + z_2(\xx)^2 + z_3(\xx)^2 = (1 - \bb^2) x_1^2 + x_2^2 + x_3^3$, which implies
    \begin{equation}\label{eq: comparison_x_z}
        (1 - \bb^2)\vert \xx\vert^2  \leq \vert {\boldsymbol{z}(\xx)} \vert^2 \leq \vert \xx\vert^2\,.
    \end{equation}

     We are  left showing that there exists $K_{\bb} = K_{\bb}(\bb, \varepsilon, C_{\bar f_{\bb}}) \gg 1$ such that for all $K > K_{\bb}$ the barriers defined as in \eqref{def: barriers_infty} are spacelike hypersurfaces with the desired sign of the mean curvature on $\{ \vert \xx \vert \geq Km\}$. 
      Using \eqref{def: C_f}, \eqref{eq: w} and \eqref{eq: comparison_x_z} and denoting by $\partial_i = \frac{\partial}{\partial x_i}$ the partial derivative in the isotropic coordinates, we estimate
    \begin{align}\label{computation: derivative_w_bK} \nonumber
        \vert \partial_i w^{\pm}_{\bb, K, \varepsilon}(\xx) \vert & = \left\vert \bb \,\delta_{i1} + m\,\frac{\partial}{\partial x^i}\bar f_{\bb}\left( \frac{\xx}{m} \right)  \mp K^{2 \varepsilon} \,\varepsilon \, m^{1+\varepsilon}\frac{ x_i (1 - \bb^2 \delta_{i1})}{\vert \zz (\xx) \vert^{2 +\varepsilon}} \right\vert  \\
        & 
        \leq \bb \, \delta_{i1} + \frac{C_{\bar f_{\bb}}  m \log\left( \frac{|\xx|}{m} \right)}{\vert \xx \vert} + \frac{K^{2 \varepsilon} \,\varepsilon \, m^{1+\varepsilon}}{(1 - \bb^2)^{\frac{1 + \varepsilon}{2}} \vert \xx \vert^{1 + \varepsilon}}
        \\ 
        \nonumber
        & 
        \leq \bb \, \delta_{i1} + \frac{C_{\bar f_{\bb}} \log K}{K} + \frac{\varepsilon }{(1 - \bb^2)^{\frac{1 + \varepsilon}{2}} K^{1 - \varepsilon}}\,.
        \end{align}
        Hence, there exists a sufficiently large $K_{\bb} = K_{\bb}(\bb, \varepsilon, C_{\bar f_{\bb}})$ such that for all $K > K_{\bb}$ it holds
        \[ \sum_{i = 1}^3 (\partial_i w^{\pm}_{\bb, K, \varepsilon}(\xx))^2 < \bb^2 + \frac{1 - \bb^2}{2} < 1\,,
        \]
        which ensures that the spacelike condition \eqref{spacelikeCond_iso} is satisfied on $\{ \vert \xx \vert \geq Km\}$, as $\bb \in [0, 1)$. Moreover, it yields $\vv_t(w^{\pm}_{\bb, K, \varepsilon}) < 
        \sqrt{\frac{2}{1 - \bb^2}} $.

        Now, we estimate the mean curvature. In the following, we will use $C_{\bb}$ to collect all the constants depending only on $\bb$ and independent of $m$. We first need the estimate
        \begin{align}\label{estimate v_t w_bK} \nonumber
    \left\vert \vv_t(w^{\pm}_{\bb, K, \varepsilon}) - \frac{1}{\sqrt{1 - \bb^2}}\right\vert &= \left\vert \frac{1}{\sqrt{1 - \frac{N_m^2(\vert \xx \vert)}{\varphi_m^4(\vert \xx \vert)} \sum_{i= 1}^3 ( \partial_i w^{\pm}_{\bb, K, \varepsilon}(\xx) )^2}} - \frac{1}{\sqrt{1 - \bb^2}}\right\vert \\ \nonumber
    & \leq   
    \sqrt{\frac{2}{(1-\bb^2)^3}}\left( \frac{N_m^2(\vert \xx \vert)}{\varphi_m^4(\vert \xx \vert)} \sum_{i= 1}^3 ( \partial_i w^{\pm}_{\bb, K, \varepsilon}(\xx) )^2 - \bb^2\right) 
    \\ & 
    < \frac{C_{\bb} K^{\varepsilon} m \log\left( \frac{\vert \xx \vert}{m}\right)}{\vert \xx \vert}\,
\end{align}
by combining the mean value theorem with \eqref{computation: derivative_w_bK}. Moreover, we compute
\begin{equation}\label{computation_modifiedLaplacian_tilde_w}
    \left(\sum_{i=1}^3 \frac{\partial^2}{\partial x_i^2} + \frac{\bb^2}{1-\bb^2} \frac{\partial^2}{\partial x_1^2}\right)  w^{\pm}_{K^{2 \varepsilon}}\left( \frac{\vert \zz (\xx) \vert}{m} \right) = \sum_{i=1}^3 \frac{\partial^2}{\partial z_i^2} w^{\pm}_{K^{2 \varepsilon}}\left( \frac{\vert \zz (\xx) \vert}{m} \right)  = \mp \frac{K^{2\varepsilon} \varepsilon (1- \varepsilon) m^{ \varepsilon}}{ \vert \zz (\xx) \vert^{2 + \varepsilon}}\,,
\end{equation}
where we used the change of coordinates \eqref{def: changeCoords_xz}  and the computation \eqref{computation: Laplace_w}. Finally, considering the mean curvature of $w_{\bb, K, \varepsilon}^{\pm}$ given by \eqref{eq: MeanCurv_isotropic} and estimating errors by using \eqref{eq: w}, \eqref{estimates N phi}, \eqref{estimate christoffel}, \eqref{def: C_f} and \eqref{estimate v_t w_bK}
we find that on $\{ \vert \xx \vert \geq K m\}$ there holds
\small
\begin{align}\label{eq: Hw_1} \nonumber
    & \left\vert N_m \vv_t(w_{\bb, K, \varepsilon}^{\pm})  g^{ij}  \DD_{ij} w_{\bb, K, \varepsilon}^{\pm} - \frac{m}{\sqrt{1 - \bb^2}}  \sum_{i=1}^3 \frac{\partial^2}{\partial x_i^2} \left( \bar f_{\bb}\left( \frac{\xx}{m} \right) +  w^{\pm}_{K^{2 \varepsilon}}\left( \frac{\vert \zz (\xx) \vert}{m} \right)\right) - \frac{\bb  \delta^{ij} \Gamma_{ij}^1}{\sqrt{1 - \bb^2}}\right\vert  \\
    & \qquad \qquad \qquad  \leq \frac{C_{\bb} K^{2\varepsilon} m^2 \log^2 \left( \frac{\vert \xx \vert}{m}\right)}{\vert \xx \vert^{3}}
\end{align}
\begin{align}\label{eq: Hw_2}
\nonumber
    & \left\vert N_m^3 \vv_t^3(w_{\bb, K, \varepsilon}^{\pm})  \DDD^iw_{\bb, K, \varepsilon}^{\pm} \DDD^jw_{\bb, K, \varepsilon}^{\pm}  \DD_{ij} w_{\bb, K, \varepsilon}^{\pm} - \frac{\bb^2 m}{(1-\bb^2)^{\frac{3}{2}}}  \frac{\partial^2}{\partial x_1^2} \left( \bar f_{\bb}\left( \frac{\xx}{m} \right) +  w^{\pm}_{K^{2 \varepsilon}}\left( \frac{\vert \zz (\xx) \vert}{m} \right)\right) - \frac{\bb^3 \Gamma_{11}^1 }{(1 - \bb^2)^{\frac{3}{2}}} \,\right\vert \\
    & \qquad \qquad \qquad \qquad \quad \leq \frac{C_{\bb} K^{2\varepsilon} m^2 \log^2 \left( \frac{\vert \xx \vert}{m}\right)}{\vert \xx \vert^{3}}
\end{align}
\begin{align}\label{eq: Hw_3}
      \Bigg\vert \vv_t(w_{\bb, K, \varepsilon}^{\pm})\left(1 + \vv^2_t(w_{\bb, K, \varepsilon}^{\pm}) \right)  
    \dd N_m(\DDD w_{\bb, K, \varepsilon}^{\pm})  -  \frac{(2 - \bb^2)\bb m}{(1-\bb^2)^{\frac{3}{2}}} \frac{x_1}{|\xx|^3}\Bigg\vert 
     \leq  \frac{C_{\bb} K^{\varepsilon} m^2 \log \left( \frac{\vert \xx \vert}{m}\right)}{\vert \xx \vert^{3}}
\end{align}
\normalsize
Now, combining \eqref{eq: Hw_1}, \eqref{eq: Hw_2}, \eqref{eq: Hw_3} with \eqref{eq: decay_MC_f_b}, we find
\begin{align}\label{computation: meanCurv_barriers_infty}
    \nonumber
     \Bigg\vert  & \HH(w_{\bb, K, \varepsilon}^{\pm})  \pm  \frac{K^{2\varepsilon} \varepsilon (1- \varepsilon) m^{1 + \varepsilon}}{ \vert \zz (\xx) \vert^{2 + \varepsilon}} \Bigg\vert \\ \nonumber
    & \leq \Bigg\vert \frac{m}{\sqrt{1 - \bb^2}} \left(\sum_{i=1}^3 \frac{\partial^2}{\partial x_i^2} + \frac{\bb^2}{1-\bb^2} \frac{\partial^2}{\partial x_1^2}\right) \left(  \bar f_{\bb}\left( \frac{\xx}{m} \right) + w^{\pm}_{K^{2 \varepsilon}}\left( \frac{\vert \zz (\xx) \vert}{m} \right)\right) \pm  \frac{K^{2\varepsilon} \varepsilon (1- \varepsilon) m^{1 + \varepsilon}}{ \vert \zz (\xx) \vert^{2 + \varepsilon}}   \\ \nonumber
    & \quad +\frac{\bb}{(1-\bb^2)^{\frac{3}{2}}} \left[    \frac{(2 - \bb^2) mx_1}{|\xx|^3}- (1- \bb^2) \delta^{ij} \Gamma_{ij}^1 + \bb^2   \Gamma_{11}^1 \right]  \Bigg\vert + \frac{C_{\bb} K^{2\varepsilon} m^2 \log^2 \left( \frac{\vert \xx \vert}{m}\right)}{\vert \xx \vert^{3}} \\ \nonumber
    & \leq  \Bigg\vert \frac{m}{\sqrt{1 - \bb^2}} \left(\sum_{i=1}^3 \frac{\partial^2}{\partial x_i^2} + \frac{\bb^2}{1-\bb^2} \frac{\partial^2}{\partial x_1^2}\right) \left(  \bar f_{\bb}\left( \frac{\xx}{m} \right) + w^{\pm}_{K^{2 \varepsilon}}\left( \frac{\vert \zz (\xx) \vert}{m} \right)\right) \pm  \frac{K^{2\varepsilon} \varepsilon (1- \varepsilon) m^{1 + \varepsilon}}{ \vert \zz (\xx) \vert^{2 + \varepsilon}}   \\
    & \quad + \frac{\bb(1 + \bb^2)m}{(1-\bb^2)^{\frac{3}{2}}} \frac{x_1}{|\xx|^3} \Bigg\vert + \frac{C_{\bb} K^{2\varepsilon} m^2 \log^2 \left( \frac{\vert \xx \vert}{m}\right)}{\vert \xx \vert^{3}} \;\leq \, \frac{C_{\bb} K^{2\varepsilon} m^2 \log^2 \left( \frac{\vert \xx \vert}{m}\right)}{\vert \xx \vert^{3}} \,,
\end{align}
where the appropriate cancellations follow from  \eqref{computation_modifiedLaplacian_tilde_w} and the definition of $\bar f_{\bb}$ via \eqref{def: tilde_f}. Again, we highlight that the constant $C_{\bb}$ is independent of $K$, $\varepsilon$ and $m$. 

The estimate \eqref{computation: meanCurv_barriers_infty} implies, in particular, that 
\[ \HH(w_{\bb, K, \varepsilon}^{+})(\xx) < 0 \,, \qquad \HH(w_{\bb, K, \varepsilon}^{-})(\xx) > 0\]
on $\{|\xx| \geq K\}$ for any $K > K_{\bb}$, where $K_{\bb} = K_{\bb}(\bb, \varepsilon, C_{\bar f_{\bb}}) \gg 1$ is sufficiently large. This concludes the proof.
\end{proof}


\vspace{0.3cm}
\section{Cauchy time functions}\label{sec: Cauchy times}

Before introducing a new Cauchy time function in the maximally extended Schwarzschild spacetime, which satisfies the crucial conditions \eqref{cond: Bartnik_RiemNorm}, we first discuss why such \emph{time-dependent Bartnik conditions} do not hold when computed with respect to the Kruskal-Szekeres coordinates \eqref{def: KruskalExtension}.

\vspace{0.2cm}

\noindent \textbf{Time-dependent Bartnik conditions fail on the Kruskal-Szekeres coordinates.}
Consider now the time function $\TT$ of Kruskal-Szekeres coordinates \eqref{def: KruskalExtension} and the associated time-dependent Riemannian norm as defined in \eqref{def: timeDep_RiemNorm}. A suitable local orthonormal frame 
is given by 
\[ \left\{  
\xi_m^{-\frac{1}{2}}\partial_{\TT},\; \xi_m^{-\frac{1}{2}} \partial_{\XX}, \; r^{-1}\sqrt{g_{\nS^2}^{II}} \,\partial_I   \right\}, \quad I=1, 2\,.\]
The lapse function is given by  
$$\alpha_{\TT} = \sqrt{\xi_m} = (2m)^{\frac{3}{2}} \,\frac{e^{-\frac{r}{4m}}}{\sqrt{r}}$$
and it is possible to bound it uniformly $0 < \alpha_{\TT} < \infty$ provided that we restrict to $\{ r > \kappa_0 m  \}$ for some $\kappa_0 > 0$. This is not a restrictive condition, as we can employ the barriers at the crushing singularity as in \cref{prop: chrushingSing_barriers} for a lower bound on the radius.  
On the other hand, the norm ${}^{\TT}\| \log \alpha_{\TT} \|_1$ does not have an uniform bound in the region $\{ r > \kappa_0 m\}$. In fact, we can compute 
\begin{align*}
    \dd \log \alpha_{\TT} &= \alpha_{\TT}^{-1} \frac{\partial \alpha_{\TT}}{\partial r}\left( \frac{\partial r}{\partial \TT} \dd \TT +  \frac{\partial r}{\partial \XX} \dd \XX\right)   =  \frac{\xi_m^{\frac{3}{2}}}{8m^2}\left( 1 + \frac{2m}{r}\right) \left( \TT \dd \TT -  \XX \dd \XX\right)
\end{align*}
where we used \eqref{eq: Kruskal_derivative_lapse} and  \eqref{eq: Kruskal_partial_r_TX}; 
we then obtain
\begin{equation}\label{eq: Bart time-dep fails on Kruskal time}
\begin{split}
    {}^{\TT} \| \olNabla \log \alpha_{\TT} \|^2 &= \sup_{\{ r > \kappa_0 m\}}  \left[ \vert  \dd \log \alpha_{\TT}(\xi_m^{-\frac{1}{2}}\partial_{\TT})\vert^2 + \vert \dd \log \alpha_{\TT}(\xi_m^{-\frac{1}{2}}\partial_{\XX})\vert^2 \right] \\
    & = \sup_{\{ r > \kappa_0 m\}} \; \frac{\xi_m^2(r)}{16m^4}\left( 1 + \frac{2m}{r}\right)^2\,\left( \TT^2 + \XX^2 \right)
    \end{split}
\end{equation}
which fails to be uniformly bounded on region of $\{ r > \kappa_0 m \}$.

\vspace{0.3cm}
\noindent
\textbf{A new Cauchy time function in the maximally extended Schwarzschild.}
While the Kruskal-Szekeres coordinates \eqref{def: KruskalExtension} are necessary to cover the maximally extended Schwarzschild spacetime, they do not capture the asymptotic flatness of the spacetime, which is instead evident in Schwarzschildean \eqref{def: Schw_Schwarzschildean} and isotropic \eqref{def: Schw_Isotropic} coordinates. Motivated by this observation, we consider a foliation of the Kruskal extension $\left( \mathfrak{L}_m, \mathfrak{g}_m\right)$ by the family of spacelike hypersurfaces
\begin{equation}\label{level sets_tildeT}
    \left\{\TT = a \sqrt{1 + \XX^2}, \quad |a| < 1 \right\} \,.
\end{equation}
Since $\lim\displaylimits_{\XX  \to +\infty} \tfrac{ a \sqrt{1 + \XX^2}}{\XX} = a$ and $\lim\displaylimits_{\XX  \to -\infty} \tfrac{ a \sqrt{1 + \XX^2}}{\XX} = -a$, they are asymptotic to the surfaces
\begin{align*}
    \left\{\TT = a  \XX, \quad |a| < 1 \right\} \, = \, \left\{t = 4m \, \arctanh{a}, \quad |a| < 1 \right\} \,, & \text{ if } \XX > 0 \\
    \left\{\TT = -a  \XX, \quad |a| < 1 \right\} \, = \, \left\{t = 4m \, \arctanh{(-a)}, \quad |a| < 1 \right\} \,, & \text{ if } \XX < 0
\end{align*} 
representing time slices of the static time. The time function having \eqref{level sets_tildeT} as time slices turns out to be not only a Cauchy function but also to satisfy the desired  conditions \eqref{cond: Bartnik_RiemNorm}, as it is established in the following.

\vspace{0.1cm}

\begin{proposition}\label{prop: NewTime}
    The function
    \begin{equation}\label{def: newTime}
        \widetilde{\TT} \coloneqq 4m \arctanh{\left( \frac{\TT}{\sqrt{1 + \XX^2}}\right)} \; \in \left( -\infty, + \infty \right)
    \end{equation}
    is a well-defined time function in the maximally extended Schwarzschild $\left( \LL_m, \mathfrak{g}_m\right)$. In particular, it is a Cauchy time function. 
    
    Moreover, the \emph{time-dependent Bartnik conditions} are satisfied with respect to $\widetilde{\TT}$ in the region $\{ r > \kappa_0m\}$ for any $\kappa_0 > 0$, with a constant $C = C(\kappa_0) > 0$, that is,
    \begin{equation}\label{eq: spacetimeConds_tildeT}
        {}^{\widetilde{\TT}} \| \log \alpha_{\widetilde{\TT}}\| + {}^{\widetilde{\TT}} \| \dd \log \alpha_{\widetilde{\TT}}\| \,m + {}^{\widetilde{\TT}} \| \nabla \widetilde{\mathcal{T}}\| \,m + {}^{\widetilde{\TT}} \| \nabla^2\widetilde{\mathcal{T}}\| \,m^2 \leq C
    \end{equation}
    where $\alpha_{\widetilde{\TT}}$ is the lapse function of $\widetilde{\TT}$  
    and $\widetilde{\mathcal{T}}$ is the future-pointing unit normal to the $\widetilde{\TT}$-time slices.
\end{proposition}

\vspace{0.1cm}

\begin{proof}
    By definition, $\widetilde{\TT}$ is a smooth function of the previous coordinates $(\TT, \XX)$ and is well defined in the maximally extended Schwarzschild spacetime $(\LL_m, \mathfrak{g}_m)$, since it holds $\TT^2 < 1 + \XX^2$ and the function $\arctanh$ is smooth on $(-1, 1)$. It is a time function, as its differential reads
    \begin{align}\label{differential T tilde}
        \dd \widetilde{\TT} = \frac{4m}{1 - \frac{\TT^2}{1 + \XX^2}} \frac{1}{\sqrt{1 + \XX^2}} \left( \dd \TT - \frac{\XX \TT}{1 + \XX^2} \dd \XX\right)
    \end{align}
    and consequently 
    \begin{align}\label{eq: tildeT_timeFct} \nonumber
        \mathfrak{g}_m\left(\dd \widetilde{\TT}, \dd \widetilde{\TT} \right) & = 16m^2 \,\xi_m^{-1} \frac{1 + \XX^2}{(1 + \XX^2 - \TT^2)^2} \left(-1 + \frac{\XX^2 \TT^2}{(1 + \XX^2)^2} \right) \\ 
        &= - \frac{16m^2 \,\xi_m^{-1} }{1 + \XX^2 - \TT^2} \left(  1 + \frac{\TT^2}{ (1 + \XX^2 - \TT^2)(1 + \XX^2)} \right) < 0
    \end{align}
    where we used that 
        $(1 + \XX^2 - \TT^2)(1 + \XX^2) = (1 + \XX^2)^2 - \XX^2 \TT^2 - \TT^2$,
    combined with the positivity of the quantities $\xi_m > 0$ and $1 + \XX^2 - \TT^2 > 0$.

    The time function $\widetilde{\TT}$ is Cauchy as we can verify that each time slice $\mathcal{S}_{\widetilde{\TT}}$ is a Cauchy surface.  Let  $s \mapsto \gamma(s) \in \LL_m$ be an inextendible (differentiable) causal curve and, without loss of generality, let us assume that it is future-directed with respect to the time-orientation induced by $\partial_{\TT}$, namely at each point $\gamma(s)$ it holds
    \begin{equation}\label{eq: features gamma}
    \begin{split}
        \xi_m \left[ - \dd \TT(\dot{\gamma}(s))^2 + \dd \XX(\dot{\gamma}(s))^2\right] + r^2 g_{\nS^2}\left(\dot{\gamma}(s), \dot{\gamma}(s)\right)   = \mathfrak{g}\left(\dot{\gamma}(s), \dot{\gamma}(s)\right) \leq 0\,, \\
         - \dd \TT ( \dot{\gamma}(s) ) = \xi_m^{-1} \,\mathfrak{g}\left(\partial_{\TT}, \dot{\gamma}(s)\right) < 0\,.
    \end{split}
    \end{equation}
    We first show that the time function $\widetilde{\TT}$ is strictly increasing along the curve and hence, if $\gamma$ intersects a time slice $\mathcal{S}_{\widetilde{\TT}}$, then this happens only once. At some point $\gamma(s)$, we compute
    \begin{align*}
        \frac{d}{ds} (\widetilde{\TT} \circ \gamma(s))  = \dd \widetilde{\TT}(\dot{\gamma}(s)) & = \frac{4m}{1 - \frac{\TT^2}{1 + \XX^2}} \frac{1}{\sqrt{1 + \XX^2}} \left( \dd \TT(\dot{\gamma}(s)) - \frac{\XX \TT}{1 + \XX^2} \dd \XX(\dot{\gamma}(s))\right) \\
        & > \frac{4m}{1 - \frac{\TT^2}{1 + \XX^2}} \frac{1}{\sqrt{1 + \XX^2}} \left( \dd \TT(\dot{\gamma}(s)) - \left|\dd \XX(\dot{\gamma}(s))\right|\right) \, \geq \, 0\,,
    \end{align*}
    where we used \eqref{differential T tilde} and the properties \eqref{eq: features gamma}, combined with the positivity of the quantities $\xi_m > 0$ , $(1+\XX^2)^2 - \XX^2 > 0$ and $1 + \XX^2 - \TT^2 > 0$. Now, let us prove that $\gamma$ intersects any time slice $\mathcal{S}_{\widetilde{\TT}}^a = \left\{ \widetilde{\TT} = 4m \arctanh(a) \right\} = \left\{ \TT = a \sqrt{1 + \XX^2}\right\}$ for some constant $0 < a < 1$, which, without loss of generality, we assume to be positive. Since it is well-known that the hypersurface $\{ \TT = 0\}$ is Cauchy, we can consider its unique intersection $(0, \XX_0, \theta_0, \varphi_0) = \gamma(0)$ with $\gamma$, which we set at $s = 0$, for convenience. Then, for any $s \geq 0$, the curve $\gamma(s)$ must be contained in the region $\left\{ \XX - \XX_0 \leq \TT\, , \; -(\XX - \XX_0)\leq \TT \right\}$
    between the null hypersurfaces passing through $\gamma(0)$. Now, if by contradiction $\gamma(s),\, s \geq 0$ does not intersect $\mathcal{S}_{\widetilde{\TT}}^a$, then it must be fully contained in the region $\left\{ \XX - \XX_0 \leq \TT\, , \; -(\XX - \XX_0)\leq \TT \,,\; \TT \leq a \sqrt{1 + \XX^2}\right\}$, which is compact. In particular, there exists a finite $ \widetilde{\TT}_{\mathrm{m}} = \max \widetilde{\TT}(\gamma(s)) 
    \in \left[ 0, 4m \arctanh(a) \right]$. Considering the unique point of intersection between $\gamma$ with the time slice $\{ \widetilde{\TT} = \widetilde{\TT}_{\mathrm{m}}\} $ and the correspondent future-directed causal tangent vector, the curve can be extended. This contradicts the inextendibility of $\gamma$.

    Finally, to prove the validity of the \emph{time-dependent Bartnik conditions}, let us compute the lapse $\alpha_{\widetilde{\TT}}$, the future-pointing unit normal $\widetilde{\mathcal{T}}$ and the associated orthonormal frame. From the computation \eqref{eq: tildeT_timeFct}, we directly obtain
    \begin{align} \label{def: lapse_tildeT}
        \alpha_{\widetilde{\TT}}^2 &= -  \mathfrak{g}_m\left(\dd \widetilde{\TT}, \dd \widetilde{\TT} \right)^{-1} = \frac{\xi_m}{16m^2} \;\frac{1 + \XX^2 - \TT^2}{  1 + f(\XX, \TT) }\,, \\ \label{def: f_tildeT}
         f(\XX, \TT) & \coloneqq \frac{\TT^2}{ (1 + \XX^2 - \TT^2)(1 + \XX^2)} \geq 0\,.
    \end{align} 
    Now, we can compute 
\begin{align*}
   -\alpha_{\widetilde{\TT}} \,\dd \widetilde{\TT} &= -\alpha_{\widetilde{\TT}} \, \frac{4m \, \sqrt{1+ \XX^2}}{1 + \XX^2 - \TT^2}\left(\dd \TT - \frac{\XX \TT}{1 + \XX^2} \dd \XX \right) \\
   & = - \sqrt{\frac{\xi_m}{1 + f}} \,\frac{1}{\sqrt{1 + \XX^2 - \TT^2}} \left( \sqrt{1+ \XX^2} \, \dd \TT - \frac{\XX \TT}{\sqrt{1 + \XX^2}} \dd \XX\right) 
\end{align*} 
and this provides the unitary future-directed normal to the time slices
\begin{equation}\label{def: normal_tildeT}
    \widetilde{\mathcal{T}} =  \frac{1}{\sqrt{\xi_m (1 + f)}} \,\frac{1}{\sqrt{1 + \XX^2 - \TT^2}} \left( \sqrt{1+ \XX^2} \, \partial_{\TT} + \frac{\XX \TT}{\sqrt{1 + \XX^2}} \partial_{\XX}\right)\,,
\end{equation}
for $f$ defined in \eqref{def: f_tildeT}. 
To find the other elements of the adapted frame, we consider a time slice $\{\TT = a \sqrt{1 + \XX^2}\,, \;|a | < 1\}$ and its tangent vectors given by $\partial_{\XX} + \tfrac{a \XX}{\sqrt{1 + \XX^2}} \partial_{\TT} = \partial_{\XX} + \tfrac{\TT \XX}{1 + \XX^2} \partial_{\TT}, \partial_{I}$. Computing the norm
\begin{align*}
    \mathfrak{g}_m\left( \partial_{\XX} + \frac{ \XX \TT}{1 + \XX^2} \partial_{\TT}, \, \partial_{\XX} + \frac{\XX \TT}{1 + \XX^2} \partial_{\TT}\right) 
    &= \xi_m \, \frac{1 + \XX^2 - \TT^2}{1 + \XX^2} \,(1 + f)\,,
\end{align*}
    we obtain an adapted frame $\{ \widetilde{\mathcal{T}}, \widetilde{V}_{\XX}, \widetilde{V}_{I} \}$, with $\widetilde{\mathcal{T}}$ as in \eqref{def: normal_tildeT} and 
\begin{align} \label{def: VX_tildeT}
    \widetilde{V}_{\XX} &=  \frac{1}{\sqrt{\xi_m (1 + f)}} \,\frac{1}{\sqrt{1 + \XX^2 - \TT^2}} \left(  \frac{\XX \TT}{\sqrt{1 + \XX^2}} \partial_{\TT} + \sqrt{1+ \XX^2} \, \partial_{\XX}\right)\,, \\ \label{def: VI_tildeT}
    \widetilde{V}_{I} &= \frac{\sqrt{g_{\nS^{2}}^{II}}}{r} \, \partial_{I} \,, \qquad \qquad  I \in \{ \theta, \varphi \}\,.
\end{align}
The existence of a constant $C = C(\kappa_0) > 0$ satisfying \eqref{eq: spacetimeConds_tildeT} in the region $\{ r > \kappa_0 m\}$ can be proven by lengthy explicit computations with the adapted frame described above.
Such computations are partially displayed  in \cref{Appendix: Bartink conds computations} and fully carried out in the doctoral thesis of the author \cite[Appendix A.1]{CogoPhD}.
\end{proof}

\vspace{0.1cm}

\subsection{Bartnik interior condition via monotonicity}

In the following we show that a relaxed version of the \emph{Bartnik interior condition} \eqref{cond: Bartink_Interior} is satisfied with respect to the time function $\widetilde{\TT}$ of \cref{prop: NewTime}. This is achieved by first showing that, given a fixed $\kappa_0 > 0$ and a point $\p$ lying on a large radius, the quantity we need to bound in $\{ r > \kappa_0 m\}$ has a maximum for points at the intersection of $\{ r = \kappa_0 m\}$ and the null hypersurface passing through $\p$; notice that $\{ r = \kappa_0 m\}$, for $0 < \kappa_0 < \tfrac{3}{2}$, is a barrier near the crushing singularity, as in \cref{prop: chrushingSing_barriers}).  In turn, such quantity on the intersection points can be expressed as a function of only $\TT(\p)$ and we are able to show that this decreases monotonically as $\TT(\p) \to \infty$.

\vspace{0.1cm}

\begin{proposition}\label{prop: monotonicity}
    For any radius $\kappa_0m \leq 2m$ and any radius $Km > 2m$,  
\    there exists a constant $C = C({\kappa_0, K}) > 0$
such that for any $\p \in \{ r = Km\} \subset \LL_m$, it holds
    \begin{equation}\label{difference_Ttilde_IntCond}
   \begin{split}
    \sup_{\{\kappa_0 m < r < Km\} \setminus I^+_{\p}} \;\widetilde{\TT} - \widetilde{\TT}(\p) \; &\leq\, C\,m \qquad \,\text{ if } \widetilde{\TT}(\p) \geq 0\,, \\
    \sup_{\{\kappa_0 m < r < Km\} \setminus I^-_{\p}} \;  \widetilde{\TT}(\p) - \widetilde{\TT} \; &\leq\, C \,m \qquad\, \text{ if } \widetilde{\TT}(\p) \leq 0\,.
\end{split}
\end{equation}
where $\widetilde{\TT}$ is the time function defined in \cref{prop: NewTime} and where we denoted 
\begin{equation}\label{def: I_p}
    I^+_{\p} \coloneqq \bigcup_{\overline{\p} \in \{ r = Km\} \cap \left\{ \widetilde{\TT} = \widetilde{\TT}(\p)\right\}} I^+(\overline{\p}) \,, \qquad I^-_{\p} \coloneqq \bigcup_{\overline{\p} \in \{ r = Km\} \cap \left\{ \widetilde{\TT} = \widetilde{\TT}(\p)\right\}} I^-(\overline{\p}) \,,
\end{equation}
for $I^{\pm}(\overline{\p})$  \emph{causal past and future} of the point $ \overline{\p}$.
\end{proposition}

\begin{proof}
Consider a point $\p \in \{ r = Km\}$, which can be taken with $\XX(\p) > 0$. For a fixed $0 <\kappa_0 \leq 2$ and $K > 2$, which without loss of generality, satisfies $K > 2\big(1 + \mathrm e^{-\frac{K}{2}}\big)$, we set the constants
$$ \varrho_{\kappa_0} \coloneqq \varrho_m(\kappa_0m)\; \in (-1, 0]\,, \qquad \quad    \varrho_{K} \coloneqq 
\varrho_m(Km) \; \in (1, +\infty)\,,$$
for $\varrho_m$ as defined in \eqref{def: changeCoord_Kruskal_r}. Now, let $\TT_{\p} \coloneqq \TT(\p)$, then  \eqref{def: changeCoord_Kruskal_r} yields
\begin{equation}\label{eq: X_p}
    \XX_{\p} \coloneqq\XX(\p) = \sqrt{\TT_{\p}^2 + \varrho_{K} }\,
\end{equation} 
and, by the definition of $\widetilde{\TT}$ in \eqref{def: newTime}, we get 
\begin{equation}\label{eq: tildeTT_p} 
    \widetilde{\TT}(\p) = 4m \arctanh \left( \frac{\TT_{\p} }{\sqrt{1 + \TT_{\p}^2 +  \varrho_{K}}} \right)\,.
\end{equation}
We can assume $\TT_{\p} \geq 0$ and hence $ \widetilde{\TT}(\p) \geq 0$, as the negative case can just be discussed using the symmetries of the Kruskal extension. In order to detect the points involved in the estimate \eqref{difference_Ttilde_IntCond}, we notice that we can just consider  the region $\big\{ \widetilde{\TT} \geq  \widetilde{\TT}(\p)\big\}$ and 
we identify the set $\{ r = Km\} \cap \big\{  \widetilde{\TT} = \widetilde{\TT}(\p) \big\} = \big\{ \TT = \TT_{\p}, \, \XX = \pm \XX_{\p} \big\}$. This determines $I^+_{\p}$ as defined in \eqref{def: I_p}, by considering the null hypersurfaces $\left\{\TT = (\TT_{\p} + \XX_{\p}) \pm \XX \right\}$ passing through such points. Hence, up to restricting to the region $\big\{ \widetilde{\TT} \geq  \widetilde{\TT}(\p)\big\}$, we have
\begin{equation}\label{region_sup}
    \{\kappa_0 m < r < Km\} \setminus I^+_{\p} = \left\{ \TT \leq \big(\TT_{\p} + \XX_{\p}\big) - \XX,  \TT \leq \big(\TT_{\p} + \XX_{\p}\big) + \XX , \TT^2 < \XX^2 - \varrho_{\kappa_0}\right\} 
\end{equation}
and we claim that the supremum of $\widetilde{\TT}$ (and consequently of $\widetilde{\TT} - \widetilde{\TT}_{\p}$) in this region is attained at the points
\begin{equation}\label{points q}
     \q \in \left\{\TT = (\TT_{\p} + \XX_{\p}) - \XX \right\} \cap \left\{  \TT^2 = \XX^2 - \varrho_{\kappa_0} \right\} \,.
\end{equation} 
Notice that by symmetry, we can show the above claim just on $\{ \XX \geq 0\}$. In the region \eqref{region_sup},
considering the points in $\big\{ \XX \geq \XX(\q)\big\}$, there holds $\TT \leq (\TT_{\p} + \XX_{\p}) - \XX \leq (\TT_{\p} + \XX_{\p}) - \XX(\q) = \TT(\q)$; hence, due to the monotonicity of $\arctanh$, we have $\widetilde{\TT} \leq 4m\arctanh\left( \frac{\TT(\q)}{1 + \XX(\q)}\right) = \widetilde{\TT}(\q)$. On the other hand, regarding the points in $\big\{ 0 \leq  \XX \leq \XX(\q)\big\}$, it follows by monotonicity that $\frac{\TT}{\sqrt{1 + \XX^2}} \leq \frac{\sqrt{\XX^2 - \varrho_0}}{\sqrt{1 + \XX^2}} \leq \frac{\sqrt{\XX^2(\q) - \varrho_0}}{\sqrt{1 + \XX^2(\q)}} = \frac{\TT(\q)}{\sqrt{1 + \XX^2(\q)}}$ and consequently, $\widetilde{\TT} \leq \widetilde{\TT}(\q)$. 

Thus far, we have shown that 
\begin{equation}\label{BartnikInterior_TpTq}
    \sup_{\{\kappa_0 m < r < Km\} \setminus I^+_{\p}} \;\widetilde{\TT} - \widetilde{\TT}(\p)  \, = \, \widetilde{\TT}(\q) - \widetilde{\TT}(\p)
\end{equation}
for $\q$ as in \eqref{points q}. Let us now prove that this quantity can be bounded by a constant independent of $\p$ and $\q$. We first compute
\begin{equation}\label{eq: X_q T_q}
    \XX(\q) = \frac{\TT_{\p} + \XX_{\p}}{2} + \frac{\varrho_{\kappa_0}}{2\left( \TT_{\p} + \XX_{\p} \right)}\,, \quad \TT(\q) = (\TT_{\p} + \XX_{\p}) - \XX(\q) = \frac{\TT_{\p} + \XX_{\p}}{2} - \frac{\varrho_{\kappa_0}}{2\left( \TT_{\p} + \XX_{\p} \right)}\,.
\end{equation}
and in turn, by \eqref{def: newTime} and \eqref{eq: X_p}, 
we obtain
\begin{align*}
    \widetilde{\TT}(\q) & = 4m \arctanh \left( \frac{\TT(\q)}{\sqrt{1 + \XX(\q)^2}} \right) 
    = 4m \arctanh \left( \frac{\TT_{\p}  + \sqrt{\TT_{\p}^2 + \varrho_{K}} - \frac{\varrho_{\kappa_0}}{\TT_{\p} + \sqrt{\TT_{\p}^2 + \varrho_{K}}}}{\sqrt{4 +  \left(\TT_{\p} +\sqrt{\TT_{\p}^2 + \varrho_{K}} + \frac{\varrho_{\kappa_0}}{\TT_{\p} + \sqrt{\TT_{\p}^2 + \varrho_{K}}}\right)^2}} \right).
\end{align*}
Notice that both $\widetilde{\TT}({\p})$ as in \eqref{eq: tildeTT_p} and $\widetilde{\TT}({\q})$ are just functions of $\TT_{\p}$ for fixed $\kappa_0, K$, corresponding to fixed constants $\varrho_{\kappa_0}$ and $\varrho_{K}$. Hence, setting $x \coloneqq \TT_{\p}$, we write the quantity we want to estimate as a function $h_{\kappa_0, K}(x) \coloneqq \frac{\widetilde{\TT}({\q}) - \widetilde{\TT}({\p})}{4m}$, which reads
\begin{align*}
    h_{\kappa_0, K}(x) = & \arctanh \left( \frac{x  + \sqrt{x^2 + \varrho_{K}} - \frac{\varrho_{\kappa_0}}{x  + \sqrt{x^2 + \varrho_{K}}}}{\sqrt{4 +  \left(x +\sqrt{x^2 + \varrho_{K}} + \frac{\varrho_{\kappa_0}}{x  + \sqrt{x^2 + \varrho_{K}}}\right)^2}} \right) 
    - \arctanh \left( \frac{x }{\sqrt{1 + x^2 + \varrho_{K}}} \right).
\end{align*}
and let us show that it is monotonically decreasing. By explicit computations, fully carried out in \cite[Proof of Proposition 2.5.3]{CogoPhD}, one finds
\begin{equation}
    \frac{d}{dx}  h_{\kappa_0, K}(x)  = \frac{\left(1 + \frac{x}{\sqrt{x^2 + \varrho_{K}}}\right) \left( 1 +  \frac{\varrho_{\kappa_0}}{\left( x + \sqrt{x^2 + \varrho_{K}}\right)^2}\right)}{\sqrt{4 +  \left(x +\sqrt{x^2 + \varrho_{K}} + \frac{\varrho_{\kappa_0}}{x  + \sqrt{x^2 + \varrho_{K}}}\right)^2}} - \frac{1}{\sqrt{1 + x^2 + \varrho_{K}}} \;\leq 0 \,.
\end{equation}
Thus, we can estimate the desired quantity \eqref{BartnikInterior_TpTq} as follows
\small
\begin{align*}
    &\frac{\widetilde{\TT}({\q}) - \widetilde{\TT}({\p})}{4m}
     = \arctanh \left( \frac{\TT_{\p}  + \sqrt{\TT_{\p}^2 + \varrho_{K}} - \frac{\varrho_{\kappa_0}}{\TT_{\p}  + \sqrt{\TT_{\p}^2 + \varrho_{K}}}}{\sqrt{4 +  \left(\TT_{\p} +\sqrt{\TT_{\p}^2 + \varrho_{K}} + \frac{\varrho_{\kappa_0}}{\TT_{\p}  + \sqrt{\TT_{\p}^2 + \varrho_{K}}}\right)^2}} \right) -\arctanh \left( \frac{\TT_{\p} }{\sqrt{1 + \TT_{\p}^2 + \varrho_{K}}} \right) \\
   & \qquad \leq \lim\displaylimits_{\TT_{\p} \to 0}  \left(  \arctanh \left( \frac{\TT_{\p}  + \sqrt{\TT_{\p}^2 + \varrho_{K}} - \frac{\varrho_{\kappa_0}}{\TT_{\p}  + \sqrt{\TT_{\p}^2 + \varrho_{K}}}}{\sqrt{4 +  \left(\TT_{\p} +\sqrt{\TT_{\p}^2 + \varrho_{K}} + \frac{\varrho_{\kappa_0}}{\TT_{\p}  + \sqrt{\TT_{\p}^2 + \varrho_{K}}}\right)^2}} \right) - \arctanh \left( \frac{\TT_{\p} }{\sqrt{1 + \TT_{\p}^2 + \varrho_{K}}} \right) \right) \\
   & \qquad 
   = \arctanh\left( \sqrt{\frac{\varrho_{K} \left(1- \frac{\varrho_{\kappa_0}}{\varrho_{K}} \right)^2}{4 + \varrho_{K} \left(1 + \frac{\varrho_{\kappa_0}}{\varrho_{K}} \right)^2}}\right)\,.
\end{align*}
\normalsize
Now, defining $C = C(\kappa_0, K) \coloneqq 4 \arctanh\left( \sqrt{\frac{\varrho_{K} \left(1- \frac{\varrho_{\kappa_0}}{\varrho_{K}} \right)^2}{4 + \varrho_{K} \left(1 + \frac{\varrho_{\kappa_0}}{\varrho_{K}} \right)^2}}\right)$, we conclude.

\end{proof}

\subsection{Gluing of time functions}\label{sec: gluing}
In the following, we will glue on large annuli in the asymptotically flat ends the time function $\widetilde{\TT}$ introduced in \cref{prop: NewTime} with $t -u_{\bb}$, where $t$ is the static time and $u_{\bb}$ is as in \cref{prop: barriers_infty}.

\vspace{0.2cm}
\noindent
\textbf{Gluing $\widetilde{\TT}$ to the static time.}
We start by gluing the new time function $\widetilde{\TT}$ defined in \eqref{def: newTime} with the static time $t$, by showing that the corresponding time slices can be glued into spacelike hypersurfaces. Consider now a constant $ c \in [ 0, 1) $ and the level set $\{ \widetilde{\TT} = 4m \arctanh c \} \subset \left\{ \TT \geq 0\right\}$. The case $ c \in (-1, 0]$ and the corresponding level sets $\{ \widetilde{\TT} = 4m \arctanh c \} \subset \left\{ \TT \leq 0\right\}$ follows by symmetry. We then perform the appropriate change of coordinates in the \emph{region} $\mathrm{I}  = \{ \TT^2 < \XX^2, \, \XX > 0 \}$ and \emph{region} $\mathrm{I'} = \{ \TT^2 < \XX^2, \, \XX < 0 \}$, where in particular the Schwarzschildean radius satisfies $r > 2m$, as follows
\begin{equation}\label{level set region I}
    \left\{ \widetilde{\TT} = 4m \arctanh c  \right\} \cap \{ \TT^2 < \XX^2, \, \XX > 0 \} = \left\{  t =  4m \arccosh \sqrt{\frac{1}{1-c^2} \left( 1 + c^2 \varrho_m^{-1}(r) \right) }\right\}\,
\end{equation}
\begin{equation}\label{level set region I'}
    \left\{ \widetilde{\TT} = 4m \arctanh c  \right\} \cap \{ \TT^2 < \XX^2, \, \XX < 0 \} = \left\{  t =  - 4m \arccosh \sqrt{\frac{1}{1-c^2} \left( 1 + c^2 \varrho_m^{-1}(r) \right) }\right\}
\end{equation}
where we used the definition of $\widetilde{\TT}$ in \eqref{def: newTime}, which yields $\left\{ \widetilde{\TT} = 4m \arctanh c  \right\} = \{ \TT = c \sqrt{1 + \XX^2} \}$,  combined with \eqref{def: changeCoord_Kruskal_r} and \eqref{def: changeCoord_Kruskal_rt_TTXX}, which imply that on the points of the level sets there holds
\[ \varrho_m(r) =  \XX^2 (1-c^2) - c^2 =   (1-c^2) \varrho_m(r) \cosh^2 \left( \frac{t}{4m}\right) - c^2. \]
We will then glue the hypersurface \eqref{level set region I} and  \eqref{level set region I'} respectively to the level set of the static time 
\begin{equation}\label{level set static}
     \left\{ t = 4m \arctanh c  \right\} = \left\{ t = 4m \arccosh \sqrt{\frac{1}{1 - c^2}}  \right\}
\end{equation}
\begin{equation*}
      \left\{ t = -4m \arctanh c  \right\} = \left\{ t = -4m \arccosh \sqrt{\frac{1}{1 - c^2}}  \right\}.
\end{equation*}
Since the two cases differ just by the sign, we will deal only with \emph{region} $\mathrm{I}$.
Consider the function in Schwarzschildean coordinates
\[ \mathrm F_1(r) \coloneqq \varsigma_1(r) \left( 4m \arccosh{ \sqrt{\frac{1}{1-c^2}}}\right) + (1-\varsigma_1(r)) \left( 4m  \arccosh{\sqrt{\frac{ 1 + c^2 \,\varrho_m^{-1}(r) }{1-c^2} }} \right)\,\]
where  $\mathcal{\varsigma}_1: (0, +\infty) \to \nR$ is a cut-off function
\begin{equation}\label{def: cutoff1}
\varsigma_1(s) = 
    \begin{cases}
        1 & \quad \text{ if } s < 16 m   \\
        0 & \quad \text{ if } s > 32 m 
    \end{cases}
    \,,\qquad \qquad  \vert \varsigma_1'(s) \vert < \frac{1}{8 m } \,.
\end{equation}
Since  $\{ t = \mathrm F_1(r) \} $ coincides with the level set \eqref{level set region I} in $\{ r \leq 16m\}$ and with the level set  \eqref{level set static} in $\{ r \geq 32m \}$, we just need to check that it is a spacelike hypersurface in $\{ 16 m < r < 32m \}$.  To check the spacelike condition \eqref{spacelikeCond_Schw}, we  compute
    \begin{align}\label{estimate_gluing1} \nonumber
    \vert h_m \partial_r \mathrm F_1 \vert 
    & =  4m \, h_m \,\vert\varsigma_1'(r) \vert \log\left( \frac{\sqrt{1 + c^2 \varrho_m^{-1}(r)} + c \sqrt{1 + \varrho_m^{-1}(r)}}{1 + c} \right) \\ 
    &  \qquad + \left( 1 - \varsigma_1(r) \right) \frac{c \, \varrho_m^{-1}(r) }{ \sqrt{1 + \varrho_m^{-1}(r)} \sqrt{1 + c^2 \varrho_m^{-1}(r)}} 
    \;\leq \frac{15}{16}  \sqrt{1 + \frac{e^{-8}}{7}} + \frac{e^{-8}}{7} \; < 1
\end{align}
using the identities $ \varrho_m'(r) 
= \frac{r}{4m^2} e^{\frac{r}{2m}} >0$ and $\varrho_m^{-1}(r) \varrho_m'(r) = \frac{1}{2m \, h_m}$.
This proofs the claim.

\vspace{0.5cm}
\noindent
\textbf{Gluing the static time slices to the coordinate boosts admitting barriers.} In the following, it is convenient to work with the isotropic coordinates and, for simplicity, in \emph{region} $\mathrm{I}$. By changing the radial coordinate $r = 32m$, found in the previous step, via \eqref{changeCoord_Iso_Schw} into $\vert \xx \vert = \frac{31 + 8 \sqrt{15}}{2}m \ll 48 m$, we consider the following Lipschitz compactly supported 
function 
\begin{equation*}
    F_{\bb} (\xx) = 
    \begin{cases}
        0 & \quad \text{ if } \; \vert \xx \vert \leq 48m \\
         \frac{\bb + 1}{2} \left( 1-   \frac{48 m}{\vert \xx \vert} \right) x_1 & \quad \text{ if } \; 48m \leq \vert \xx \vert \leq 48m \left( 1 + \frac{2\bb}{1 - \bb}\right) \\
                \bb x_1 & \quad \text{ if } \;  48 m \left( 1 + \frac{2\bb}{1 - \bb}\right) \leq  \vert \xx \vert \leq 192 m \left( 1 + \frac{2\bb}{1 - \bb}\right)
            \\ \frac{\bb + 1}{2} \left(  \left( 1 + \frac{4\bb}{1 - \bb}\right) \frac{192 m}{\vert \xx \vert} -1  \right) x_1 & \quad \text{ if } \;  192 m \left( 1 + \frac{2\bb}{1 - \bb}\right) \leq  \vert \xx \vert \leq 192 m \left( 1 + \frac{4\bb}{1 - \bb}\right)
            \\ 0 & \quad \text{ if } \;   \vert \xx \vert \geq 192 m \left( 1 + \frac{4\bb}{1 - \bb}\right)
    \end{cases}
\end{equation*}
 which has Lipschitz constant $\operatorname{Lip}(F_{\bb}) = \tfrac{\bb + 1 }{2} < 1$. Consider now for sufficiently small parameters $0 < \delta < 1$ a sequence of mollified functions $F_{\bb}^{\delta} \in C^{\infty}_0 (\nR^3)$ with compact support on the region 
$\left\{ 48m - 2\delta <  \vert \xx \vert < 192m \left(1 + \tfrac{4\bb}{1 - \bb} \right) + 2 \delta\right\}$ and satisfying $\operatorname{Lip}(F^{\delta}_{\bb}) = \tfrac{1 + \bb }{2}$. Finally, considering radial cut-off functions in the spirit of \eqref{def: cutoff1}, we can glue $\{ t = 0 \}$ to $\{ t = F^{\delta}_{\bb}(\xx)\}$ in the region $\{ 36m < \vert \xx \vert < 40m\}$ and $\{ t = F^{\delta}_{\bb}(\xx)\}$ to $\{ t = \bb x_1\}$ in the region $\left\{ 96m \left( 1 + \frac{2 \bb}{1 - \bb} \right) < \vert \xx \vert < 144m \left( 1 + \frac{2 \bb}{1 - \bb} \right) \right\}$. Since in both these regions the function $F_{\bb}$ is smooth, the sequence of mollified functions and their gradients uniformly converge as $\delta \to 0$; this ensures that for sufficiently small $\delta = \delta(\bb) > 0 $, the resulting glued hypersurface satisfies the spacelike condition. Since this can be done replacing the zero level set with any constant, we have glued the static time to the time function $t - \bb x_1$.

Finally, let us consider the function $u_{\bb}$ and $K_{\bb} = K_{\bb}(\bb, \varepsilon, C_{f_{\bb}} ) \gg 1$ as in \cref{prop: barriers_infty}. We can now assume, without loss of generality that $\frac{K_{\bb}}{4} > 144 \left( 1 + \frac{2 \bb}{1 - \bb} \right)$ and consider a cut-off function $\mathcal{\varsigma}_2: (0, +\infty) \to \nR$ is taken so that
\begin{equation*}\label{def: cutoff3}
\varsigma_3(s) = 
    \begin{cases}
        1 & \quad \text{ if } s < \frac{K_{\bb}}{2} m   \\
        0 & \quad \text{ if } s > K_{\bb} m 
    \end{cases}
    \,,\qquad \qquad  \vert \varsigma_3'(s) \vert < \frac{4}{K_{\bb} m }\,.
\end{equation*}
Up to taking $K_{\bb}$ big enough still depending on $\bb$ and $C_{f_{\bb}}$, it is easy to see that the graph of the function
\[ F_2(\xx) = \varsigma_2(\vert \xx \vert) \bb x_1 + (1 - \varsigma_2(\vert \xx \vert)) u_{\bb}(\xx) \]
is spacelike. 
In conclusion, we have glued the static time to the time function $t - u_{\bb}$ with $u_{\bb}$ as in \cref{prop: barriers_infty}.

\vspace{0.1cm}

\subsection{Boosted
time function satisfying Bartnik 
conditions and admitting barriers}

In the following theorem, we claim the existence of a time function satisfying both the \emph{time-dependent Bartnik conditions} \eqref{cond: Bartnik_RiemNorm} and a relaxed version of the \emph{Bartnik interior condition}. Moreover, we ensure that its time slices belong to the class $\mathcal{B}_{\bb}$ and admit lower and upper barriers for large radii. Such function is obtained by the gluing of $t - u_{\bb}(\xx)$ and $\widetilde{\TT}$, as in \cref{sec: gluing}.

\vspace{0.1cm}

\begin{theorem}\label{thm: global_time_function}
    For any $\bb \in [0, 1)$, the maximally extended Schwarzschild spacetime admits a time function $\ttt_{\bb} \in C^{\infty}(\mathfrak{L}_m)$ such that the following hold simultaneously.
    \begin{itemize}
        \item [(i)]  The \emph{Bartnik's spacetime conditions} \eqref{cond: Bartnik_RiemNorm} with norm \eqref{def: timeDep_RiemNorm} defined with respect to $\ttt_{\bb}$  are satisfied for some constant $C_o = C_o(\kappa_0) > 0$ and $\kappa_0 > 0$ on  a region $ \{ r > \kappa_0m\}$:
        \begin{equation}\label{eq: spacetimeConds_ttt_b}
        {}^{\ttt_{\bb}} \| \log \alpha_{\ttt_{\bb}}\| + {}^{\ttt_{\bb}} \| \dd \log \alpha_{\ttt_{\bb}}\| m +  {}^{\ttt_{\bb}} \| \nabla \mathcal{T}_{\bb}\| m + {}^{\ttt_{\bb}} \| \nabla^2 \mathcal{T}_{\bb}\| m^2 \leq C_o\,,
    \end{equation}
    where $\alpha_{\ttt_{\bb}}$ is the lapse function of $\ttt_{\bb}$ as defined in \eqref{def: lapse} and $\mathcal{T}_{\bb}$ is the future-pointing unit normal to the time slices $\mathcal{S}_{\bb}$.
        \item [(ii)] Each time slice is a coordinate boost of parameter $\bb \in [0, 1)$ as in \cref{def: coordinate_boost}.
         Moreover, the spacelike hypersurface $\mathcal{S}_{\bb} = \{ \ttt_{\bb} = 0\} \in \mathcal{B}_{\bb}$, given a fixed $0 < \varepsilon < 1$, admits barriers $\mathcal{W}_{\bb, K, \varepsilon}^{\pm}$ on $\{ \vert \xx \vert \geq Km\}$ for any $K \geq K_{\bb}$, with $K_{\bb} = K_{\bb}(\bb, \varepsilon, C_{\bar f_{\bb}})$ as in \cref{prop: barriers_infty}. Moreover, is invariant under the isometry $(\XX, \TT) \mapsto (-\XX, -\TT)$. 
        \item [(iii)] A relaxed version of \emph{Bartnik interior condition}
        holds in $ \{ r > \kappa_0m\}$ for any $\kappa_0 > 0$ and any  $K \geq K_{\bb}$, that is there exists a constant $C_{K} = C_{K}(\bb,   C_{\bar f_{\bb}}, \kappa_0, K) > 0$  such that  
        \begin{equation}
    \begin{split}\label{cond: Bartink_Interior_t_b}
   \sup_{\p \in \{ r = \widetilde{K}m\}} \;\; \sup_{\q \in \{\kappa_0m < r < \widetilde{K}m\} \setminus I^+_{\p}} \;\ttt_{\bb}(\q) - \ttt_{\bb}(\p) \; &\leq C_{K} m \qquad \text{ if } \ttt_{\bb}(\p) \geq 0\,, \\
   \sup_{\p \in \{ r = \widetilde{K}m\}} \;\; \sup_{\q \in \{\kappa_0m < r < \widetilde{K}m\} \setminus I^-_{\p}} \;  \ttt_{\bb}(\p) - \ttt_{\bb}(\q) \; &\leq C_{K} m \qquad \text{ if } \ttt_{\bb}(\p) \leq 0\,.
\end{split}
\end{equation}
where $\widetilde{K}m$ is the Schwarzschildean radius corresponding to the isotropic $Km$ via \eqref{changeCoord_Iso_Schw} and 
\begin{equation}\label{def: I_p_general}
    I^+_{\p} \coloneqq \bigcup_{\overline{\p} \in \{ r = \widetilde{K} m\} \cap \left\{ \ttt_{\bb} = \ttt_{\bb}(\p)\right\}} I^+(\overline{\p}) \,, \qquad I^-_{\p} \coloneqq \bigcup_{\overline{\p} \in \{ r = \widetilde{K} m\} \cap \left\{ \ttt_{\bb} = \ttt_{\bb}(\p)\right\}} I^-(\overline{\p}) \,.
\end{equation}
    \end{itemize}
\end{theorem}

\vspace{0.1cm}

\begin{proof}
    Due to the symmetries of the spacetime, in the following we will work in $\{ \XX \geq 0\}$.
    We claim that the desired time function is the one explicitly constructed by a gluing argument in \cref{sec: gluing}, which coincides with $\ttt_{\bb} = \widetilde{\TT}$ in $\{ r < 16m\}$ and with $\ttt_{\bb} = t - u_{\bb}$ in $\{ \vert \xx \vert > K_{\bb}m \}$, where $u_{\bb}$ and $K_{\bb}m \gg 16m$ are as in \cref{prop: barriers_infty}.
    \begin{itemize}
        \item [\emph{(i)}] In the region $\{ r \leq 16 m\}$, the claim directly follows from the property \eqref{eq: spacetimeConds_tildeT} of $\widetilde{\TT}$ proven in \cref{prop: NewTime}. On the other hand, in $\{ r > 16m\}$ such a feature holds due to  \cref{prop: Bartik_spactime_cond_graphs_exterior}: given the explicit gluing of time functions performed in \cref{sec: gluing}, in this region the time slices of $\ttt_{\bb}$ are given by graphs of a smooth function with bounded derivatives of any order, as a consequence of the definition \eqref{S_b u_b admitting barriers} of $u_{\bb}$
        via the function $\bar f_{\bb}$ satisfying \eqref{def: C_f}. 
        \item [\emph{(ii)}] This property is an immediate consequence of \cref{prop: barriers_infty}.
        \item [\emph{(iii)}] This relaxed version of the \emph{Bartnik interior condition} follow from the corresponding feature of $\widetilde{\TT}$ established in \cref{prop: monotonicity} and the fact that the time functions $\widetilde{\TT}$ and $t - u_{\bb}$ are comparable by fixed constants for points $\q \in \{16 m \leq  r \leq  \widetilde{K} m\}$. To perform the comparison, we use the identity
     $\XX(\q) = \sqrt{\TT(\q)^2 + \varrho_m(r(\q)) }$ given by \eqref{def: changeCoord_Kruskal_rt_TTXX} to obtain the following estimate on the difference between $\widetilde{\TT}$ and the static time
     \begin{align*}
         \vert & t(\q) -  \widetilde{\TT}(\q) \vert = 4m \left\vert \arctanh\left( \frac{\TT(\q)}{\XX(\q)}\right) - \arctanh\left( \frac{\TT(\q)}{\sqrt{1 + \XX(\q)^2}}\right)   \right\vert \\
         &
         \leq 4m \lim_{\TT(\q) \to +\infty} \left( \arctanh\left( \frac{\TT(\q)}{\sqrt{\TT(\q)^2 + \varrho_m(r) }}\right)  - \arctanh\left( \frac{\TT(\q)}{\sqrt{1 + \TT(\q)^2 + \varrho_m(r)}}\right) \right) \\
         & = 2m \log \left( 1 + \frac{1}{\varrho_m(r)}\right) \,\leq 2m \log \left( 1 + \frac{1}{7 \mathrm{e}^8}\right)
     \end{align*}
     Moreover, it is easy to see that
     $ \vert t(\q) - (t-u_{\bb})(\q) \vert \leq \left(\bb K + C_{\bar f_{\bb}} \log(K)\right) m$, where $C_{\bar f_{\bb}}$ is defined in \eqref{def: C_f},
     and hence
     \[\vert \ttt_{\bb}(\q) - \widetilde{\TT}(\q) \vert \leq \left[ 2 \log \left( 1 + \frac{1}{7 \mathrm{e}^8}\right) +  \bb K + C_{\bar f_{\bb}} \log(K)\right] m\,. \]
     The condition \eqref{cond: Bartink_Interior_t_b} is then satisfied with the constant 
     \[ C_K( \bb, C_{\bar f_{\bb}}, \kappa_0, K) \coloneqq C(\kappa_0, K) + 2 \log \left( 1 + \frac{1}{7 \mathrm{e}^8}\right) +  \bb K + C_{\bar f_{\bb}} \log(K)\,,\]
     where $C(\kappa_0, K)$ is the constant bounding the quantities \eqref{difference_Ttilde_IntCond} in \cref{prop: NewTime}.
    \end{itemize}
\end{proof}

As we have mentioned, a fundamental tool to prove the first claim in the above result is given by the following proposition, which shows the validity of the \emph{time-dependent Bartnik conditions} for time functions defined as suitable graphs over the static time in exterior Schwarzschild. 

\vspace{0.1cm}

\begin{proposition}\label{prop: Bartik_spactime_cond_graphs_exterior}
  Consider a spacelike hypersurface $\mathcal{S} = \{ t = \mathcal f(\xx)\}$ in a region $\{ \vert \xx \vert  > \kappa m\}$ with $\kappa > \frac{1}{2}$ of exterior Schwarzschild.
  Assume that $\mathcal f \in C^{\infty} (\nR^3 \setminus \overline{B}_{\kappa m})$ and that there exists a uniform constant $C_{\mathcal f} > 0$ such that
  \begin{equation}\label{bound_hypothesis}
       \sup_{\{ \vert \xx \vert  >  \kappa m\}} \vv_t(\mathcal f) +  \vert \partial^2 \mathcal f \vert \, m +  \vert \partial^3 \mathcal f \vert \, m^2 \leq C_{\mathcal f}\,.
  \end{equation}
  Then the function 
  \[ \tau \coloneqq t - \mathcal f(\xx)\]
  is a time function on $\{ \vert \xx \vert  >  \kappa m\}$ satisfying the conditions \eqref{cond: Bartnik_RiemNorm}, namely there exists a constant $C = C(C_{\mathcal f}, \kappa) > 0$ such that 
   \[ {}^{\tau} \| \log \alpha_{\tau}\| +  {}^{\tau} \| \dd \log \alpha_{\tau}\| m + {}^{\tau} \| \olNabla \mathcal{T}_{\mathcal f}\| m+ {}^{\tau} \| \olNabla^2 \mathcal{T}_{\mathcal f}\| m^2 \leq C\,,\] 
  where $\alpha_{\tau}$ is the lapse function of $\tau$ as defined in \eqref{def: lapse} and $\mathcal{T}_{\mathcal f}$ is the future-pointing unit normal to the time slices of $\tau$.
\end{proposition}

\begin{proof}
We first compute the lapse for the time function $\tau$, which reads
\begin{equation*}
    \alpha_{\tau}^2 = - \G_m(\dd\tau, \dd\tau)^{-1} = \frac{N_m^2}{ 1 - \frac{N_m^2}{\varphi_m^4} \sum_{i = 1}^3 (\partial_i \mathcal f)^2 } = N_m^2 \, \vv_t( \mathcal f)^2 
\end{equation*}
and thus we can estimate
\begin{equation}\label{estimate_alpha}
     \log\left( \frac{1 - \frac{1}{2  \kappa}}{1 + \frac{1}{2  \kappa}} \right) < \log \alpha_{\tau} < \log C_{\mathcal f}\,.
\end{equation}

To estimate the other quantities of \eqref{cond: Bartnik_RiemNorm} in terms of the norm \eqref{def: timeDep_RiemNorm}, we have to find an adapted orthonormal frame. We consider $\{\mathcal{T}_{\mathcal f}, V_1 , V_2, V_3 \} $,
where $\mathcal{T}_{\mathcal f}$ is the unit normal to the time slices and $\{ V_1, V_2, V_3\}$ are tangent to the hypersurface $\mathcal{S} = \{ t = \mathcal f(\xx)\}$. We can compute $-\alpha _{\tau} \dd\tau = 
-N_m \vv_t(\mathcal f) \left( \dd t - \partial_i \mathcal f \, \dd x^i \right)$
and consequently
\begin{align}\label{ONF_t}
    \mathcal{T}_{\mathcal f} = 
    N_m \vv_t(\mathcal f) \left(N_m^{-2} \,\partial_t + \varphi_m^{-4} \sum_{i =1}^3 \partial_i \mathcal f \, \partial_{x^i} \right)
      = \alpha_{\tau} \left(N_m^{-2} \,\partial_t + \varphi_m^{-4} \sum_{i =1}^3 \partial_i \mathcal f \, \partial_{x^i} \right) \,.
\end{align}
On the other hand, the tangential frame can be obtained by applying the Gram–Schmidt process to the frame $\{ \widetilde{V}_1, \widetilde{V}_2, \widetilde{V}_3 \}$, where
\begin{equation}\label{Frame_notON}
    \widetilde{V}_i = 
 \partial_{x^i} + \partial_i \mathcal f \, \partial_t  \,.
\end{equation}
In particular, we can set 
\begin{equation}\label{ONF_1}
    V_1  = \frac{1}{\sqrt{ \varphi_m^4 - N_m^2 (\partial_1 \mathcal f)^2}} \, \widetilde{V}_1 
\end{equation}
and compute 
\begin{align*}
     \mathfrak{g}_m\left( \widetilde{V}_2, V_1 \right)  = \frac{- N_m^2 (\partial_1 \mathcal f)(\partial_2 f) }{\sqrt{ \varphi_m^4 - N_m^2 (\partial_1 \mathcal f)^2}} \, , \qquad     \mathfrak{g}_m\left( \widetilde{V}_3, V_1 \right)  = \frac{- N_m^2 (\partial_1 \mathcal f)(\partial_3 \mathcal f) }{\sqrt{ \varphi_m^4 - N_m^2 (\partial_1 \mathcal f)^2}}\,.
\end{align*}
A tangent vector field to $\mathcal{S}= \{ t = \mathcal f(\xx) \}$ that is orthogonal to $V_1$ is given by $\widetilde{V}_2 -  \mathfrak{g}_m\left( \widetilde{V}_2, V_1 \right) V_1$. By normalising it, we obtain 
\begin{equation}\label{ONF_2}
\begin{split}
     V_2 & = \sqrt{\frac{\varphi_m^4 - N_m^2 (\partial_1 \mathcal f)^2}{\varphi_m^4 \left[ \varphi_m^4 - N_m^2 \left((\partial_1 \mathcal f)^2 + (\partial_1 \mathcal f)^2 \right)\right]} } 
    \; \widetilde{V}_2 \\
    & \qquad + 
   \frac{N_m^2 (\partial_1 \mathcal f)(\partial_2 \mathcal f)}{\sqrt{\varphi_m^4 \left[ \varphi_m^4 - N_m^2 \left((\partial_1 \mathcal f)^2 + (\partial_1 \mathcal f)^2 \right)\right] \left[ \varphi_m^4 - N_m^2 (\partial_1 \mathcal f)^2\right]}}
    \; \widetilde{V}_1\,.
\end{split}
\end{equation}
Finally, we compute 
\begin{align*}
     \mathfrak{g}_m\left( \widetilde{V}_3, V_2 \right) 
    &= - \frac{ \varphi_m^2 N_m^2 \, (\partial_2 \mathcal f)(\partial_3 \mathcal f)}{\sqrt{ \left[\varphi_m^4 - N_m^2 \left((\partial_1 \mathcal f)^2 + (\partial_2 \mathcal f)^2 \right) \right]\left[ \varphi_m^4 - N_m^2 (\partial_1 \mathcal f)^2\right]}}
\end{align*}
and we consider the tangent vector field to $\mathcal{S}= \{ t = \mathcal f(\xx) \}$ that is orthogonal to $V_1$ and $V_2$, given by $\widetilde{V}_3 -  \mathfrak{g}_m\left( \widetilde{V}_3, V_2 \right) V_2 -  \mathfrak{g}_m\left( \widetilde{V}_3, V_1 \right) V_1$. By normalizing it, we obtain
\begin{equation}\label{ONF_3}
\begin{split}
 V_3 & = \sqrt{\frac{\varphi_m^4 - N_m^2 \left((\partial_1 \mathcal f)^2  + (\partial_2 \mathcal f)^2\right)}{\varphi_m^4 \left[ \varphi_m^4 - N_m^2 \sum_{i = 1}^3(\partial_i \mathcal f)^2 \right]} } \; \widetilde{V}_3 \\
 & \qquad +  \frac{  N_m^2 \, (\partial_2 \mathcal f)(\partial_3 \mathcal f)}{\sqrt{ \varphi_m^4 \left[ \varphi_m^4 - N_m^2 \sum_{i = 1}^3(\partial_i \mathcal f)^2 \right]\left[\varphi_m^4 - N_m^2 \left((\partial_1 \mathcal f)^2 + (\partial_2 \mathcal f)^2 \right) \right] }} \; \widetilde{V}_2 \\
 & \qquad + \frac{N_m^2 (\partial_1 \mathcal f)   (\partial_3 \mathcal f)}{\sqrt{\varphi_m^4 \left[ \varphi_m^4 - N_m^2 \sum_{i = 1}^3(\partial_i \mathcal f)^2 \right] \left[\varphi_m^4 - N_m^2 \left((\partial_1 \mathcal f)^2 + (\partial_2 \mathcal f)^2 \right) \right]}  }  \; \widetilde{V}_1
\end{split}
\end{equation}
The frame defined by \eqref{ONF_t}, \eqref{ONF_1}, \eqref{ONF_2} and \eqref{ONF_3} is orthonormal. 
Since the tangent frame $\{ V_1, V_2, V_3 \}$  is given by linear combinations of the frame $ \{ \widetilde{V}_1, \widetilde{V}_2, \widetilde{V}_3 \}$, as defined in \eqref{Frame_notON}, with coefficients bounded by the hypothesis \eqref{bound_hypothesis}, it will be sufficient to verify the desired conditions for the latter. 

Observe now that $2\alpha^{-1} \dd \alpha = \alpha^{-2} \dd\alpha^2$, so we compute
\small
\begin{align*}
    \alpha_{\tau}^{-2}\dd\alpha_{\tau}^2 & =   \frac{2}{N_m} \, \dd N_m + \frac{\alpha_{\tau}^2}{N_m^2} \Bigg( \frac{2N_m\sum_{i = 1}^3 (\partial_i \mathcal f)^2 }{\varphi_m^4} \dd N_m  - \frac{4N_m^2 \sum_{i = 1}^3 (\partial_i \mathcal f)^2 }{\varphi_m^5} \dd\varphi_m + \frac{2N_m^2}{\varphi_m^4}  \partial_i \mathcal f\,\partial^2_{ij}\mathcal f \, dx^j \Bigg)
\end{align*}
\normalsize
where the repeated indices have to be understood as in the Einstein summing convention with respect to the standard Euclidean metric. 
For the following estimates, it will be useful to notice that by the definition \eqref{def: v} of $\vv_t$, the bound in the hypotheses ensures, in particular, that
\[ \frac{N_m^2}{\varphi_m^4} \sum_{i = 1}^3 (\partial_i \mathcal f)^2 < 1 - \frac{1}{C_{\mathcal f}^2}\,.\]
We then find that
\begin{align}\label{estimate_alpha1_T}
\nonumber
    \vert &  \mathfrak{g}_m\left(\alpha_{\tau}^{-2} \dd\alpha_{\tau}^2, -\alpha_{\tau} \dd \tau \right)\vert  =  
     N_m \vv_t(\mathcal f) \,\varphi_m^{-4} \, \partial_k \mathcal f \, \Bigg\vert \frac{2}{N_m} \frac{x^k}{\vert \xx \vert} \frac{4m}{\left( 2 \vert \xx \vert + m\right)^2} + \\ \nonumber
    & \; \; + \frac{\alpha_{\tau}^2}{N_m^2} \Bigg( \frac{2N_m \sum_{i = 1}^3 (\partial_i \mathcal f)^2 }{\varphi_m^4} \frac{x^k}{\vert \xx \vert} \frac{4m}{\left( 2\vert \xx \vert + m\right)^2}  + \frac{4N_m^2 \sum_{i = 1}^3 (\partial_i \mathcal f)^2 }{\varphi_m^5} \frac{x^k}{\vert \xx \vert} \frac{m}{2\vert \xx \vert^2}  +\frac{2N_m^2}{\varphi_m^4}  \partial_i \mathcal f\,\partial^2_{ik}\mathcal f \Bigg) \Bigg\vert \\
    & \,\leq 
    C_{\mathcal f} \sqrt{1 - \frac{1}{C_{\mathcal f}}} \left( \frac{8     }{
    \left( 4  \kappa^2 - 1\right) \, m} + \frac{C_{\mathcal f}^3 (2  \kappa + 1)^2}{(2  \kappa - 1)^3 m}  + \frac{2C_{\mathcal f}}{ \kappa^2 m} + \frac{2C_{\mathcal f}^2}{m}\right)
\end{align}
and similarly, we compute
\begin{align}\label{estimate_alpha1_V}
\nonumber
    \vert & \alpha_{\tau}^{-2} \dd\alpha_{\tau}^2(\widetilde{V}_i)\vert  \leq 
     \Bigg\vert \frac{2}{N_m} \frac{x^i}{\vert \xx \vert} \frac{4m}{\left( 2 \vert \xx \vert + m\right)^2} +  \frac{\alpha_{\tau}^2}{N_m^2} \Bigg( \frac{2N_m \sum_{i = 1}^3 (\partial_i \mathcal f)^2 }{\varphi_m^4} \frac{x^i}{\vert \xx \vert} \frac{4m}{\left( 2\vert \xx \vert + m\right)^2} \\ 
    & \qquad    + \frac{4N_m^2 \sum_{i = 1}^3 (\partial_i \mathcal f)^2 }{\varphi_m^5} \frac{x^i}{\vert \xx \vert} \frac{m}{2\vert \xx \vert^2} +\frac{2N_m^2}{\varphi_m^4}  \partial_j \mathcal f\,\partial^2_{ij}\mathcal f \Bigg) \Bigg\vert \, \leq \frac{C(C_{\mathcal{f}}, \kappa)}{m}
\end{align}
Combining \eqref{estimate_alpha}, \eqref{estimate_alpha1_T}, \eqref{estimate_alpha1_V} and the bounds on the coefficients of the tangential frame $V_1, V_2, V_3$ with respect to $\widetilde{V}_1, \widetilde{V}_2, \widetilde{V}_3$ we obtain that for some $C = C(C_{\mathcal f},  \kappa) > 0$, it holds
\begin{equation}\label{eq: est 1 t_b}
    {}^{\tau}\| \log \alpha_{\tau}\| + {}^{\tau}\| \,\dd \log \alpha_{\tau}\| \, m \leq C\,.
\end{equation}
For the sake of the following computations, we first observe that 
\begin{align*}
    \vert \partial_i \alpha_{\tau} \vert &= \vert \dd \alpha_{\tau} (\partial_{x^i}) \vert  = \vert  \alpha_{\tau} \vert  
    \vert \alpha_{\tau}^{-1}\dd \alpha_{\tau} (\widetilde{V}_i) \vert 
    \leq \frac{C(C_{\mathcal f},  \kappa)}{m}
\end{align*}
We then compute
\begin{align*}
    \nabla_{\partial_{x^i}} \mathcal{T}_{\mathcal f} 
    & = \left( \partial_i \alpha_{\tau} N_m^{-2} - \alpha_{\tau} \frac{2 \partial_i N_m}{N_m^3}  + \alpha_{\tau} N_m^{-2} \Gamma_{it}^t + \alpha_{\tau} \varphi_m^{-4} \partial_j \mathcal f \, \Gamma_{ij}^t \right) \partial_t \\
    & \qquad + \left(  \partial_i \alpha_{\tau} \varphi_m^{-4} \partial_k \mathcal f - \alpha_{\tau} \frac{4 \partial_i \varphi_m \, \partial_k \mathcal f}{\varphi_m^5} - \frac{\partial_{ik}^2f}{\varphi_m^4}  \alpha_{\tau} N_m^{-2} \Gamma_{it}^k + \alpha_{\tau} \varphi_m^{-4} \partial_j \mathcal f \, \Gamma_{ij}^k  \right) \partial_{x^k}
\end{align*}
and similarly
\begin{align*}
    \nabla_{\partial_{t}} \mathcal{T}_{\mathcal f} 
    & = \alpha_{\tau} \left(   N_m^{-2} \Gamma_{tt}^t + \varphi_m^{-4} \partial_j \mathcal f \, \Gamma_{tj}^t \right) \partial_t  +   \alpha_{\tau} \left(    N_m^{-2} \Gamma_{tt}^k + \varphi_m^{-4} \partial_j \mathcal f \,\Gamma_{tj}^k  \right) \partial_{x^k}\,.
\end{align*}
Since it is now easy to see that all the objects appearing in the coefficients  
as well as their derivatives, are bounded by hypothesis, we conclude that for some constant $C = C(C_{\mathcal f}, \kappa) > 0 $ we have
\begin{equation}\label{eq: est 2 t_b}
    {}^{\tau}\| \nabla \mathcal{T}_{\mathcal f}\| \, m + {}^{\tau}\| \nabla^2 \mathcal{T}_{\mathcal f}\| \, m^2 \leq C\,.
\end{equation}
The estimates \eqref{eq: est 1 t_b} and \eqref{eq: est 2 t_b} yield the desired result.
\end{proof}


\vspace{0.1cm}

\section{Boosted maximal hypersurfaces in Schwarzschild}\label{sec: height estimate and main proof}

\subsection{Boundary data and a priori 
estimates}

For the proof of the main \cref{thm: Main intro}, we will consider a sequence of solutions $\{M_k\}_{k=1}^{\infty}$ to Dirichlet problems with prescribed boundaries $\{ \Sigma_k\}_{k=1}^\infty$, which have to lie on the hypersurface $\mathcal{S}_{\bb}$.
In order to get uniform estimates in $k$, we will need to apply \cite[Theorem 3.1-$(iii)$]{Bart}, which requires, in addition to the properties of $\ttt_{\bb}$, the conditions \eqref{cond: boundary grad} on the boundary data. In the following we show that a natural explicit choice of $\{ \Sigma_k\}_{k=1}^\infty$  data satisfies such conditions.

\vspace{0.1cm}

\begin{lemma}\label{lemma: exhaustion}
    The spacelike hypersurface $\mathcal{S}_{\bb} = \{ \ttt_{\bb} = 0\}$ can be exhausted by a sequence of bounded connected sets $\{\mathcal{U}_k\}_{k=1}^{\infty}$ with boundaries $\partial \mathcal{U}_k = \Sigma_k = \Sigma_k^{\mathrm{I}} \sqcup \Sigma_k^{\mathrm{I}'}$ consisting in two connected components $\Sigma_k^{\mathrm{I}}$ in \emph{region}  $\mathrm{I}$ and $\Sigma_k^{\mathrm{I}'}$ in \emph{region}  $\mathrm{I}'$, given by $\Sigma_k = \{ r = km\} \cap \mathcal{S}_{\bb}$ such that there exists a uniform constant $C_{\Sigma} > 0$ independent of $k$ satisfying
    \begin{equation}\label{eq: condition_MC_exhaustion}
        {}^{\ttt_{\bb}}\| \vec{\mathcal{H}}(\Sigma_k)  \| \leq \frac{C_{\Sigma}}{m}
    \end{equation}
    where ${}^{\ttt_{\bb}}\| \cdot \| $ is the time-dependent Riemannian norm defined in \eqref{def: timeDep_RiemNorm} with respect to the time $\ttt_{\bb}$ of \cref{thm: global_time_function} and $\vec{\mathcal{H}}(\Sigma_k)$ is the spacetime mean curvature \eqref{def: spacetimeMeanCurv} of $\Sigma_k$.
\end{lemma}

\vspace{0.1cm}

\begin{proof}
    Consider the coordinate spheres $\{ r = km\} \cap \{ t = 0\} $ in Schwarzschildean coordinates for $k > K_{\bb}$ with $K_{\bb}$ as in \cref{prop: barriers_infty} and let
    $$\Sigma_k \coloneqq\{ (t, r, \theta_I) \, :\, r = km,\, t = u_{\bb}(\xx(r, \theta_I))\} \subset \mathcal{S}_{\bb}$$
    for $u_{\bb}$ defined in \cref{prop: barriers_infty}. We claim that these  boundaries of the exhaustion $\mathcal{S}_{\bb} \cap \{ r < km\} = \mathcal{U}_k$ satisfy the desired features. The tangent vectors to $\Sigma_k$ are spanned by $X_I = \partial_I + \partial_I u_{\bb} \partial_t$ and induced metric
    \[ (g_{\Sigma_k})_{IJ} =  (g_k)_{IJ} - h_m\, \partial_I u_{\bb} \, \partial_J u_{\bb}\,, \quad  (g_{\Sigma_k})^{IJ} = (g_k)^{IJ}  + \frac{h_m\, (g_k)^{IK}(g_k)^{JL} \partial_K u_{\bb} \,\partial_L u_{\bb}}{(1 - h_m \vert \dd u_{\bb} \vert^2_{g_k}) }\,,  \]
    where $g_k = (km)^2 g_{\nS^2}$ is the induced metric on $\{ r = km\} \cap \{ t = 0\} $.
    On the other hand, the normal bundle to $\Sigma_k$ is spanned by the vector fields
    \begin{align*}
        \mathcal{T}_{\bb} &= \frac{ \partial_t + h_m^2 \partial_r u_{\bb} \partial_r + h_m (g_k)^{IJ}\partial_J u_{\bb} \,\partial_I}{\sqrt{h_m} \sqrt{1 - h_m^2 (\partial_r u_{\bb})^2-h_m \vert \dd u_{\bb} \vert^2_{g_k}}} \,, \\
        V_r &= \frac{\sqrt{h_m} \left( \left(1 -h_m \vert \dd u_{\bb} \vert^2_{g_k}\right)\partial_r + \partial_r u_{\bb} \partial_t  + h_m \partial_r u_{\bb} (g_k)^{IJ}\partial_J u_{\bb} \,\partial_I\right)}{\sqrt{1-h_m \vert \dd u_{\bb} \vert^2_{g_k}}\sqrt{1 - h_m^2 (\partial_r u_{\bb})^2-h_m \vert \dd u_{\bb} \vert^2_{g_k}}}
    \end{align*}
        which are elements of an orthonormal adapted frame in the sense of \eqref{def: timeDep_RiemNorm}.
    We then compute
    \[ \nabla_{X_I} X_J = \nabla_{\partial_I} \partial_J + \partial_{IJ} u_{\bb} \, \partial_t + \partial_I u_{\bb} \partial_J u_{\bb} \nabla_{\partial t} \partial_t = \partial_{IJ} u_{\bb} \, \partial_t + (\Gamma_{IJ}^r + \partial_I u_{\bb} \partial_J u_{\bb} \Gamma_{tt}^r) \, \partial_r + \Gamma_{IJ}^K \, \partial_K .  \]
    To compute the time-dependent Riemannian norm \eqref{def: timeDep_RiemNorm}, we need to check the boundedness independently of $k$ of the following two quantities
    \begin{align*}
        &\left|\G_m \left( \vec{\mathcal{H}}(\Sigma_k),  \mathcal{T}_{\bb}\right)\right|=  \left| -(g_{\Sigma_k})^{IJ} \G_m \left(\nabla_{X_I} X_J, \mathcal{T}_{\bb} \right) \right| \\
        & \quad \leq   \frac{\big|(g_{\Sigma_k})^{IJ} \big| \,\sqrt{h_m}}{\sqrt{1 - h_m^2 (\partial_r u_{\bb})^2-h_m \vert \dd u_{\bb} \vert^2_{g_k}}} \left(\left\vert \partial_{IJ} u_{\bb} \right\vert + \left\vert \partial_r u_{\bb}(\Gamma_{IJ}^r + \partial_I u_{\bb} \partial_J u_{\bb} \Gamma_{tt}^r) \right\vert +\left\vert \Gamma_{IJ}^K \partial_K u_{\bb} \right\vert\right) \,,
    \end{align*}
    \begin{align*}
         \left|\G_m \left( \vec{\mathcal{H}}(\Sigma_k),  V_r\right)\right| & =  \left| -(g_{\Sigma_k})^{IJ} \G_m \left(\nabla_{X_I} X_J, V_r \right) \right| \\
        &  \leq \frac{\big|(g_{\Sigma_k})^{IJ} \big| \,\sqrt{h_m}}{\sqrt{1 - h_m^2 (\partial_r u_{\bb})^2} \sqrt{1 - h_m^2 (\partial_r u_{\bb})^2-h_m \vert \dd u_{\bb} \vert^2_{g_k}}} \Bigg( h_m \left\vert \partial_r u_{\bb}\right\vert \left\vert \partial_{IJ} u_{\bb} \right\vert  \\ 
        & \qquad \qquad \qquad \qquad  +\left\vert \frac{1}{h_m}(\Gamma_{IJ}^r + \partial_I u_{\bb} \partial_J u_{\bb} \Gamma_{tt}^r) \right\vert + \left\vert \partial_r u_{\bb}\right\vert \left\vert \Gamma_{IJ}^K \partial_K u_{\bb} \right\vert\Bigg) \,.
    \end{align*}
    Since $\partial_r u_{\bb} = \frac{\partial x_i}{\partial r} \partial_i u_{\bb}$ is bounded, while $\partial_r u_{\bb} = \frac{\partial x_i}{\partial \theta_I} \partial_i u_{\bb} = O(km)$ and $\partial_{IJ} u_{\bb} = O(km)$, the above quantities decay to $0$ like $(km)^{-1}$ as $k \to \infty$. This concludes the proof.
\end{proof}

In the following, we prove the a priori height estimate for compact maximal hypersurfaces with boundary data that are invariant under the isometry $(\XX, \TT) \mapsto (-\XX, -\TT)$. In particular, this implies that the maximal hypersurface itself is invariant under such isometry and hence passes through the bifurcation sphere $\{ \XX = \TT = 0\}$. This is a direct consequence of the \emph{uniqueness of solution to the compact Dirichlet problem}, which is ensured on globally hyperbolic spacetimes satisfying the so-called \emph{Timelike Convergence Condition} $\mathfrak{Ric}(V, V) \geq 0$ for any timelike vector field $V$ (in our case, we have $\mathfrak{Ric} = 0$). The proof of such statement can be found e.g. in  \cite[Section 1]{Bartnik_NoCMC}.

\vspace{0.1cm}

\begin{theorem}[\textbf{A priori height estimate}]\label{thm: height estimate}
    Consider a boost parameter $\bb \in [0, 1)$ and the maximally extended Schwarzschild spacetime endowed with the time function $\ttt_{\bb}$ as in \cref{thm: global_time_function}. Let $M \hookrightarrow \left( \mathfrak{L}_m, \mathfrak{g}_m \right)$ be a compact spacelike acausal\footnote{An hypersurface is said to be \emph{acausal} if no two points on it can be connected via a causal curve.} hypersurface with vanishing mean curvature $\HH(M) = 0$ and with boundary $\partial M = \Sigma$ for some $\Sigma \subset \mathcal{S}_{\bb}$ consisting of two connected components and invariant under the isometry $(\XX, \TT) \mapsto (-\XX, -\TT)$.
    Then there exists 
    a constant $C_{\bb} =  C_{\bb}\left(\varepsilon, \bb, C_o, K_{\bb}, C_{K_{\bb}},   C_{\bar f _{\bb}}\right) > 0$, with $C_o$ as in \cref{thm: global_time_function}-$(i)$, $K_{\bb}$ as in \cref{prop: barriers_infty} and $C_{K_{\bb}}$ as in \cref{thm: global_time_function}-$ (iii)$, such that 
    \begin{equation}\label{eq: sup_estimate}
         \sup_{M} \big\vert u \big\vert \leq C_{\bb} m
    \end{equation}
    where $u$ is the height function \eqref{def: heightFct} of  $M$ with respect to $\ttt_{\bb}$. 
    
    Moreover, on $\big\{\vert \xx \vert > C_{\bb}^{\frac{1}{\varepsilon}}m\big\}$ in \emph{region} $\mathrm{I}$ the hypersurface $M$ lies between the barriers $\mathcal{W}^{\pm}_{\bb, C_{\bb}, \varepsilon}$ as given in \cref{prop: barriers_infty}. In \emph{region} $\mathrm{I}'$ the same holds by symmetry.
\end{theorem}

\vspace{0.1cm}

 The proof is based on \cite[Theorem 5.3]{Bart} with major differences, as   we  work with the spatial isotropic coordinates instead of finding explicit coordinates on the time slices of $\ttt_{\bb}$. Moreover, the asymptotic decays of geometric quantities on the time slices are more delicate in our case: for example, the decay of the second fundamental form of time slices of $\ttt_{\bb}$ decay as $O\left( \frac{\log(\vert \xx \vert)}{\vert \xx \vert^2}\right)$, rather than the rate $O\left(\frac{1}{\vert \xx \vert^2}\right)$. This has consequences, among others, on the fact that the test function involved in the integral estimate has to be heavily modified.
   
   In the following, for the sake of simplicity, we will ignore the dependence of any constant on $0 < \kappa_0 < \tfrac{3}{2}$ since due to \cref{prop: chrushingSing_barriers}, compact maximal hypersurfaces $M$ with boundaries on $\partial M \subset \mathcal{S}_{\bb}$ satisfy $M \subset \{ r > \kappa_0m\}$ due to \cref{prop: chrushingSing_barriers}.

\vspace{0.1cm}

\begin{proof}
    We first notice that once the height estimate \eqref{eq: sup_estimate} is established, the second part of the statement is easily ensured by the value attained by the barriers $\mathcal{W}_{\bb, C_{\bb}, \varepsilon}^{\pm}$ on $\big\{ \vert \xx \vert = C_{\bb}^{\frac{1}{\varepsilon}} m \big\}$ as in   \eqref{def: barriers_K} and the maximum principle in the static asymptotically flat ends. In fact, without loss of generality, we have $C_{\bb} > K_{\bb}$.

     To show \eqref{eq: sup_estimate},
    let us fix $0 < \varepsilon < 1$ and recall the height control of the barriers $\mathcal{W}_{\bb, K_{\bb}, \varepsilon}^{\pm}$ on $\{ \vert \xx \vert = K_{\bb} m \}$
    as in \eqref{def: barriers_K}. 
    By the staticity of the asymptotically flat ends, exactly as in \cite[Corollary 5.2]{Bart}, it is easy to show that 
    \[ \text{either } \;\;\sup_M |u| \leq  K_{\bb}^{\varepsilon} m \qquad \text{ or }\;\;  \sup_M |u|  \leq \sup_{M \cap \{ r \leq \widetilde{K}_{\bb} m\}} |u| \]
    for the Schwarzschildean radius $\widetilde{K}_{\bb} m$ corresponding to the isotropic radius $K_{\bb}m$ via \eqref{changeCoord_Iso_Schw}.
    Since the first option concludes the proof, we will assume that the second alternative holds. In the following, we will deal with the estimate on the supremum of $u$ by employing the upper barrier $\mathcal{W}^{+}_{\bb, K_{\bb}, \varepsilon}$. The estimate on the infimum can be treated analogously, by repeating the argument with the lower barrier $\mathcal{W}^{-}_{\bb, K_{\bb}, \varepsilon}$.

    From now on, we will work in isotropic coordinates in \emph{region} $\mathrm{I}$; the conclusions will clearly hold on \emph{region} $\mathrm{I}'$ as well, simply due to
    the symmetries of the Schwarzschild spacetime.  Consider now the time function on $\{ \vert \xx \vert > K_{\bb}m\}$ given by
    \begin{equation}
        \widetilde{\ttt}_{\bb} \coloneqq \ttt_{\bb} - m \,w^{+}_{K_{\bb}^{2 \varepsilon}}\left( \frac{\vert{\zz (\xx)}\vert}{m} \right)
    \end{equation}
     with $ w^{\pm}_{ K^{2 \varepsilon}}$ as in \eqref{eq: w}--\eqref{eq: w_-}  and the coordinates $\zz$ as in \eqref{def: changeCoords_xz}.
     We then define the set of points of $M$ lying above the barrier $\widetilde{M} = M \cap \{ \vert \xx \vert > K_{\bb}m\} \cap \{ \widetilde{\ttt}_{\bb} \geq 0\}$. Let  $\widetilde{u} \coloneqq \widetilde{\ttt}_{\bb} |_M$ the corresponding height function. Finally, let $\widetilde{\vv} = \vv_{\widetilde{\ttt}_{\bb}}$ be the quantity defined in \eqref{def: v} with respect to $\widetilde{\ttt}_{\bb}$, by $\widetilde{\K}^o$ and by $\widetilde{\HH}^o$ the second fundamental form and the mean curvature of the time slices of $\widetilde{\ttt}_{\bb}$, respectively.
    Then the boundary $\partial \widetilde{M} = \{ \widetilde{u} = 0\} \sqcup \mathfrak{C}_{\bb}$ where $\mathfrak{C}_{\bb} \subset M \cap \{ \vert \xx \vert = K_{\bb}m\}$.
    We now consider the equation \eqref{eq: MeanCurv_forIntEstimate} with respect to the time function $\widetilde{\ttt}_{\bb}$, which reads
    \[ 0 = \operatorname{div}_M (\alpha_{\widetilde{\ttt}_{\bb}} \, \operatorname{grad}_M \widetilde{u} ) + \operatorname{div}_M \widetilde{\mathcal{T}}_{\bb}\]
    where $\widetilde{\mathcal{T}}_{\bb}$ is the future-pointing unit normal to the time slices of $\widetilde{\ttt}_{\bb}$. Multiplying by $\alpha_{\widetilde{\ttt}_{\bb}} \psi\big( \frac{ \widetilde{u}}{m}\big)$, where $\psi \in C^1((0, + \infty)) \cap  C^0([0, + \infty))$, with $\psi(0) = 0$, $\psi(s), \psi'(s) > 0$ and integrating by parts, we obtain
    \begin{equation}
    \begin{split}\label{computation: integral_est}
        0 
        & =  - \int_{\mathfrak{C}_{\bb}} \alpha_{\widetilde{\ttt}_{\bb}}^2 \psi\Big( \frac{ \widetilde{u}}{m}\Big) \, \langle \operatorname{grad}_M \widetilde{u},\nu_{\bb} \rangle \; \dd \sigma_{\mathfrak{C}_{\bb}} + \int_{\widetilde{M}} \alpha_{\widetilde{\ttt}_{\bb}} \psi\Big( \frac{ \widetilde{u}}{m}\Big) \widetilde{\vv}\,
       \widetilde{\HH}^o\; \dd V_{\scriptscriptstyle M}  \\
        & \qquad  - \int_{\widetilde{M}}  \psi\Big( \frac{ \widetilde{u}}{m}\Big) (\widetilde{\vv}^2 - 1) 
        \left( \frac{\psi'\big( \frac{ \widetilde{u}}{m}\big)}{ m\psi\big( \frac{ \widetilde{u}}{m}\big)} + \frac{\alpha_{\widetilde{\ttt}_{\bb}}}{\widetilde{\vv} + 1}\widetilde{\HH}^o -  \alpha_{\widetilde{\ttt}_{\bb}} \widetilde{\K}^o_{ \widetilde{u}} + \widetilde{\mathcal{T}}_{\bb}( \alpha_{\widetilde{\ttt}_{\bb}})\right) \; \dd V_{\scriptscriptstyle M} \,,
        \end{split}
    \end{equation}
     where $\widetilde{\K}^o_{\widetilde{u} } = \widetilde{\K}^o\left(\frac{\operatorname{grad}_{\widetilde{\ttt}_{\bb}} \widetilde{u}}{\vert \operatorname{grad}_{\widetilde{\ttt}_{\bb}} \widetilde{u} \vert}, \frac{\operatorname{grad}_{\widetilde{\ttt}_{\bb}} \widetilde{u}}{\vert \operatorname{grad}_{\widetilde{\ttt}_{\bb}} \widetilde{u} \vert} \right)$ for $\operatorname{grad}_{\widetilde{\ttt}_{\bb}}$ gradient on the time slices of $\widetilde{\ttt}_{\bb}$  and the vector field $\nu_{\bb} $ is the unit normal to $\mathfrak{C}_{\bb}$ in $\widetilde{M}$ pointing towards the asymptotic flat end. 
      In the above computation we used \eqref{eq: alpha_gradu_v} and \eqref{eq: divT} and the identity $1 = \widetilde{\vv}- \frac{\widetilde{\vv}^2 - 1}{\widetilde{\vv} + 1}$. Moreover, we remark that the integration by parts has to be understood as performed on the truncated domain $\{ \widetilde{u} = \delta \}$ with the limit $\delta \to 0$.

      In the following, we will use the notation $C$ to collect all the constants depending only on $ \varepsilon, \bb, K_{\bb}, C_{\bar f_{\bb}}$ and independent of $m$.
      To estimate the last term in \eqref{computation: integral_est}, we first need to control  $\alpha_{\widetilde{\ttt}_{\bb}}$ and $\widetilde{\K}^o$. Since the time function $\ttt_{\bb} = t - \bb x_1 - m \bar f_{\bb}\left( \frac{\xx}{m} \right)$ on $\{ \vert \xx \vert \geq K_{\bb} m\}$, we compute
      \[ \dd \widetilde{\ttt}_{\bb}  = \dd t - \bb \dd x_1 - m\,\frac{\partial}{\partial x^i}\bar f_{\bb}\left( \frac{\xx}{m} \right) \dd x_i + K^{2 \varepsilon} \,\varepsilon \, m^{1+\varepsilon}\frac{ x_i (1 - \bb^2 \delta_{i1})}{\vert \zz (\xx) \vert^{2 +\varepsilon}} \dd x_i\,.  \]
      By the definition \eqref{def: lapse}, using the properties of $\bar f_{\bb}$ as in \eqref{def: C_f} and the estimates \eqref{estimates N phi}, we obtain
      \[ \left\vert \alpha_{\widetilde{\ttt}_{\bb}}^2 -\frac{1}{1 - \bb^2} \right\vert \leq C \frac{m \log \left( \frac{\vert \xx \vert }{m}\right)}{\vert \xx \vert} \,, \qquad   \left\vert \partial_i \alpha_{\widetilde{\ttt}_{\bb}} \right\vert \leq C \frac{m \log \left( \frac{\vert \xx \vert }{m}\right)}{\vert \xx \vert^2} \,. \]
      Moreover, since the tangent vectors are spanned by  $X_i = \partial_i + \left( \bb \delta_{i1} + \zeta_i(\xx)  \right)\partial_t$ with $\vert \zeta_i(\xx) \vert \leq C \frac{m \log \left( \frac{\vert \xx \vert }{m}\right)}{\vert \xx \vert}$ and by the control on the Christoffel symbols as in \eqref{estimate christoffel}, we can also estimate
      \[ \vert \widetilde{\K}^o(X_i, X_j) \vert = \vert \alpha_{\widetilde{\ttt}_{\bb}} \dd \widetilde{\ttt}_{\bb} \left( \olNabla_{X_i} X_j\right) \vert \leq C \frac{m \log \left( \frac{\vert \xx \vert }{m}\right)}{\vert \xx \vert^2}\,. \]
      Finally, we have
      \[ \left\vert \widetilde{\mathcal{T}}_{\bb}( \alpha_{\widetilde{\ttt}_{\bb}}) \right\vert = \left\vert \alpha_{\widetilde{\ttt}_{\bb}} \, \G_m \left( \dd \widetilde{\ttt}_{\bb}, \partial_i \alpha_{\widetilde{\ttt}_{\bb}} \dd x_i  \right)\right \vert \leq C \frac{m \log \left( \frac{\vert \xx \vert }{m}\right)}{\vert \xx \vert^2}\,. \]
      Hence, by recalling the decay control on $\widetilde{\HH}^o$ computed in \eqref{computation: meanCurv_barriers_infty}, we finally get 
       \begin{equation}\label{estimate_C1}
        \left| \frac{\alpha_{\widetilde{\ttt}_{\bb}}}{\widetilde{\vv} + 1}\widetilde{\HH}^o -  \alpha_{\widetilde{\ttt}_{\bb}} \widetilde{\K}^o_{ \widetilde{u}} + \widetilde{\mathcal{T}}_{\bb}( \alpha_{\widetilde{\ttt}_{\bb}})\right| \leq  C_1 \frac{ m\log\left( \frac{\vert \xx \vert }{m}\right)}{\vert \xx \vert^{2}}\,,
    \end{equation}
    for some constant $C_1 = C_1( \varepsilon, \bb, K_{\bb}, C_{\bar f_{\bb}}) > 0$.

    To estimate the boundary term in \eqref{computation: integral_est}, using
    the identity \eqref{eq: alpha_gradu_v}, we observe that it is necessary to control $\mathcal{\widetilde{\vv}}$ at points on $\mathfrak{C}_{\bb}$. This can be done using \cite[Theorem 3.1-$(iv)$]{Bart}, up to showing that the hypotheses are satisfied. To this end, we first claim that for any $\p \in \mathfrak{C}_{\bb}$ there holds
    \begin{equation}\label{eq: oscillation u}
        \sup_M u - u(\p) \leq (C_{K_{\bb}} + C_{\mathfrak{C}_{\bb}})m\,,
    \end{equation}
    where $C_{K_{\bb}}$ is the constant in \eqref{cond: Bartink_Interior_t_b} and $C_{\mathfrak{C}_{\bb}} = C_{\mathfrak{C}_{\bb}}(\bb, K_{\bb}, C_{\bar f_{\bb}}) > 0$ is a constant controlling the oscillation of $u$ on $\mathfrak{C}_{\bb}$.  The definition of $C_{K_{\bb}}$ ensures $\sup_M u - \sup_{\mathfrak{C}_{\bb}} u \leq C_{K_{\bb}} m$ and hence we are left with finding the constant $ C_{\mathfrak{C}_{\bb}}$. Since $M$ is an acausal spacelike hypersurface, we can write $M \cap \{ \vert \xx \vert > K_{\bb} m\} =\{ t = U(\xx) \}$ as a graph of a smooth function $U$ in isotropic coordinates over the domain $\{ t = 0, \, K_{\bb} m < \vert \xx \vert \leq k_{\Sigma} m \}$, where $k_{\Sigma} m$ is the isotropic radius on which the boundary $\Sigma$ lies. Now, the spacelike condition implies \begin{equation}\label{gradient_euclidean_bound}
         \sum_{i = 1}^3 (\partial_i U)^2 \leq \varphi_m^4 \,N_m^{-2} \leq \tfrac{\left( 2K_{\bb}+1\right)^6}{(2K_{\bb})^4 (2K_{\bb}-1)^2} \leq 2\,,
     \end{equation}
     where we have assumed without loss of generality that $K_{\bb} \geq 6$. Considering any two points $  \yy, \zz \in\{ t = 0, \, \vert \xx \vert = K_{\bb} m \}$ and the shortest curve connecting them, we find 
     \[ \vert U(\zz) - U(\yy) \vert \leq 2\pi K_{\bb} m\,. \]
     Moreover, given the control on the derivative of $\bar f_{\bb}$ as in \eqref{def: C_f}, it holds
     \[ \quad \vert\bb z_1 - \bb y_1 \vert \leq 2 \bb K_{\bb} m\,, \qquad \;\;\left\vert \bar f_{\bb}\left( \frac{\zz}{m}\right) - \bar f_{\bb}\left(\frac{\yy}{m}\right) \right\vert \leq  C_{\bar f_{\bb}} \pi \log(K_{\bb}) m\,.   \]
     Thus, by defining $C_{\mathfrak{C}_{\bb}}=C_{\mathfrak{C}_{\bb}}(\bb, K_{\bb}, C_{\bar f_{\bb}}) \coloneqq 2(\pi + \bb) K_{\bb} + C_{\bar f_{\bb}} \pi \log(K_{\bb}) $, we have shown \eqref{eq: oscillation u}.
     We then define 
      the following quantity $\mu$. Since otherwise we could directly conclude, we assume that it satisfies the following inequality
    \begin{equation}\label{def: mu_m}
        \mu  m \coloneqq \sup_M u -  (C_{K_{\bb}} + C_{\mathfrak{C}_{\bb}})  m - K_{\bb}^{\varepsilon} m > 2K_{\bb} m \,.
    \end{equation}
    Combining this lower bound on $\mu$ with the properties of $\ttt_{\bb}$ in \cref{thm: global_time_function} and $\sup_{\partial M} u = 0$ and by the symmetry of $M$, we see that for any $\p \in \mathfrak{C}_{\bb} \subset \{ \vert \xx \vert = K_{\bb} m\}$ it holds
   \begin{equation}\label{estimate_u_p}
        u(\p) - \sup_{\partial M} u \,\geq\,  \sup_M u -  (C_{K_{\bb}} + C_{\mathfrak{C}_{\bb}}) m  \,=\, \mu m  + K_{\bb} ^{\varepsilon} m\, >\, m \,>\, 0\,.
    \end{equation}
      Now, the control provided by \eqref{estimate_u_p}, combined with the property \eqref{eq: spacetimeConds_ttt_b} of $\ttt_{\bb}$ as in \cref{thm: global_time_function}, ensures that the hypothesis of \cite[Theorem 3.1-$(iv)$]{Bart} are satisfied, and  this yields  the pointwise estimate  for any $\p \in \mathfrak{C}_{\bb}$
     \[ \vv_{\ttt_{\bb}}(M)\big|_{\p} \leq C\left(\varepsilon, \bb, K_{\bb}, C_{\bar f_{\bb}}, C_{K_{\bb}}, C_o\right)\,. \]
     Notice that $\widetilde{\ttt}_{\bb} - \ttt_{\bb} = - m \,w^{+}_{K_{\bb}^{2 \varepsilon}}\left( \frac{\vert{\zz (\xx)}\vert}{m} \right)$ and therefore on $\{ \vert \xx \vert= K_{\bb} m\}$  the quantity $\vv_{\ttt_{\bb}}$ is comparable with $\widetilde{\vv}$ due to \cite[Lemma 3.3]{Bart} by a fixed constant depending on $\varepsilon, \bb$ and $K_{\bb}$. This, combined with \eqref{eq: alpha_gradu_v}, yields the estimate
     \begin{align}\label{estimate_C2}
    \left\vert\int_{\mathfrak{C}_{\bb}} \alpha_{\widetilde{\ttt}_{\bb}}^2 \psi\Big( \frac{ \widetilde{u}}{m}\Big) \, \langle \operatorname{grad}_M \widetilde{u},\nu_{\bb} \rangle \; \dd \sigma_{\mathfrak{C}_{\bb}} \right\vert & \,\leq \, \sup_{\widetilde{M}} \psi\Big( \frac{ \widetilde{u}}{m}\Big) \int_{\mathfrak{C}_{\bb}} \vert \alpha_{\widetilde{\ttt}_{\bb}} \vert \sqrt{\widetilde{\vv}^2 - 1} \; \dd \sigma_{\mathfrak{C}_{\bb}} 
         \,\leq \,C_2 m^2\, \sup_{\widetilde{M}}\psi\Big( \frac{ \widetilde{u}}{m}\Big) 
     \end{align}
     for a constant  $C_2 = C_2\left(\varepsilon, \bb, K_{\bb}, C_{\bar f_{\bb}}, C_{K_{\bb}}, C_o\right) > 0$.

     Finally, to estimate the second term in \eqref{computation: integral_est}, we consider  $\widetilde{M}$ as a graph in isotropic coordinates over a domain $\widetilde{\Omega}  \subset \{t = 0, \,K_{\bb} m \leq \vert \xx \vert \leq k_{\Sigma} m \}$. Since this will be a portion of the aforementioned graph via the function $U$,  the induced metric is given by
    $g_{ij} - N_m^2 \,\partial_i U \, \partial_j U$ and hence the volume form reads
    $ \dd V_{\scriptscriptstyle M} = \sqrt{\operatorname{det} g_{ij}} \vv_t^{-1}$, where $\vv_t$ is the quantity defined in \eqref{def: v} with respect to the static time $t$.
    Again, by  \cite[Lemma 3.3]{Bart}, which ensures that $\vv_t$ and $\mathcal{\widetilde{v}}$ are comparable by a constant depending on $\varepsilon, \bb, K_{\bb}, C_{\bar f_{\bb}}$,  combined with the computation \eqref{computation: meanCurv_barriers_infty}, we can estimate
    \begin{align}\label{estimate_C3}
        \int_{\widetilde{M}} \alpha_{\widetilde{\ttt}_{\bb}} \psi\Big( \frac{ \widetilde{u}}{m}\Big) \widetilde{\vv}\,
       \widetilde{\HH}^o\; \dd V_{\scriptscriptstyle M} \leq  -   C_3 m^{1 + \varepsilon}  \int_{\widetilde{\Omega}} \frac{\psi\big( \frac{ \widetilde{u}}{m}\big)}{\vert \xx \vert^{2 + \varepsilon}} \; \dd \xx
    \end{align}
    for $C_3 = C_3\left( \varepsilon, \bb, K_{\bb}, C_{\bar f_{\bb}}\right) > 0$\,.
    We can now plug \eqref{estimate_C1}, \eqref{estimate_C2} and \eqref{estimate_C3} into \eqref{computation: integral_est} and we obtain the following
    \begin{equation}\label{refinedEstimate_C123}
        C_3  m^{1 + \varepsilon} \int_{\widetilde{\Omega}} \frac{\psi\big( \frac{ \widetilde{u}}{m}\big)}{\vert \xx \vert^{2 + \varepsilon}} \, \dd \xx \leq C_2m^2\sup_{\widetilde{M}} \psi\Big( \frac{ \widetilde{u}}{m}\Big) -  \int_{\widetilde{M}}  \frac{\psi\big( \frac{ \widetilde{u}}{m}\big) }{m}(\widetilde{\vv}^2 - 1) 
        \left( \frac{\psi'\big( \frac{ \widetilde{u}}{m}\big)}{\psi\big( \frac{ \widetilde{u}}{m}\big)} -   \frac{ C_1 m\log\left( \frac{\vert \xx \vert }{m}\right)}{\vert \xx \vert^{2}} \right) \dd V_{ \scriptscriptstyle M}  \;\;
    \end{equation}

In order to neglect the second term on the right-hand side of \eqref{refinedEstimate_C123}, we need a suitable choice of function $\psi \in  C^1((0, + \infty)) \cap  C^0([0, + \infty))$, as follows
\begin{equation}
    \begin{split}
       \log\left(\psi\left(s \right)\right)=
       \begin{cases}
            4C_1 \frac{ \log(\mu + 1) \log s}{(\mu  + 1)^2} \, & \qquad \text{ if } 0 < s \leq 1 \\
           4 C_1 \left[ \frac{\log \left(\mu + 2 -s \right) + 1}{\mu + 2 -s} - \frac{\log(\mu + 1) + 1}{\mu + 1}\right]  \, & \qquad \text{ if } 1 \leq s \leq \mu   \\
           C_1 \log 2 \left( s - \mu\right) + 4C_1 \left[ \frac{\log 2 + 1}{2} - \frac{\log(\mu + 1) + 1}{\mu + 1} \right]\, & \qquad \text{ if }  \mu  \leq s \\
        \end{cases}
    \end{split}
    \end{equation}
    with  $\psi(0) = 0$, so that 
    \begin{equation}
    \begin{split}
       \frac{\psi'\left(s\right)}{\psi\left(s\right)} = \frac{d}{ds}\log(\psi\left(s \right))=
       \begin{cases}
            \frac{4C_1  \log(\mu + 1)}{(\mu  + 1)^{2} s} \, & \qquad \qquad \text{ if } 0 < s \leq 1 \\
           4 C_1 \frac{\log\left( \mu  + 2 - s\right)}{\left( \mu  + 2 - s \right)^2} \, & \qquad  \qquad \text{ if } 1 \leq s \leq \mu  \\
           C_1 \log 2\, & \qquad \qquad  \text{ if }  \mu \leq s  \\
        \end{cases}
    \end{split}
    \end{equation}
    We claim the non-negativity of the quantity $ \frac{\psi'\big( \frac{ \widetilde{u}}{m}\big)}{\psi\big( \frac{ \widetilde{u}}{m}\big)} -   \frac{ C_1 m\log\left( \frac{\vert \xx \vert }{m}\right)}{\vert \xx \vert^{2}} \geq 0$, as a consequence of the estimate on any $\q \in \widetilde{M}$ 
    \begin{equation}\label{eq: estimate tilde u}
         \vert \xx(\mathrm{q})\vert \geq \frac{1}{2} \left( \mu m  + 2 m - \widetilde{u}(\q) \right)\,.
    \end{equation}
To show the validity of \eqref{eq: estimate tilde u}, we estimate $\widetilde{u}(\q)$ for any $\q \in \widetilde{M}$, comparing the quantities involved with the ones of a point $\p \in \{ \vert \xx \vert = K_{\bb} m\}$ having the same angular coordinates, as follows
 \begin{align*}
     \widetilde{u}(\q) & = t(\q) - \bb x_1(\q) - m \bar f_{\bb}\left( \frac{\xx(\q)}{m} \right) - m \,w^{+}_{K_{\bb}^{2 \varepsilon}}\left( \frac{\vert{\zz (\q)}\vert}{m} \right) \\
     & \geq t(\p) - \frac{(2K_{\bb} + 1)^3}{4 K_{\bb}^2 (2K_{\bb} -1)}\left( \vert \xx (\q) \vert - K_{\bb} m\right) - \bb x_1(\p) - \bb \left( \vert \xx (\q) \vert - K_{\bb} m\right) \\
     & \qquad -  m \bar f_{\bb}\left( \frac{\xx(\p)}{m} \right) - C_{\bar f_{\bb}} \frac{\log(K_{\bb})}{K_{\bb}} \left( \vert \xx (\q) \vert - K_{\bb} m\right)  -  m \,w^{+}_{K_{\bb}^{2 \varepsilon}}\left( \frac{\vert{\zz (\p)}\vert}{m} \right)  \\
     & \geq u(\p) + \left( \frac{(2K_{\bb} + 1)^3}{4 K_{\bb}^2 (2K_{\bb} -1)} + \bb + C_{\bar f_{\bb}} \frac{\log(K_{\bb})}{K_{\bb}}   \right) \left( K_{\bb} m - \vert \xx (\q) \vert \right)  - K_{\bb}^{\varepsilon} m \\
     & \geq \mu m + 2m - 2\vert \xx (\q) \vert
 \end{align*}
 where we estimated the quantity $\frac{(2K_{\bb} + 1)^3}{4 K_{\bb}^2 (2K_{\bb} -1)} + \bb + C_{\bar f_{\bb}} \frac{\log(K_{\bb})}{K_{\bb}}    \leq 2 $ up to enlarging $K_{\bb} \gg 6$, depending on the fixed parameter and constant $\bb, C_{\bar f_{\bb}}$.
 Notice that in the above computations we have used again \eqref{gradient_euclidean_bound}, the control on the derivative of $\bar f_{\bb}$ as in \eqref{def: C_f} and the estimate \eqref{estimate_u_p}.
 Now, \eqref{eq: estimate tilde u} implies that
\begin{align*}
    C_1 \frac{m\log\left( \frac{\vert \xx \vert}{m}\right)}{\vert \xx \vert^2} & \leq 4C_1 \frac{\log\left( \mu  + 2 - \frac{ \widetilde{u}}{m}\right)}{m\left( \mu  + 2 - \frac{ \widetilde{u}}{m} \right)^2} \leq 4C_1 \frac{\log\left( \mu  + 1\right)}{\left( \mu  + 1\right)^2  \widetilde{u}}   = \frac{\psi'\big( \frac{ \widetilde{u}}{m}\big)}{m \psi\big( \frac{ \widetilde{u}}{m}\big)} \qquad   & \text{ if } \, 0 < \widetilde{u} \leq m \\
    C_1 \frac{m\log\left( \frac{\vert \xx \vert}{m}\right)}{\vert \xx \vert^2} & \leq 4C_1 \frac{\log\left( \mu  + 2 - \frac{ \widetilde{u}}{m}\right)}{m\left( \mu  + 2 - \frac{ \widetilde{u}}{m} \right)^2}  = \frac{\psi'\big( \frac{ \widetilde{u}}{m}\big)}{m \psi\big( \frac{ \widetilde{u}}{m}\big)} \qquad   & \text{ if } \, m \leq \widetilde{u} \leq \mu m \\
    C_1 \frac{m\log\left( \frac{\vert \xx \vert}{m}\right)}{\vert \xx \vert^2} & \leq C_1 \frac{\log 2}{4 m} \leq \frac{\psi'\big( \frac{ \widetilde{u}}{m}\big)}{m\psi\big( \frac{ \widetilde{u}}{m}\big)} \qquad &  \text{ if } \, \widetilde{u} \geq \mu m\,,
\end{align*}
 where the last inequality is valid due to  $\vert \xx \vert \geq K_{\bb} m > 2m$.  Thus, the second term on the right-hand side of \eqref{refinedEstimate_C123} can be neglected.

 Finally, we see that
    \[ \widetilde{u} \leq \sup_M u + K_{\bb}^{\varepsilon} m = \mu m  + (C_{K_{\bb}} + C_{\mathfrak{C}_{\bb}})  m+ 2K_{\bb}^{\varepsilon} m \]
    and consequently, by the definition of $\psi$, it holds  
    \[ \sup_M \psi\Big( \frac{ \widetilde{u}}{m}\Big) \leq \psi\left(\mu  + (C_{K_{\bb}} + C_{\mathfrak{C}_{\bb}})  + 2K_{\bb}^{\varepsilon} \right) \leq  C_4\,,\]
    where $C_4 = C_4\left(\varepsilon, \bb,  K_{\bb},   C_{\bar f _{\bb}},  C_{K_{\bb}}, C_o\right) = \mathrm{e}^{2C_1} 2^{C_1 \left(\mu  + (C_{K_{\bb}} + C_{\mathfrak{C}_{\bb}})  + 2K_{\bb}^{\varepsilon} +2 \right) } $.
 Since the estimate \eqref{eq: estimate tilde u} implies that $\widetilde{u} \geq \mu m + 2 m- \mu m= 2 m$ for $\vert \xx \vert \leq \tfrac{\mu m }{2}$, we conclude that on $\left\{ K_{\bb} m \leq \vert \xx \vert \leq \tfrac{\mu m}{2}\right\}$ we have $\psi\big( \frac{ \widetilde{u}}{m}\big) \geq 1$ and hence 
 \eqref{refinedEstimate_C123} gives 
    \begin{equation*}
        C_3 m^{1 + \varepsilon} \int\displaylimits_{K_{\bb}m \leq \vert \xx \vert \leq \tfrac{\mu m}{2}} \frac{1}{\vert \xx \vert^{2+ \varepsilon}} \dd \xx \leq  C_2 m^2 \sup_{\widetilde{M}} \psi\Big( \frac{ \widetilde{u}}{m}\Big) \leq C_2 C_4 m^2\,. 
    \end{equation*}
    After integrating, we obtain the following estimate on $\mu$
    \[ \frac{4\pi C_3}{1 - \varepsilon} \left[\left( \frac{\mu}{2} \right)^{1 - \varepsilon}  - K_{\bb}^{1 - \varepsilon}\right] \leq  C_2  C_4\,,\]
    which yields
    \[ \sup_M u = \left(\mu + C_{K_{\bb}}+ K_{\bb}^{\varepsilon}\right)m \leq  C_{\bb} m\]
    for a constant $C_{\bb} =  C_{\bb}\left(\varepsilon, \bb,  K_{\bb},   C_{\bar f _{\bb}},  C_{K_{\bb}}, C_o\right) > 0$ and we conclude the proof.
 \end{proof}

\vspace{0.1cm}

\subsection{Proof of the main results}

The proof of the main result is along the lines of \cite[proof of Theorem 5.4]{Bart}.

\begin{proof}[Proof of \cref{thm: Main intro}]
Consider the boundaries $\{\Sigma_k \}_{k = 1}^{\infty}$ of an exhaustion 
as in \cref{lemma: exhaustion} of the spacelike hypersurface $\mathcal{S}_{\bb} = \{ \ttt_{\bb} = 0\}$ given in \cref{thm: global_time_function}. For any fixed $k$, the existence of a solution to the Dirichlet problem associated to the maximal surface equation with boundary data $\Sigma_k$, reading
    \begin{equation}\label{Dirichlet_exhaustion}
            \HH(M_k) = 0 \,, \qquad \partial M_k = \Sigma_k
    \end{equation}
    is ensured by 
    \cite[Theorem 4.3]{Bart}
    combined with the presence of the barriers $\B^{\pm}$ at the crushing singularity as in \cref{prop: chrushingSing_barriers}.   
    In particular, the Dirichlet problem is solved inside the compact region defined by the domain of dependence \footnote{\noindent The \emph{domain of dependence} of a spacelike hypersurface $\mathcal{S} \hookrightarrow (\LL, \mathfrak{g})$ is 
$\mathcal{D}(\mathcal{S}) = \overline{\mathcal{D}^+(\mathcal{S})} \cup \overline{\mathcal{D}^-(\mathcal{S})}$,
where $\mathcal{D}^+(\mathcal{S}) $ is the set of points $p \in \mathcal{L}$ such that every past-directed inextendible curve passing through $p$ intersect $\mathcal{S}$. The definition of $\mathcal{D}^-(\mathcal{S}) $ is analogous, considering future-directed inextendible causal curves.} and the causal past (future) \footnote{The \emph{causal past} $I^-(\mathcal{S})$ of a spacelike hypersurface $\mathcal{S} \hookrightarrow (\LL, \mathfrak{g})$ is the collection of points $p \in \mathcal{L}$ such that there exists a past-directed causal curve from a point of $\mathcal{S}$ to $p$. The definition of \emph{causal future} $I^+(\mathcal{S})$ is analogous, considering future-directed causal curves.}
of the future (past) barriers, that is  $\mathfrak{U}_k \coloneqq\mathcal{D}(\mathcal{U}_k) \cap I^-(\B^+) \cap I^+(\B^-)$. 
    The Banach space $C^{l, \alpha}(\mathcal{U}_k)$ of $l$-differentiable $\alpha$-H\"older functions over $\mathcal{U}_k$ is well defined once one introduces zero-shift coordinates in $\mathfrak{U}_k$ as in the proof of \cite[Theorem 4.1]{Bart}.
    Then, the existence of a smooth hypersurface
    $$M_k = \{ (\zz, u_k(\zz)) \,: \, \zz \in \mathcal{U}_k , \, u_k \in C^{\infty}(\mathcal{U}_k)\} \subset \mathfrak{U}_k \subset \LL_m$$ solving the problem \eqref{Dirichlet_exhaustion} relies on the \emph{Leray--Schauder fixed point theorem} in the form of  \cite[Theorem 4.1]{Bart} and on standard elliptic regularity theory.

    The height estimate in \cref{thm: height estimate}  is uniform in $k$ for the sequence  $\{M_k\}_{k=1}^{\infty}$, with respect to the time function $\ttt_{\bb}$ of \cref{thm: global_time_function}. 
    The properties of $\ttt_{\bb}$ in \cref{thm: global_time_function}-$(i)$, combined with the fact that all $\Sigma_k$ lie in the zero-time slice $\mathcal{S}_{\bb}$ and satisfy \eqref{eq: condition_MC_exhaustion}, enable to apply \cite[Theorem 3.1-(\textit{iv})]{Bart}.
    This provides a uniform control on the geometric quantity $\vv_{\ttt_{\bb}}(M_k)$, which in turn, by the \emph{De Giorgi--Nash theorem} \cite[Theorem 13.1]{GT} gives a uniform control in $k$
on the coefficients and the ellipticity of the operator \eqref{eq: max surf eq}. 
Now, the \emph{interior Schauder estimate} \cite[Theorem 6.2]{GT} ensures that for some $0 < \alpha < 1 $ the $C^{2, \alpha}$-norm
    of $\{u_k\}_{k=1}^{\infty}$ is controlled uniformly in $k$ on each compact subset of $\mathcal{S}_{\bb}$. Hence, by \emph{Arzelà--Ascoli}, one can extract a subsequence of graphs over any such compact subset, converging 
    to a maximal hypersurface.
    By a standard diagonal argument on a exhaustion of $\mathcal{S}_{\bb}$ and standard elliptic regularity theory, we obtain
    a smooth spacelike entire hypersurface $M_{\bb} \hookrightarrow \left( \LL_m, \mathfrak{g}_m \right)$.

By construction, the hypersurface $M_{\bb}$ lies between the barriers $\mathcal{W}^{\pm}_{\bb, C_{\bb}, \varepsilon}$ as given in \cref{prop: barriers_infty} on $\big\{\vert \xx \vert > \C m\big\}$ for $\C \coloneqq C_{\bb}^{\frac{1}{\varepsilon}}$ in \emph{region} $\mathrm{I}$; the same holds by symmetry in \emph{region} $\mathrm{I}'$. This follows from the fact that this property is satisfied by all the solutions $M_k$ due to the a priori height estimate in \cref{thm: height estimate}. Since $M_{\bb}$ is a graph over $\mathcal{S}_{\bb}$  
and hence it is acausal, writing  $M \cap \{ \vert \xx \vert > \C m\} =\{ t = U(\xx) \}$ as a graph of a smooth function $U$ in isotropic coordinates over the domain $\{ t = 0, \,  \vert \xx \vert > \C m  \}$ gives $| U(\xx) - u_{\bb}(\xx) | \leq \C^\varepsilon m \,| \xx |^{-\varepsilon}$. This proves that $M_{\bb} \in \mathcal{B}_{\bb}$.

By the same uniqueness argument displayed right before \cref{thm: height estimate}, since all the hypersurfaces $\{M_k\}_{k= 1}^\infty$ are invariant under the isometry $(\XX, \TT) \mapsto (-\XX, -\TT)$ and pass through the bifurcation sphere, we obtain $\{ \XX = \TT = 0\} \in M_{\bb}$, as desired.

\end{proof}
 
\begin{proof}[Proof of \cref{coroll: main}]
    Given the existence of the hypersurface $M_{\bb}$ in \cref{thm: Main intro}, consider its intersection to  \emph{region} $\mathrm{I}$, 
    corresponding to the exterior Schwarzschild.
    We have that the spacelike hypersurfaces for $c \in \nR$ 
    $$M_{\bb}^c \coloneqq \{ (\xx, t + c) \in (\LL_m)_{\mathrm{ext}} \, : \, (\xx, t) \in M_{\bb} \}$$
provide a foliation of (boosted) maximal hypersurfaces of the whole \emph{region} $\mathrm{I}$. In fact, it is easy to see that $M_{\bb}^c$ are spacelike and $\HH(M_{\bb}^c) = 0$, since the foliation flows in the direction of the Killing vector field $\partial_t$, where $t$ is the static time. 
\end{proof}


\vspace{0.3cm}

\begin{appendices}

\crefalias{section}{appendix}
   \section{Time-dependent Bartnik conditions}\label{Appendix: Bartink conds computations}
  
     Here we display some selected computations that show the validity of the \emph{time-dependent Bartnik conditions} \eqref{eq: spacetimeConds_tildeT} with respect to the new time function $\widetilde{\TT}$ of  \cref{prop: NewTime} in the region $\{ r > \kappa_0m\}$. The full set of computations can be found in the doctoral thesis of the author \cite[Appendix A.1]{CogoPhD}.

     The crucial observation emerging from these computations is that all the contributions involved can be reduce to radial monotonic quantities, mostly by using the relation between the Schwarzschildean coordinate $r \in (0, +\infty)$ with respect to the Kruskal ones \eqref{def: changeCoord_Kruskal_r}.
We then get estimates that are only dependent on a lower bound on the radial coordinate $r > \kappa_0 m$. We highlight that since $\xi_m(\kappa_0 m) = 32\kappa_0^{-1} e^{-\frac{\kappa_0}{2}} \, m^2\, $
    and the Kruskal--Szekeres coordinates $\TT, \XX$ are scaling invariant, all the estimates will scale as expected.

\addtocontents{toc}{\protect\setcounter{tocdepth}{0}}

\subsection{Preliminary computations}
Employing the identity \eqref{def: changeCoord_Kruskal_r} one can express the differential of the Schwarzschildean radial function in terms of the Kruskal coordinates as follows 
 \begin{equation}\label{dr_dT_dX}
     \dd r = \frac{8m^2}{r} e^{-\frac{r}{2m}} (\XX \dd \XX - \TT \dd \TT) = \frac{\xi_m}{4m}(\XX \dd \XX - \TT \dd \TT)\,.
 \end{equation}
 Now, comparing \eqref{dr_dT_dX} with $\dd r = \tfrac{\partial r}{\partial \TT} \dd \TT + \tfrac{\partial r}{\partial \XX} \dd \XX$, we find
\begin{equation}\label{eq: Kruskal_partial_r_TX}
    \frac{\partial r}{\partial \TT} = - \frac{8m^2}{r} e^{-\frac{r}{2m}} \,\TT = - \frac{\xi_m}{4m} \TT \,, \qquad \frac{\partial r}{\partial \XX} = \frac{8m^2}{r} e^{-\frac{r}{2m}} \,\XX = \frac{\xi_m}{4m} \XX \,.
\end{equation}
Finally, we compute
 \begin{equation}\label{eq: Kruskal_derivative_lapse}
     \frac{d \xi_m}{dr} = -\frac{\xi_m}{2m}\left( 1 + \frac{2m}{r}\right)
 \end{equation} 
 and, using all the above, we can write the \emph{Christoffel symbols} of the Kruskal metric as follows
 \allowdisplaybreaks
\begin{align}\label{Christoffel_Kruskal}
\nonumber
     \Gamma_{\TT \TT}^{\TT} = \Gamma_{\TT \XX}^{\XX} = \Gamma_{\XX \XX}^{\TT} = &  \,\frac{\xi_m^{-1}}{2} \frac{\partial \xi}{\partial r} \frac{\partial r}{\partial \TT} = \frac{\TT}{16 m^2} \left( 1 + \frac{ 2m}{r} \right) \xi_m \\ \nonumber
    \Gamma_{\TT \TT}^{\XX} = \Gamma_{\TT \XX}^{\TT} = \Gamma_{\XX \XX}^{\XX} = & \, \frac{\xi_m^{-1}}{2} \frac{\partial \xi}{\partial r} \frac{\partial r}{\partial \XX} = -\frac{\XX}{16 m^2} \left( 1 + \frac{ 2m}{r} \right) \xi_m \\ \nonumber
    \Gamma_{\TT I}^{\TT} = \Gamma_{\TT I}^{\XX} = \Gamma_{\XX I}^{\TT} = \Gamma_{\XX I}^{\XX} = & \Gamma_{\XX \XX}^{I} = \Gamma_{\TT \TT}^{I} = \Gamma_{\XX \TT}^{I} = \, 0 \\ 
    \Gamma_{IJ}^{\TT} = &  \frac{\xi_m^{-1}}{2} 2r \,\frac{\partial r}{\partial \TT} (g_{\nS^{n-1}})_{IJ} =- \frac{\TT}{4m\, r} \G_{IJ} \, \delta_{IJ}
    \\ \nonumber
    \Gamma_{IJ}^{\XX} = & -\frac{\xi_m^{-1}}{2} 2r \, \frac{\partial r}{\partial \XX} (g_{\nS^{n-1}})_{IJ} = - \frac{\XX}{4m\, r} \G_{IJ}\,\delta_{IJ}
    \\ \nonumber
    \Gamma_{I \TT}^{J} = & -\delta_I^J \,\frac{8m^2}{r^2} \, e^{-\frac{r}{2m}} \, \TT = - \delta_I^J \frac{\xi}{4m\, r} \TT\,
    \\ \nonumber
    \Gamma_{I \XX}^{J} = & \delta_I^J \,\frac{8m^2}{r^2} \, e^{-\frac{r}{2m}} \, \XX =  \delta_I^J \frac{\xi}{4m\, r} \XX\,.
\end{align}

\subsection{Estimates on $\log \alpha_{\widetilde{\TT}}$}

    \noindent
     The lapse computed in \eqref{def: lapse_tildeT} can be written via the identity \eqref{def: changeCoord_Kruskal_r} as follows
    \begingroup
     \allowdisplaybreaks
\begin{align*}
    \alpha_{\widetilde{\TT}}^2 &= \frac{\xi_m}{16m^2} \frac{(1 + \XX^2 - \TT^2)}{1 + f}
     = \frac{1}{1 + f} \left(  1 +\frac{2m}{r} ( \mathrm{e}^{-\frac{r}{2m}} -1) \right)  \, \geq \, 0  \,,
\end{align*}
\endgroup
with $f = f(\XX, \TT) \geq 0$ as in \eqref{def: f_tildeT}.
Now, we can easily compute the upper bound
\begin{equation}\label{comp: lapse_upper}
    \alpha_{\widetilde{\TT}}^2 < 1
\end{equation}
where we used that $\frac{1}{1 + f} \leq 1$ combined with $\lim_{r \to +\infty} \left(  1 +\frac{2m}{r} ( \mathrm{e}^{-\frac{r}{2m}} -1) \right) = 1$ and 
\begin{equation}
\label{derivative_xi_1XT}
     \frac{d}{dr} \frac{\xi_m (1 + \XX^2 - \TT^2)}{16m^2}  = 
    \frac{d}{dr}\left(  1 +\frac{2m ( \mathrm{e}^{-\frac{r}{2m}} -1)}{r}  \right) 
    = \frac{2m \,\mathrm{e}^{-\frac{r}{2m}}}{r^2} \left( \mathrm{e}^{\frac{r}{2m}} - \left( \frac{r}{2m} + 1\right)\right) >  0 .
\end{equation}
To get a lower bound, we plug \eqref{def: changeCoord_Kruskal_r} in \eqref{def: f_tildeT}, so that we estimate
 \[ \frac{1}{1 + f} = \frac{1}{1 + \frac{\TT^2}{1 + \XX^2}\frac{1}{\left(  1 + \left( \frac{r}{2m} -1\right) \mathrm{e}^{\frac{r}{2m}}\right)}} > \frac{1}{1 + \frac{1}{\left(  1 + \left( \frac{r}{2m} -1\right) \mathrm{e}^{\frac{r}{2m}}\right)}}\,,\]
 where the right hand side can be checked to be monotone increasing in $r > 0$.
 Hence, we obtain
 \begin{equation}\label{comp: lapse_lower}
     \alpha_{\widetilde{\TT}}^2 > \frac{1}{1 + \frac{1}{\left(  1 + \left( \frac{\kappa_0}{2} -1\right) \mathrm{e}^{\frac{\kappa_0}{2}}\right)}}  \left(  1 +\frac{2}{\kappa_0} ( \mathrm{e}^{-\frac{\kappa_0}{2}} -1) \right) > 0\,.
 \end{equation}
 Combining \eqref{comp: lapse_lower} with \eqref{comp: lapse_upper}, we conclude that on $\{ r > \kappa_0m\}$, there exists a constant $C = C(\kappa_0) > 0$ such that 
 \begin{equation}\label{A_Bound_log_alpha}
     {}^{\widetilde{\TT}} \| \log \alpha_{\widetilde{\TT}}\| \leq C\,.
 \end{equation}

\subsection{Estimates on $ \dd \log \alpha_{\widetilde{\TT}}$}

    Observing that $ \dd \log \alpha_{\widetilde{\TT}}  
    = \frac{1}{2}\alpha_{\widetilde{\TT}}^{-2} \dd\alpha_{\widetilde{\TT}}^2$ and using \eqref{dr_dT_dX} and \eqref{eq: Kruskal_derivative_lapse}, we compute
\begin{align*}
    \alpha_{\widetilde{\TT}}^{-2} \dd\alpha_{\widetilde{\TT}}^2 & = \left( \frac{16m^2 (1 + f)}{\xi_m \, (1 + \XX^2 - \TT^2)}\right) \; \dd \left( \frac{\xi_m}{16m^2} \,  \frac{(1 + \XX^2 - \TT^2)}{1 + f}\right) \\
    & = \left( \frac{2}{1 + \XX^2 - \TT^2} - \frac{\xi_m  \left( 1 + \frac{2m}{r}\right) }{8m^2}\right) \left( \XX \dd \XX - \TT \dd \TT\right) - \frac{\left( \partial_{\XX}f \dd \XX + \partial_{\TT}f \dd \TT\right)}{(1 + f)}\,, 
\end{align*}
where, by the definition of $f=f(\XX, \TT)$ in \eqref{def: f_tildeT}, it holds
\begin{align}\label{derivative_f_TX}
    \partial_{\TT} f  = 2 \, \frac{f}{\TT} \, \left( 1 +  (1+ \XX^2) f\right)\,, \qquad
    \partial_{\XX} f  = -2\,\frac{f^2}{\TT^2} \, \XX\left( (1 + \XX^2 - \TT^2) +  (1 + \XX^2)\right).
\end{align}

We now have to apply $\dd \log \alpha_{\widetilde{\TT}}$ to the elements of the adapted orthonormal frame $\{ \widetilde{\mathcal{T}}, \widetilde{V}_{\XX}, \widetilde{V}_{I} \}$, $I \in \{\theta, \varphi\}$ given by \eqref{def: normal_tildeT}, \eqref{def: VX_tildeT} and \eqref{def: VI_tildeT}. To compute  $ \dd \log \alpha_{\widetilde{\TT}}(\widetilde{\mathcal{T}})$, we first calculate
\begin{align*}
    \left( \XX \dd \XX - \TT \dd \TT\right) \left( \sqrt{1+ \XX^2} \, \partial_{\TT} + \frac{\XX \TT}{\sqrt{1 + \XX^2}} \partial_{\XX}\right) = \frac{\XX^2 \TT}{\sqrt{1 + \XX^2}} - \TT \sqrt{1 + \XX^2} =
    - \frac{\TT}{\sqrt{1 + \XX^2}}
\end{align*}
\begin{align*}
    \left( \partial_{\XX}f \dd \XX + \partial_{\TT}f \dd \TT\right) &\left( \sqrt{1+ \XX^2} \, \partial_{\TT} + \frac{\XX \TT}{\sqrt{1 + \XX^2}} \partial_{\XX}\right)  = \partial_{\TT}f \, \sqrt{1 + \XX^2} + \frac{\partial_{\XX}f \, \XX \TT}{\sqrt{1 + \XX^2}} \\
    & = \frac{2}{1 + \XX^2 - \TT^2} \,\frac{\TT}{\sqrt{1 + \XX^2}} \left( 1 + f  - \frac{\TT^2 \XX^2}{(1 + \XX^2)^2}\right)
\end{align*}
and hence we obtain
\begin{align} \nonumber 
    \big\vert \alpha_{\widetilde{\TT}}^{-2}\dd\alpha_{\widetilde{\TT}}^2\big(\widetilde{\mathcal{T}}\big) \big\vert &=  \frac{1}{\sqrt{\xi_m (1 + f) (1 + \XX^2 - \TT^2)}} \, \frac{|\TT|}{\sqrt{1 + \XX^2}} \,\Bigg\vert \left( \frac{\xi_m  }{8m^2} \left( 1 + \frac{2m}{r}\right) - \frac{2}{1 + \XX^2 - \TT^2} \right) \\ \nonumber
    & \qquad -\frac{2}{(1 + f) (1 + \XX^2 - \TT^2)}  \left( 1 + f  - \frac{\TT^2 \XX^2}{(1 + \XX^2)^2}\right)\Bigg\vert \\ \label{estimate: dlogalpha_T}
    & \leq \frac{1}{\sqrt{\xi_m (1 + \XX^2 - \TT^2)}} \, \left( \frac{\xi_m \left(1 + \frac{2m}{r} \right)}{8m^2}  + \frac{6}{1 + \XX^2 - \TT^2}    \right)   \; \leq  \frac{C}{m}
\end{align}
where we used that $\TT^2 < 1 + \XX^2$,  $f \geq 0$, the monotonicity given by \eqref{derivative_xi_1XT} and the fact that 
\begin{equation}
    \frac{1}{1 + \XX^2 - \TT^2} = \frac{1}{1 + \left(\frac{r}{2m} -1 \right) e^{\frac{r}{2m}}} \leq \frac{1}{1 + \left(\frac{\kappa_0}{2} -1 \right) e^{\frac{\kappa_0}{2}}}\,,
\end{equation}
again by monotonicity on $\{ r > \kappa_0m\}$. 

Estimating $ \alpha^{-2}\dd\alpha^2\big(\widetilde{V}_{\XX}\big) $ is analogous, while it is easy to see that 
$\dd\alpha_{\widetilde{\TT}}^2\big(\widetilde{V}_{I}\big)= 0.$ 
We then conclude that there exists a constant $C = C(\kappa_0) > 0$ such that in the region $\{ r > \kappa_0m\}$ it holds 
 \begin{equation}\label{A_Bound_d_log_alpha}
     {}^{\widetilde{\TT}} \| \dd \log \alpha_{\widetilde{\TT}}\| \leq C\,.
 \end{equation}

\subsection{Estimates on $ \nabla \widetilde{\mathcal{T}}$}
Here we show the estimate on  $\mathfrak{g}_m\left(\nabla_{\widetilde{V}_{\XX}} \widetilde{\mathcal{T}} \normalsize, \widetilde{V}_{\XX}\right)$. To this end, we first calculate
\begin{align*}
    \nabla_{\partial_{\XX}} \left(\frac{\XX \TT}{1 + \XX^2} \partial_{\XX} + \partial_{\TT} \right) & = \frac{\TT (1-\XX^2)}{(1 + \XX^2)^2}  \partial_{\XX} + \frac{\XX \TT}{1 + \XX^2}  \left( \Gamma_{\XX \XX}^{\XX} \partial_{\XX} +  \Gamma_{\XX \XX}^{\TT} \partial_{\TT} +  \Gamma_{\XX \XX}^{I} \partial_{I}\right) \\
    & \qquad +  \Gamma_{\XX \TT}^{\XX} \partial_{\XX} +  \Gamma_{\XX \TT}^{\TT} \partial_{\TT} +  \Gamma_{\XX \TT}^{I} \partial_{I}
\end{align*}
\begin{align*} 
    \nabla_{ \partial_{\TT} } \left(\frac{\XX \TT}{1 + \XX^2} \partial_{\XX} + \partial_{\TT} \right) &= \frac{\XX}{1 + \XX^2} \partial_{\XX} + \frac{\XX \TT}{1 + \XX^2} \left( \Gamma_{\TT \XX}^{\XX} \partial_{\XX} +  \Gamma_{\TT \XX}^{\TT} \partial_{\TT} +  \Gamma_{ \TT \XX}^{I} \partial_{I}\right) \\
    & \qquad + \Gamma_{\TT \TT}^{\XX} \partial_{\XX} +  \Gamma_{\TT \TT}^{\TT} \partial_{\TT} +  \Gamma_{\TT \TT}^{I} \partial_{I}\,.
\end{align*}
By the identities in \eqref{Christoffel_Kruskal} and $\overline{\Gamma}_{\XX \XX}^{\XX} = - \tfrac{\XX}{\TT} \overline{\Gamma}_{\TT \TT}^{\TT}$, we obtain
\begin{align}\label{nabla_VX_T_V_X}
\nonumber
    \mathfrak{g}_m\Bigg(  \nabla_{\partial_{\XX} + \frac{\XX \TT}{1 + \XX^2} \partial_{\TT}} & \left(\frac{\XX \TT}{1 + \XX^2} \partial_{\XX} + \partial_{\TT} \right) ,\,   \partial_{\XX} + \frac{\XX \TT}{1 + \XX^2} \partial_{\TT} \Bigg) = \\ 
    & = \xi_m \Bigg[ \frac{\TT}{(1 + \XX^2)^2} + \left( 1 - \frac{\XX^2 \TT^2}{(1 + \XX^2)^2} \right) \left( 1 - \frac{\XX^2}{1 + \XX^2} \right) \overline{\Gamma}_{\TT \TT}^{\TT}  \Bigg] \\ \nonumber
    & = \xi_m \,\frac{\TT}{(1 + \XX^2)^2} \left[ 1 + \frac{\xi_m \, (1 + \XX^2 - \TT^2) (1 + f)}{16m^2} \left( 1 + \frac{2m}{r}\right)\right]\,.
\end{align}
Finally, using the identity \eqref{nabla_VX_T_V_X}, we can estimate
\begin{align}\label{bound_NablaVX_T_VX}
\nonumber
    \Bigg\vert & \mathfrak{g}_m\left(\nabla_{\widetilde{V}_{\XX}} \widetilde{\mathcal{T}}, \widetilde{V}_{\XX}\right) \Bigg\vert  = \left\vert\frac{1 + \XX^2}{\xi_m (1+f) (1 + \XX^2 - \TT^2)} \; \mathfrak{g}_m\left( \nabla_{\partial_{\XX} + \frac{\XX \TT}{1 + \XX^2} \partial_{\TT}} \widetilde{\mathcal{T}} , \partial_{\XX} + \frac{\XX \TT}{1 + \XX^2} \partial_{\TT}\right) \right\vert \\ \nonumber
    & = \left\vert\left(\frac{1 + \XX^2}{\xi_m (1+f) (1 + \XX^2 - \TT^2)}\right)^{\frac{3}{2}} \mathfrak{g}_m\left(  \nabla_{\partial_{\XX} + \frac{\XX \TT}{1 + \XX^2} \partial_{\TT}} \left(\frac{\XX \TT}{1 + \XX^2} \partial_{\XX} + \partial_{\TT} \right) ,   \partial_{\XX} + \frac{\XX \TT}{1 + \XX^2} \partial_{\TT} \right) \right\vert \\ \nonumber
    & \leq \frac{1}{\sqrt{\xi_m \left(  1 + \XX^2 - \TT^2\right)^3}}  + \sqrt{\frac{\xi_m}{ 1 + \XX^2 - \TT^2}}\frac{1}{16m^2} \left( 1 + \frac{2m}{r}\right)\, \leq  \frac{C}{m}
\end{align}
where we relied on the usual inequalities $\TT^2 < 1 + \XX^2$ and $f \geq 0$ and on the monotonicity of the quantities involved.

Similarly, we can estimate  \small $\;\mathfrak{g}_m\left(\nabla_{\widetilde{V}_{\XX}} \widetilde{\mathcal{T}}, \widetilde{\mathcal{T}}\right)$, $\;\mathfrak{g}_m\left(\nabla_{\widetilde{V}_{\XX}} \widetilde{\mathcal{T}}, \widetilde{V}_{I}\right)$,  $\;\mathfrak{g}_m\left(\nabla_{\widetilde{V}_{J}} \widetilde{\mathcal{T}}, \widetilde{V}_{\XX}\right)$, $\;\mathfrak{g}_m\left(\nabla_{\widetilde{V}_{J}} \widetilde{\mathcal{T}}, \widetilde{\mathcal{T}}\right) $, $\,\mathfrak{g}_m\left(\nabla_{\widetilde{V}_{J}} \widetilde{\mathcal{T}}, \widetilde{V}_{I}\right) $,  $\,\mathfrak{g}_m\left(\nabla_{\widetilde{\mathcal{T}}} \widetilde{\mathcal{T}}, \widetilde{V}_{\XX}\right) $,  $\,\mathfrak{g}_m\left(\nabla_{\widetilde{\mathcal{T}}} \widetilde{\mathcal{T}}, \widetilde{\mathcal{T}}\right) $, $\,\mathfrak{g}_m\left(\nabla_{\widetilde{\mathcal{T}}} \widetilde{\mathcal{T}}, \widetilde{V}_{I}\right) $ \normalsize 
and hence we conclude that there exists a constant $C = C(\kappa_0) > 0$ such that in the region $\{ r > \kappa_0m\}$ it holds 
\begin{equation}\label{A_Bound_Nabla_T}
     {}^{\widetilde{\TT}} \| \nabla \widetilde{\mathcal{T}}\| \leq  \frac{C}{m}\,.
 \end{equation}

\subsection{Estimates on $ \nabla^2 \widetilde{\mathcal{T}}$}
Here we show how to estimate $\nabla^2_{\widetilde{V}_{\XX} \widetilde{V}_{\XX}} \widetilde{\mathcal{T}}$ against any vector field of the adapted frame.
By standard computations and the values of the Christoffel symbols \eqref{Christoffel_Kruskal}, one can verify that \small $\nabla_{\widetilde{V}_{\XX}} \widetilde{\mathcal{T}} = \mathfrak{g}_m \left( \nabla_{\widetilde{V}_{\XX}} \widetilde{\mathcal{T}}, \widetilde{V}_{\XX}\right) \widetilde{V}_{\XX}$ \normalsize and \small $\nabla_{\widetilde{V}_{\XX}} \widetilde{V}_{\XX} = - \mathfrak{g}_m \left( \nabla_{\widetilde{V}_{\XX}}  \widetilde{\mathcal{T}}, \widetilde{V}_{\XX} \right) \widetilde{\mathcal{T}}$. \normalsize Hence, it holds
\begin{align}\label{observations Hessian tilde t}
    \nabla&^2_{\widetilde{V}_{\XX} \widetilde{V}_{\XX}}  \widetilde{\mathcal{T}} = \nabla_{\widetilde{V}_{\XX}} \left(\nabla_{\widetilde{V}_{\XX}} \widetilde{\mathcal{T}} \right) - \nabla_{\nabla_{\widetilde{V}_{\XX}} \widetilde{V}_{\XX}} \widetilde{\mathcal{T}} \\ \nonumber
    & = \widetilde{V}_{\XX} \left(\mathfrak{g}_m \left( \overline{\nabla}_{\widetilde{V}_{\XX}} \widetilde{\mathcal{T}}, \widetilde{V}_{\XX}\right) \right) \widetilde{V}_{\XX}  - \mathfrak{g}_m \left( \overline{\nabla}_{\widetilde{V}_{\XX}} \widetilde{\mathcal{T}}, \widetilde{V}_{\XX}\right) \mathfrak{g}_m  \left( \overline{\nabla}_{\widetilde{V}_{\XX}}  \widetilde{\mathcal{T}}, \widetilde{V}_{\XX} \right) \widetilde{\mathcal{T}}  + \mathfrak{g}_m  \left( \overline{\nabla}_{\widetilde{V}_{\XX}}  \widetilde{\mathcal{T}}, \widetilde{V}_{\XX} \right)  \overline{\nabla}_{\widetilde{\mathcal{T}}} \widetilde{\mathcal{T}} 
\end{align}
Now, computing the second and third terms on the right-hand side against any vector of the adapted frame gives quantities that can be estimated just by \eqref{A_Bound_Nabla_T}. We are left with estimating  \small $\widetilde{V}_{\XX} \left(\mathfrak{g}_m\left( \nabla_{\widetilde{V}_{\XX}} \widetilde{\mathcal{T}}, \widetilde{V}_{\XX}\right) \right)$. \normalsize Using the explicit identity \eqref{nabla_VX_T_V_X}, we compute
\begin{align*}
   & \partial_{\XX}  \Bigg(\mathfrak{g}_m\left( \nabla_{\widetilde{V}_{\XX}} \widetilde{\mathcal{T}}, \widetilde{V}_{\XX}\right) \Bigg) = \\
   &= \partial_{\XX} \left[\frac{ \TT }{\sqrt{1 + \XX^2}} \, \frac{1}{ \sqrt{ \xi_m \left((1+f) (1 + \XX^2 - \TT^2)\right)^3}} \left( 1 + \frac{\xi_m \, (1 + \XX^2 - \TT^2) (1 + f)\left( 1 + \frac{2m}{r}\right)}{16m^2} \right) \right] \\
   & = - \frac{\TT \XX}{(1 + \XX^2)^{\frac{3}{2}}} \frac{1}{ \sqrt{ \xi_m \left((1+f) (1 + \XX^2 - \TT^2)\right)^3}} \left( 1 + \frac{\xi_m \, (1 + \XX^2 - \TT^2) (1 + f)\left( 1 + \frac{2m}{r}\right)}{16m^2} \right) \\
   & \qquad + \frac{ \TT }{\sqrt{1 + \XX^2}} \Bigg[ \left( - \frac{\frac{d \xi_m}{dr} \frac{\partial r}{\partial \XX}}{ 2 \sqrt{\left( \xi_m (1+f) (1 + \XX^2 - \TT^2)\right)^3}} -  \frac{ \left(\frac{d \xi_m}{dr} \left(1 + \frac{2m}{r} \right) - \frac{4m}{r^2} \xi_m \right) \frac{\partial r}{\partial \XX}}{ 32m^2 \sqrt{ \xi_m (1+f) (1 + \XX^2 - \TT^2)}}\right)  \\
   & \qquad - \frac{\partial_{\XX}f (1+ \XX^2 - \TT^2) + 2\XX(1 +f)}{2 \sqrt{\xi_m \left( (1+f) (1 + \XX^2 - \TT^2)\right)^3}} \left(  \frac{3}{(1+f) (1 + \XX^2 - \TT^2)} + \frac{ \xi_m \left( 1 + \frac{2m}{r}\right)}{16m^2}\right) \Bigg]
\end{align*}
and 
\begin{align*}
   & \partial_{\TT}  \Bigg(\mathfrak{g}_m\left( \nabla_{\widetilde{V}_{\XX}} \widetilde{\mathcal{T}}, \widetilde{V}_{\XX}\right) \Bigg) = \\
   &= \partial_{\TT} \left[\frac{ \TT }{\sqrt{1 + \XX^2}} \, \frac{1}{ \sqrt{ \xi_m \left((1+f) (1 + \XX^2 - \TT^2)\right)^3}} \left( 1 + \frac{\xi_m \, (1 + \XX^2 - \TT^2) (1 + f)\left( 1 + \frac{2m}{r}\right)}{16m^2} \right) \right] \\
   & = \frac{1 }{\sqrt{1 + \XX^2}} \frac{1}{ \sqrt{ \xi_m \left((1+f) (1 + \XX^2 - \TT^2)\right)^3}} \left( 1 + \frac{\xi_m \, (1 + \XX^2 - \TT^2) (1 + f)\left( 1 + \frac{2m}{r}\right)}{16m^2} \right) \\
   & \qquad + \frac{ \TT }{\sqrt{1 + \XX^2}} \Bigg[ \left( - \frac{\frac{d \xi_m}{dr} \frac{\partial r}{\partial \TT}}{ 2 \sqrt{\left( \xi_m (1+f) (1 + \XX^2 - \TT^2)\right)^3}} -  \frac{ \left(\frac{d \xi_m}{dr} \left(1 + \frac{2m}{r} \right) - \frac{4m}{r^2} \xi_m \right) \frac{\partial r}{\partial \TT}}{ 32m^2 \sqrt{ \xi_m (1+f) (1 + \XX^2 - \TT^2)}}\right)  \\
   & \qquad - \frac{\partial_{\TT}f (1+ \XX^2 - \TT^2) - 2\TT(1 +f)}{2 \sqrt{\xi_m \left( (1+f) (1 + \XX^2 - \TT^2)\right)^3}} \left(  \frac{3}{(1+f) (1 + \XX^2 - \TT^2)} + \frac{ \xi_m \left( 1 + \frac{2m}{r}\right)}{16m^2}\right) \Bigg]
\end{align*}
Now, we notice that \eqref{eq: Kruskal_partial_r_TX} implies 
$$\frac{\partial r}{\partial \XX} + \frac{\TT \XX}{1 + \XX^2} \frac{\partial r}{ \partial \TT} = \frac{\xi_m \, \XX \left( 1 + \XX^2 - \TT^2 \right) }{4m \, \left( 1 + \XX^2\right)}  \,,$$
while  \eqref{derivative_f_TX} gives  
\begin{align*}
    \partial_{\XX} f + \frac{\TT \XX}{ 1 + \XX^2} \partial_{\TT} f & = 2f \XX \left( -\frac{2 f (1 + \XX^2)}{\TT^2} + 2f + \frac{1}{1 + \XX^2} \right) = - \frac{2 f \XX}{1 + \XX^2} \,.
\end{align*}
Combining these identities with the two above and \eqref{eq: Kruskal_derivative_lapse}, we obtain
\begin{align*}
     \Bigg\vert &  \widetilde{V}_{\XX} \left(\G_m\left( \nabla_{\widetilde{V}_{\XX}} \widetilde{\mathcal{T}}, \widetilde{V}_{\XX}\right) \right) \Bigg\vert \\
     & = \sqrt{\frac{1 + \XX^2}{\xi_m (1+f) (1 + \XX^2 - \TT^2)}} \left[  \partial_{\XX} \Bigg(\mathfrak{g}_m\left( \overline{\nabla}_{\widetilde{V}_{\XX}} \widetilde{\mathcal{T}}, \widetilde{V}_{\XX}\right) \Bigg) + \frac{\TT \XX}{1 + \XX^2} \partial_{\TT} \Bigg(\mathfrak{g}_m\left( \overline{\nabla}_{\widetilde{V}_{\XX}} \widetilde{\mathcal{T}}, \widetilde{V}_{\XX}\right) \Bigg)\right] \\
     &= \frac{|\TT|| \XX|}{16m^2 (1 +f) (1 + \XX^2)}  \left( \frac{1 + \frac{2m}{r}}{(1+f)(1 + \XX^2 - \TT^2)} + \xi_m \frac{ \left( 1 + \frac{2m}{r} \right)^2  + \frac{8m^2}{r^2} }{16m^2} \right) \\
     & \quad + \frac{|\TT|| \XX|}{\xi_m (1 + \XX^2- \TT^2) \sqrt{(1 + f)^3}(1+\XX^2)} \left(  \frac{3}{(1+f) (1 + \XX^2 - \TT^2)} + \frac{ \xi_m \left( 1 + \frac{2m}{r}\right)}{16m^2}\right) \leq \frac{C}{m^2}\,.
\end{align*}
Combining this last estimate with the observations on the identity \eqref{observations Hessian tilde t}, we get
\begin{equation*}
    \left\vert \mathfrak{g}_m \left( \nabla^2_{\widetilde{V}_{\XX} \widetilde{V}_{\XX}} \widetilde{\mathcal{T}}, \widetilde{V}_{\XX} \right) \right\vert + \left\vert \mathfrak{g}_m \left( \nabla^2_{\widetilde{V}_{\XX} \widetilde{V}_{\XX}} \widetilde{\mathcal{T}}, \widetilde{\mathcal{T}}\right) \right\vert + \left\vert \mathfrak{g}_m \left( \nabla^2_{\widetilde{V}_{\XX} \widetilde{V}_{\XX}} \widetilde{\mathcal{T}}, \widetilde{V}_{I} \right) \right\vert \leq \frac{C}{m^2}\,.
\end{equation*}

Analogous computations give estimates on \small $\nabla^2_{\widetilde{V}_{\XX} \widetilde{\mathcal{T}}} \widetilde{\mathcal{T}}$, $\nabla^2_{\widetilde{V}_{I} \widetilde{V}_{\XX}} \widetilde{\mathcal{T}}$, $\nabla^2_{\widetilde{V}_{I} \widetilde{\mathcal{T}}} \widetilde{\mathcal{T}}$, $\nabla^2_{\widetilde{V}_{I} \widetilde{V}_J} \widetilde{\mathcal{T}}$   \normalsize and \small $\nabla^2_{\widetilde{\mathcal{T}} \widetilde{\mathcal{T}}} \widetilde{\mathcal{T}}$ \normalsize
and hence we conclude that there exists a constant $C = C(\kappa_0) > 0$ such that in the region $\{ r > \kappa_0m\}$ it holds 
\begin{equation}\label{A_Bound_Nabla2_T}
     {}^{\widetilde{\TT}} \| \nabla^2 \widetilde{\mathcal{T}}\| \leq  \frac{C}{m^2}\,.
 \end{equation}

\end{appendices}

\vspace{1cm}

\addtocontents{toc}{\protect\setcounter{tocdepth}{2}}

\printbibliography[heading=bibintoc,title={References}]

\vspace{1cm}

\end{document}